\documentclass[11pt]{article}
\renewcommand{\title}[1]{\begin{center}\textbf{\large #1}\end{center}}
\renewcommand{\author}[1]{\begin{center}#1\end{center}}
\renewcommand{\date}[1]{\begin{center}#1\end{center}}
\usepackage{amssymb, amsmath, euscript, amsthm, verbatim, array, url, mathrsfs, amscd, sw, default, color, graphicx, mathtools, slashed}
\usepackage[all]{xy}

\pdfoutput=1
\begin{document}


\title{The Seiberg-Witten Equations on Manifolds with Boundary II: \\ Lagrangian Boundary Conditions for a Floer Theory}
\author{Timothy Nguyen\footnote{The author was supported by NSF grants DMS-0706967 and DMS-0805841.}}

\abstract{In this paper, we study the Seiberg-Witten equations on the product $\R \times Y$, where $Y$ is a compact $3$-manifold with boundary. Following the approach of Salamon and Wehrheim in \cite{We1} and \cite{SaWe} in the instanton case, we impose Lagrangian boundary conditions for the Seiberg-Witten equations. The resulting equations we obtain constitute a nonlinear, nonlocal boundary value problem. We establish regularity, compactness, and Fredholm properties for the Seiberg-Witten equations supplied with Lagrangian boundary conditions arising from the monopole spaces studied in \cite{N1}.  This work therefore serves as an analytic foundation for the construction of a monopole Floer theory for $3$-manifolds with boundary.}

\tableofcontents

\section{Introduction}

Consider the Seiberg-Witten equations on a $4$-manifold $X$.  These equations are a system of nonlinear partial differential equations for a connection and spinor on $X$.  When $X$ is a product $\R \times Y$, where $Y$ is a closed $3$-manifold, the Seiberg-Witten equations on $\R \times Y$ become the formal downward gradient flow of the Chern-Simons-Dirac functional on $Y$. The associated Floer theory of the Chern-Simons-Dirac functional has been extensively studied, and after setting up the appropriate structures, we obtain the monopole Floer homology groups of $Y$, which are interesting topological invariants of $Y$ \cite{KM}.

In this paper, we consider the case when $Y$ is a $3$-manifold with boundary.  To obtain well-posed equations on $\R \times Y$ in this case, we must impose boundary conditions for the Seiberg-Witten equations.  Following the approach of \cite{We1} and \cite{SaWe}, we impose Lagrangian boundary conditions, which means that at every time $t \in \R$, a solution of our equations must have its boundary value lying in a fixed Lagrangian submanifold $\fL$ of the boundary configuration space.  The resulting equations become a Floer type equation on the space of configurations whose boundary values lie in $\fL$.  Understanding the analytic underpinnings of the Seiberg-Witten equations on $\R\times Y$ with Lagrangian boundary conditions is therefore a first step in defining a monopole Floer theory for the pair $(Y,\fL)$ of a $3$-manifold $Y$ with boundary and a Lagrangian $\fL$.

This paper is the analogue of \cite{We1} for the Seiberg-Witten setting, since \cite{We1} establishes similar foundational analytical results for the anti-self-dual (ASD) equations with Lagrangian boundary conditions.  The analysis there was eventually used to construct an instantion Floer homology with Lagrangian boundary conditions in \cite{SaWe}.  This latter work was the first to construct a gauge-theoretic Floer theory using Lagrangian boundary conditions.  The original motivation for \cite{SaWe} was to prove the Atiyah-Floer conjecture in the ASD setting (see also \cite{WeAF}).  Informally, this conjecture states the following: given a homology $3$-sphere $Y$ with a Heegard splitting $H_0 \cup_\Sigma H_1$, where $H_0$ and $H_1$ are two handlebodies joined along the surface $\Sigma$, there should be a natural isomorphism between the instanton Floer homology for $Y$ and the symplectic Floer homology for the pair of Lagrangians $(\L_{H_0}, \L_{H_1})$ inside the representation variety of $\Sigma$.  Here, $\L_{H_i}$ is the moduli space of flat connections on $\Sigma$ that extend to $H_i$, $i=0,1$.  As explained in \cite{WeAF}, instanton Floer homology with Lagrangian boundary conditions is expected to serve as an intermediary Floer homology theory in proving the Atiyah-Floer conjecture, where the instanton Floer homology of the pair $([0,1]\times\Sigma, \L_{H_0}\times \L_{H_1})$ should interpolate between the two previous Floer theories.

For the Seiberg-Witten setting, one could also formulate an analogous Atiyah-Floer type conjecture, although in this case, one ends up with infinite-dimensional Lagrangians inside an infinite-dimensional symplectic quotient (we will discuss this more thoroughly later).  One could also expect (as in the instanton case) that a monopole Floer homology for a $3$-manifold $Y$ with boundary $\Sigma$ supplied with suitable Lagrangian boundary conditions  should recover the usual monopole Floer homology for closed extensions $\bar Y = Y \cup_{\Sigma} Y'$ of $Y$, where the bounding $3$-manifolds $Y'$ satisfy appropriate hypotheses.  These considerations served as our preliminary motivation for laying the foundational analysis for a monopole Floer homology with Lagrangian boundary conditions.  Recently, several other Floer theories on $3$-manifolds with boundary have been constructed, in particular, the bordered Heegaard Floer homology theory of Lipshitz-Ozsv\'ath-Thurston \cite{LOT} and the sutured monopole Floer homology theory of Kronheimer-Mrowka \cite{KM2}.  A complete construction of a monopole Floer theory with Lagrangian boundary conditions would therefore add to this growing list of Floer theories, and it would be of interest to understand what relationships, if any, exist among all these theories.\\

\noindent\textit{Basic Setup and Main Results\\}

In order to give precise meaning to the notion of a Lagrangian boundary condition for the Seiberg-Witten equations, we first explain the infinite-dimensional symplectic aspects of our problem.  We will then explain the geometric
significance of our setup and its applicability to Floer homology after a precise statement of our main results.  Recall from \cite{N1} that the boundary configuration space
$$\fC(\Sigma) = \A(\Sigma) \times \Gamma(\S_\Sigma)$$
of connections and spinors on the boundary comes equipped with the symplectic form
\begin{equation}
  \omega((a,\phi), (b,\psi)) = \int_\Sigma a \wedge b + \int_\Sigma \Re (\phi,\rho(\nu)\psi),  \qquad (a,\phi), (b,\psi) \in \Omega^1(\Sigma; i\R) \oplus \Gamma(\S_\Sigma) \label{P3:omega}
\end{equation}
on each of its tangent spaces. Here, $\S_\Sigma$ is the spinor bundle on $\R \times Y$ restricted to $\Sigma$ and $\rho(\nu)$ is Clifford multiplication by the outward normal $\nu$ to $\Sigma$.  The form $\omega$ is symplectic because it has a compatible complex structure
\begin{equation}
  J_\Sigma = (-\check{*}, -\rho(\nu)), \label{P3:JSigma}
\end{equation}
that is, $\omega(\cdot,J_\Sigma\cdot)$ is the $L^2$ inner product on
$\T_\Sigma := \Omega^1(\Sigma; i\R) \oplus \Gamma(\S_\Sigma)$ naturally induced from the Riemannin metric on differential forms and the real part of the Hermitian inner product on the space of spinors.  Here, $\check{*}$ denotes the Hodge star operator on $\Sigma$. It follows that the $L^2$ closure of the configuration space $L^2\fC(\Sigma)$ is a Hilbert manifold whose tangent spaces $L^2\T_\Sigma$ are all strongly symplectic Hilbert spaces \cite[Section A.2]{N1}.

It follows that $\omega$ induces a nondegenerate skew-symmetric form on any topological vector space $\mathcal{X}$ densely contained in $L^2\T_\Sigma$.  As in \cite{N1}, we still call $\omega|_{\X}$ a symplectic form on $\X$. We need to consider the restriction of $\omega$ to other topological vector spaces not only because we work with the smooth configuration space.  Since we will be considering Sobolev spaces on $\R \times Y$, we also need to complete the boundary configuration space in Besov spaces, these latter spaces being boundary value spaces of Sobolev spaces. More precisely, given a smooth manifold $M$ (possibly with boundary), we can consider the Sobolev spaces $H^{s,p}(M)$ and Besov spaces $B^{s,p}(M)$ on $M$, where $s \in \R$ and $1 < p < \infty$.  When $s$ is a nonnegative integer, $H^{s,p}(M)$ is just the usual space of functions that have all derivatives up to order $s$ belonging to $L^p(M)$.  The Besov spaces are defined as in \cite{N1}, and their most important feature is that for any $s > k/p$, if $N \subset M$ is a codimension $k$ submanifold of $M$, there is a continuous restriction map
\begin{align}
  r_N: H^{s,p}(M) & \to B^{s-k/p,p}(N) \nonumber \\
  f & \mapsto f|_N. \label{eq3:trace}
\end{align}
Thus, restriction to a submanifold maps a Sobolev space on $M$ into a Besov space on $N$ and decreases the order of regularity by $k/p$. (Note that for $p = 2$, Besov spaces coincide with Sobolev spaces and the above result becomes the familiar fact that a codimension $k$ restriction decreases regularity by $k/2$ fractional derivatives).

For $s > 0$ and $p \geq 2$, we have an inclusion $B^{s,p}(\Sigma) \subset L^2(\Sigma)$.  Thus, consider  $B^{s,p}\fC(\Sigma)$, the closure the smooth configuration space $\fC(\Sigma)$ in the $B^{s,p}(\Sigma)$ topology.  The symplectic form (\ref{P3:omega}) induces a (weak) symplectic form on the Banach configuration space $B^{s,p}\fC(\Sigma)$ and the smooth configuration space $\fC(\Sigma)$.  Since $\omega$ possesses a compatible complex structure $J_\Sigma$, we can define a Lagrangian subspace of $\T_\Sigma$ to be a closed subspace $L$ such that $\T_\Sigma = L \oplus J_\Sigma L$ as a direct sum of topological vector spaces (see \cite[Section A.2]{N1} for further reading). A Lagrangian submanifold of $\fC(\Sigma)$ is then a (Fr\'echet) submanifold of $\fC(\Sigma)$ for which each tangent space is a Lagrangian subspace of $\T_\Sigma$. A Lagrangian subspace (submanifold) of the $B^{s,p}(\Sigma)$ completion of $\T_\Sigma$ is defined similarly.

\begin{Def}
  Fix $p > 2$.  An $H^{1,p}$ \textit{Lagrangian boundary condition} is a choice of a closed Lagrangian submanifold $\fL$ of $\fC(\Sigma)$ whose closure $\fL^{1-2/p,p}$ in the $B^{1-2/p,p}(\Sigma)$ topology is a smoothly embedded Lagrangian submanifold of $B^{1-2/p,p}\fC(\Sigma)$.
\end{Def}

Here, the modifier $H^{1,p}$, which we may omit in the future for brevity, expresses the fact that we will be considering our Seiberg-Witten equations on $\R \times Y$ in the $H^{1,p}_\loc(\R \times Y)$ topology.  Here, $H^{s,p}_\loc(\R \times Y)$, denotes the space of functions whose restriction to any compact subset $K \subset \R \times Y$ belongs to $H^{s,p}(K)$, $s \in \R$.

The significance of a Lagrangian boundary condition is that we can impose the following boundary conditions for a $\spinc$ connection $A$ \label{p:Aconn} and spinor $\Phi$\label{p:Phi} on $\R \times Y$ of regularity $H^{1,p}_\loc(\R \times Y)$.  Namely, we require
\begin{equation}
  (A,\Phi)|_{\{t\} \times \Sigma} \in \fL^{1-2/p,p}, \qquad \forall t \in \R, \label{eq3:intro-LBC}
\end{equation}
i.e., the restriction of $(A,\Phi)$ to every time-slice $\{t\} \times \Sigma$ of the boundary lies in the Lagrangian submanifold $\fL^{1-2/p,p}$.  Here, we made use of (\ref{eq3:trace}) with $k=2$.  Note that this restriction theorem requires $p > 2$ when $s = 1$, thereby requiring $p > 2$ Sobolev spaces for the analysis on $\R \times Y$ and the subsequent use of Besov space on the boundary. We require our boundary condition to be given by a Lagrangian submanifold because it allows us to give a Morse-Novikov-Floer theoretic interpretation of the Seiberg-Witten equations supplied with the boundary conditions (\ref{eq3:intro-LBC}), a viewpoint which we discuss more thoroughly after stating our main results\footnote{The Lagrangian property is also crucial for the analytic details of the proofs of our main results, see the outline at the end of this introduction.}.

From the results of \cite{N1}, we have a natural class of Lagrangian boundary conditions.  Namely, consider a $3$-manifold $Y'$ with $\partial Y' = -\Sigma$ and such that $Y' \cup_\Sigma Y$ is a smooth Riemannian $3$-manifold.  Moreover, suppose the $\spinc$ structure $\s$ on $Y$ extends smoothly to a $\spinc$ structure $\s'$ over $Y'$.  In such a case, the
boundary configuration spaces arising from $Y$ and $Y'$ can be identified, and so can their Lagrangian submanifolds, since the symplectic forms induced on $\fC(\Sigma)$ from $Y'$ and $Y$ differ by a minus sign\footnote{Since this is the essential property for $Y'$, in actuality, one merely need that the Riemannian metric on $Y' \cup_\Sigma Y$ be continuous instead of smooth.} (the induced orientation on $\Sigma$ differ in the two cases).

Consequently, the main theorem of \cite{N1} provides us with a Lagrangian boundary condition.  Namely, define
\begin{equation}
  \L(Y',\s') := \{(B',\Psi')|_\Sigma : (B',\Psi') \in \fC(Y'),\; SW_3(B',\Psi') = 0\} \subset \fC(\Sigma)
\end{equation}
to be the space of boundary values of connections and spinors $(B',\Psi')$ belonging to the configuration space $\fC(Y')$ on $Y'$ that solve the monopole equations $SW_3(B',\Psi') = 0$ on $Y'$.  Then if
\begin{equation}
  c_1(\s') \textrm{ is non-torsion or } H^1(Y',\Sigma) = 0, \label{P3:assum}
\end{equation}
the main result of \cite{N1} is that $\L(Y',\s')$ is an $H^{1,p}$ Lagrangian boundary condition for $p > 4$.

\begin{Def}\label{DefML}
  Let $Y'$ and $\s'$ be as above, and suppose they satisfy (\ref{P3:assum}).  Then we call the Lagrangian submanifold $\L(Y',\s') \subset \fC(\Sigma)$ a \textit{monopole Lagrangian}.
\end{Def}

Our main result is that the Seiberg-Witten equations $SW_4(A,\Phi) = 0$, defined by (\ref{SW4}), supplied with a Lagrangian boundary condition arising from a monopole Lagrangian yields an elliptic boundary value problem, i.e., one for which elliptic regularity modulo gauge holds.  Here the gauge group is $\G = \Maps(\R \times Y, S^1)$ and it acts on the configuration space $\fC(\R \times Y)$ of $\spinc$ connections and spinors on $\R \times Y$ (where the $\spinc$ structure on $\R \times Y$ has been fixed and pulled back from a $\spinc$ structure on $Y$) via
$$(A,\Phi) \mapsto g^*(A,\Phi) = (A - g^{-1}dg, g\Phi).$$
Let $\G_\id$ denote the identity component of the gauge group, and let the prefix $H^{s,p}_\loc$ denote closure with respect to the $H^{s,p}_\loc(\R \times Y)$ topology.\\

\noindent \textbf{Theorem A (Regularity). }\textit{Let $p > 4$, and let $(A,\Phi) \in H^{1,p}_\loc\fC(\R \times Y)$ solve the boundary value problem
\begin{equation}
\begin{split}
  SW_4(A,\Phi) & =0 \\
  (A,\Phi)|_{\{t\}\times\Sigma} & \in \fL^{1-2/p,p}, \quad \forall t \in \R,
  \end{split}\label{BVP}
\end{equation}
where $\fL^{1-2/p,p}$ denotes the $B^{1-2/p,p}(\Sigma)$ closure of a monopole Lagrangian $\fL$.  Then there exists a gauge transformation $g \in H^{2,p}_\loc\G_\id$ such that $g^*(A,\Phi)$ is smooth.\\}

Next, we have a compactness result for sequences of solutions provided that the Lagrangian $\fL$ is invariant under the gauge group action of $\G(Y)|_{\Sigma}$.  If $\fL$ satisfies this, we say that $\fL$ is \textit{fully gauge-invariant}.  Observe that if we take $Y' = -Y$ and $\s'$ the $\spinc$ structure on $Y$, then $\L(Y',\s')$ will be a fully gauge-invariant monopole Lagrangian, provided $Y'$ and $\s'$ satisfy (\ref{P3:assum}).  In general, the condition that $\L(Y',\s')$ be fully gauge-invariant is precisely the condition that the natural restriction maps on cohomology $H^1(Y) \to H^1(\Sigma)$ and $H^1(Y') \to H^1(\Sigma)$ have equal images in $H^1(\Sigma)$. In this situation, the following theorem says that if we have a local bound on the $H^{1,p}$ ``energy" of a sequence of configurations $(A_i,\Phi_i)$, then modulo gauge, a subsequence converges smoothly on every compact subset of $\R \times Y$.  Here, the energy of a configuration on a compact set $K$ is given by the gauge-invariant norms appearing in (\ref{energybound}), where $\nabla_{A}$ denotes the $\spinc$ covariant derivative determined by the connection $A$.\\

\noindent \textbf{Theorem B (Compactness). }\textit{Let $p > 4$ and let $(A_i,\Phi_i) \in H^{1,p}_\loc\fC(\R \times Y)$ be a sequence of solutions to (\ref{BVP}), where $\fL$ is a fully gauge-invariant monopole Lagrangian. Suppose that on every compact subset $K \subset \R \times Y$, we have
\begin{equation}
  \sup_i \|F_{A_i}\|_{L^p(K)}, \|\nabla_{A_i}\Phi_i\|_{L^p(K)}, \|\Phi_i\|_{L^p(K)} < \infty. \label{energybound}
\end{equation}
Then there exists a subsequence of configurations, again denoted by $(A_i,\Phi_i)$, and a sequence of gauge transformations $g_i \in H^{2,p}_\loc\G$ such that $g_i^*(A_i,\Phi_i)$ converges in $C^\infty(K)$ for every compact subset $K \subset \R \times Y$.\\}

Both Theorems A and B apply verbatim to the periodic setting, where $\R \times Y$ is replaced with $S^1 \times Y$, in which case, we can work with $H^{k,p}(S^1 \times Y)$ spaces instead of $H^{k,p}_{\loc}(\R \times Y)$.  In fact we will prove Theorems A and B in the periodic setting, which then implies the result on $\R \times Y$ by standard patching arguments.  In the periodic setting we will also prove that the linearization of (\ref{BVP}) is Fredholm in a suitable gauge (and in suitable topologies), see Theorem \ref{ThmFredProp}. One can also prove the Fredholm property on $\R \times Y$ assuming suitable decay hypotheses at the ends, but we will not pursue that here.

Note that the requirement $p > 4$ in the above is sharp with respect to the results in \cite{N1}, in that for no value of $p \leq 4$ is it known that $\fL^{1-2/p,p}$ is a smooth Banach manifold.  Thus, our results here cannot be sharpened unless the results in \cite{N1} are also sharpened.  On the other hand, the value $p > 4$ is a priori unsatisfactory from the point of view of Floer theory.  This is because the a priori energy bounds we have on solutions to the Seiberg-Witten equations, namely the analytic and topological energy as defined in \cite{KM}, are essentially an $H^{1,2}$ control.  Therefore, Theorem B is not sufficient to guarantee  compactness results for the moduli space of solutions to (\ref{BVP}) that are of the type needed for a Floer theory. However, this is not the end the story, as can be seen in the ASD situation, where a Floer theory still exists even though the analogous regularity and compactness results are proven only for $p > 2$, which although better than $p > 4$, still misses $p = 2$.\footnote{In the ASD situation, it is possible to prove that $L^p\fL$, the $L^p$ closure of $\fL$, is a smooth submanifold of the space of $L^p$ connections, $p > 2$.  Since we have the embedding $B^{1-2/p,p}(\Sigma) \hookrightarrow L^p(\Sigma)$, the ASD equations are well-behaved for $p > 2$, and thus the analogue of Theorem A with $p > 2$ is the optimal result there.  Note however, in both the ASD and Seiberg-Witten setting, $p = 2$ can never be achieved, since then a Lagrangian boundary condition cannot be defined.  Indeed, a function belonging to $H^{1,2}$ does not have a well-defined restriction to a codimension two submanifold.}  The ASD Floer theory is possible due to the bubbling analysis carried out in \cite{We2} and the presence of energy-index formulas in \cite{SaWe}, which allow one to use the $p > 2$ analysis to understand the compactification of the space Floer trajectories between critical points. We will leave the study of the analog of such issues in the Seiberg-Witten setting for the future, namely, the study of what can happen to a sequence of solutions to (\ref{BVP}) if one is only given an $H^{1,2}$ type energy bound (more precisely, a bound on the analytic and topological energy of \cite{KM}).  At present then, our main theorems therefore serve the foundational purpose of showing that the Seiberg-Witten equations with Lagrangian boundary conditions are well-posed and satisfy a weak type of compactness. These results are key for a future construction of an associated Floer theory.\\

\noindent\textit{Geometric Origins}\\

Having stated our main results, we explain how the Seiberg-Witten equations supplied with Lagrangian boundary conditions naturally arise in trying to construct a Floer homology on a $3$-manifold with boundary.  On a product $\R \times Y$, for $Y$ with or without boundary, the Seiberg-Witten equations take the following form.

We have a decomposition
\begin{equation}
  \A(\R \times Y) = \Maps(\R,\A(Y)) \times \Maps(\R, \Omega^0(Y;i\R))\\
\end{equation}
whereby a connection $A \in \A(\R \times Y)$ on $\R \times Y$ can be decomposed as
\begin{equation}
    A = B(t) + \alpha(t)dt, \label{A}
\end{equation}
where $B(t) \in \A(Y)$ is a path of connections on $Y$ and $\alpha(t) \in \Omega^0(Y; i\R)$ is a path of $0$-forms on $Y$, $t \in \R$.  Likewise, if we write $\S^+$ for the bundle of self-dual spinors on $\R \times Y$ and $\S$ for the spinor bundle on $Y$ obtained by restriction of $\S^+$ to $\{0\} \times Y$, we can write
\begin{equation}
  \Gamma(\S^+) = \Maps(\R,\Gamma(\S)),
\end{equation}
where we have identified $\S^+$ with the pullback of $\S$ under the natural projection of $\R \times Y$ onto $Y$.  Thus, any spinor $\Phi$ on $\R \times Y$ is given by a path $\Phi(t)$ of spinors on $Y$. Altogether, the configuration space
$$\fC(\R \times Y) = \A(\R \times Y) \times \Gamma(\S^+)$$
on $\R \times Y$ can be expressed as a configuration space of paths:
\begin{align}
  \fC(\R \times Y) &= \Maps(\R, \A(Y) \times \Gamma(\S) \times \Omega^0(Y;i\R))\nonumber \\
  (A,\Phi) & \mapsto (B(t),\Phi(t),\alpha(t)). \label{paths}
\end{align}
Under this correspondence, the Seiberg-Witten equations can be written as follows \cite[Chapter 4]{KM}: \label{p:SW4}
\begin{align}
SW_4(A,\Phi) &:= \left(\frac{d}{dt}B + \left(\frac{1}{2}*_Y F_{B^t} + \rho^{-1}(\Phi\Phi^*)_0\right) - d\alpha, \frac{d}{dt}\Phi + D_B\Psi + \alpha\Phi\right) \label{eq3:introSW4} \\
&= 0,\label{SW4}
\end{align}
where we have suppressed the time-dependence from the notation.  Here, the terms appearing on the right-hand-side of (\ref{eq3:introSW4}) are defined as in \cite{N1}.  In particular, when $\alpha \equiv 0$, i.e. when $A$ is in temporal gauge, then the Seiberg-Witten equations (\ref{SW4}) are equivalent to the equations
\begin{equation}
  \begin{split}
    \frac{d}{dt}B & = -\left(\frac{1}{2}*_Y F_{B^t} + \rho^{-1}(\Phi\Phi^*)_0\right)\\
    \frac{d}{dt}\Phi &= -D_{B}\Phi.
  \end{split}\label{CSDflow}
\end{equation}
For $(B,\Psi) \in \A(Y)\times\Gamma(\S)$, we have the three-dimensional Seiberg-Witten map
\begin{equation}
  SW_3(B,\Psi) := \left(\frac{1}{2}*_Y F_{B^t} + \rho^{-1}(\Psi\Psi^*)_0, D_{B}\Psi\right) \in \Omega^1(Y; i\R) \times \Gamma(\S), \label{SW3}
\end{equation}
which we may think of as a vector field on
$$\fC(Y) = \A(Y) \times \Gamma(\S),$$
the configuration space on $Y$, since $T_{\c}\fC(Y) = \Omega^1(Y; i\R) \times \Gamma(\S)$. Thus, (\ref{CSDflow}) is formally the downward ``flow" of $SW_3$.  Moreover, the equations (\ref{SW4}) and (\ref{CSDflow}) are equivalent modulo gauge transformations.

When $Y$ is closed, the vector field $SW_3$ is the gradient of a functional, the \textit{Chern-Simons-Dirac functional} $CSD: \fC(Y) \to \R$.  This functional is defined by
\begin{align*}
  CSD(B,\Psi) = -\frac{1}{2}\int_Y (B - B_0) \wedge (F_{B} + F_{B_0}) + \frac{1}{2}\int_Y (D_B\Psi,\Psi),
\end{align*}
where $B_0$ is any fixed reference connection.  Thus, since $SW_3(B,\Psi)$ is the gradient of $CSD$ at $(B,\Psi)$, (\ref{CSDflow}) is formally the downward gradient flow of $CSD$.  In this way, the monopole Floer theory of $Y$ can be thought of informally as the Morse theory of the Chern-Simons-Dirac functional.  However, we ultimately must work modulo gauge, and in doing so, the Chern-Simons-Dirac functional does not descend to a well-defined function on the quotient configuration space, but it instead becomes an $S^1$ valued function.  This is because the first term of $CSD$ is a Chern-Simons term which is not fully gauge-invariant but which changes by an amount depending on the homotopy class of the gauge transformation.  In this case, when we take the differential of $CSD$ as a circle-valued function on the quotient configuration space\footnote{We ignore the singular points of the quotient configuration space in this informal discussion.}, we obtain not an exact form but a closed form.  On compact manifolds, the Morse theory for a closed form is more generally known as Morse-Novikov theory, from which the case of exact forms reduces to the usual Morse theory.  Thus, the monopole Floer theory of $Y$ is more accurately the Morse-Novikov theory for the differential of $CSD$ on the quotient configuration space.

It is this point of view which we wish to adopt in trying to generalize monopole Floer theory to manifolds with boundary.
Suppose $\partial Y = \Sigma$ is nonempty.  We now impose boundary conditions for (\ref{CSDflow}) that preserve the above Morse-Novikov viewpoint for monopole Floer theory.  From this, we are naturally led to Lagrangian boundary conditions as we now explain.  Define the following Chern-Simons-Dirac one-form $\mu = \mu_{CSD}$ on $T^*\fC(Y)$: \label{p:muCSD}
\begin{equation}
  \mu(b,\psi) = ((b,\psi),SW_3(B,\Psi))_{L^2(Y)}, \qquad (b,\psi) \in T_\c\fC(Y) = \Omega^1(Y;i\R) \oplus \Gamma(\S). \label{CSDform}
\end{equation}
Here, the above $L^2$ inner product on $\Omega^1(Y;i\R) \oplus \Gamma(\S)$ is the one induced from the Riemannian inner product on $Y$ and the real part of the Hermitian inner product on $\Gamma(\S)$. If we take the differential of $\mu$, we obtain
\begin{align}
  d_{(B,\Psi)}\mu((a,\phi),(b,\psi)) =  ((a,\phi),\H_\c(b,\psi))_{L^2(Y)} - (\H_\c(a,\phi),(b,\psi))_{L^2(Y)} \label{dmu}
\end{align}
where $\H_\c: \Omega^1(Y;i\R) \oplus \Gamma(\S) \to \Omega^1(Y;i\R) \oplus \Gamma(\S)$ is the ``Hessian" operator obtained by differentiating the map $SW_3$ (of course, $\H_\c$ is only a true Hessian when $Y$ is closed, because only then is $SW_3$ the gradient of a functional). Recall from \cite{N1} that the Hessian is given by
\begin{equation}
  \H_{(B,\Psi)} = \begin{pmatrix}
    *_Yd & 2i\Im\rho^{-1}(\cdot\Phi^*)_0\\
    \rho(\cdot)\Phi & D_{B}
  \end{pmatrix},
\end{equation}
a first order formally self-adjoint operator.  Thus, (\ref{dmu}) automatically vanishes on a closed manifold.  However, integration by parts shows that when $\partial Y = \Sigma$ is nonempty, (\ref{dmu}) defines a skew-symmetric pairing on the boundary.  A simple computation shows that this pairing is the symplectic form $\omega$ in (\ref{P3:omega}).

Define the (tangential) restriction map
\begin{align}
  r_\Sigma: \fC(Y) &\to \fC(\Sigma) \nonumber \\
  (B,\Psi) & \mapsto (B|_\Sigma,\Psi|_\Sigma) \label{P3:rSigma}
\end{align}
The above discussion shows that if we pick a submanifold $\frak{X} \subset \fC(Y)$ such that for every $(B,\Psi) \in \frak{X}$, the space $r_\Sigma(T_{(B,\Psi)}\frak{X})$ is an isotropic subspace of $T_{r_\Sigma(B,\Psi)}\fC(\Sigma)$, then $d\mu|_{\frak{X}}$ vanishes.  In particular, pick a Lagrangian submanifold $\fL \subset \fC(\Sigma)$ and define the space
\begin{equation}
  \fC(Y,\fL) = \{(B,\Psi) \in \fC(Y) : r_\Sigma(B,\Psi) \in \fL\}
\end{equation}
consisting of those configurations on $Y$ whose restriction to $\Sigma$ lies in $\fL$.  Then by the above considerations, $\mu$ is a closed $1$-form when restricted to $\fC(Y,\fL)$.

It is now possible to consider the Morse-Novikov theory for $\mu$ on $\fC(Y,\fL)$.  The resulting Floer equations, i.e., the formal downward flow of $\mu$ viewed as a vector field, are the equations (\ref{CSDflow}) supplemented with the boundary condition
\begin{equation}
  r_\Sigma(B(t),\Psi(t)) \in \fL, \qquad t \in \R. \label{LBC}
\end{equation}
Here, we choose a Lagrangian submanifold (rather than an isotropic one) because Lagrangian boundary conditions are precisely the ones that give rise to self-adjoint boundary conditions\footnote{Here in this informal discussion, we are being very loose with the precise functional analytic details, since self-adjointness requires that we find Lagrangians in the correct Hilbert space of boundary value data.  See Section \ref{SecP3Linear} and Section 22 of \cite{N0} for a more rigorous discussion of the relationship between Lagrangians and self-adjointness.}.  In other words, by supplying Lagrangian boundary conditions, the linearization of the system of equations (\ref{CSDflow}) and (\ref{LBC}) yields a time-dependent family of self-adjoint operators, from which it is then possible to compute the spectral flow of such a family (provided the requisite decay properties hold at infinity).  The existence of this spectral flow makes it possible to assign a relative grading for the chain complex generated by the critical points of our ``flow", which is an essential ingredient in the future construction of a Floer theory.  Moreover, as we shall explain below, self-adjointness also plays crucial role in obtaining auxiliary elliptic estimates needed within the proof of our main theorems.     \\

\noindent \textsl{Summary of Analytic Difficulties}\\

Let us shed some insight on the analytic difficulties that the equations (\ref{BVP}) pose, and in particular, let us compare these equations with the corresponding ASD equations studied in \cite{We1}.  In both these situations, what we essentially have is an elliptic semilinear partial differential equation, with nonlinear, nonlocal boundary conditions.  By elliptic, we mean that in a suitable gauge, the principal symbol of the equations on $\R \times Y$ are elliptic.  By nonlocal, we mean that the Lagrangian boundary condition is not given by a set of differential equations on $\Sigma$.  More precisely, a tangent space to our Lagrangian $\L$ is given not by the kernel of a differential operator but of a pseudodifferential operator (at least approximately, in the sense described in \cite{N1}).  However, let us note that what is truly nonstandard about both these boundary problems is that the nonlocal boundary conditions are imposed ``slicewise", that is, they are specified pointwise in the time variable $t \in \R$. This implies that the linearization of the boundary condition (\ref{LBC}) is determined (again, approximately) by the range of a product-type pseudodifferential operator, or more precisely, a time-dependent pseudodifferential operator on $\Sigma$, which is therefore, not pseudodifferential as an operator on $\R \times \Sigma$.  We therefore have neither a local nor a nonlocal (pseudodifferential) boundary condition in the usual sense.

However, let us point out that in many ways, the ASD situation is ``almost local" whereas in our situation, this is not at all the case.  In both the ASD and Seiberg-Witten case, the action of the gauge group gives the ``local part" of the Lagrangian (a gauge orbit and its tangent space are defined by differential equations) and dividing by the gauge group gives the remaining ``nonlocal part" of the Lagrangian.   However, in the ASD case, the Lagrangians modulo gauge are finite dimensional, whereas in the Seiberg-Witten case, they are infinite dimensional.  Indeed, in the ASD case, the Lagrangians must lie in the space of flat connections, the zero set of the moment map associated to gauge group action on the space of $SU(2)$ connections on $\Sigma$.  Hence, modulo gauge, these Lagrangians descend to Lagrangian submanifolds of the (singular) finite dimensional symplectically reduced space, the representation variety of $\Sigma$. On the other hand, because of the presence of spinors in the Seiberg-Witten case, the symplectic reduction of the moment map associated to the gauge group action on $\fC(\Sigma)$ is infinite dimensional.  Hence, the Lagrangians one must consider descend to an infinite dimensional Lagrangian submanifold of this reduced space.  For monopole Lagrangians in particular, the nonlocal part of Lagrangian is of a pseudodifferential nature \cite{N1}.  Thus, our work here requires pseudodifferential analysis whereas the ASD case does not.

Moreover, it turns out that we are led to introduce some nonstandard function spaces because of the slicewise nature of the Lagrangian boundary condition.  Such a boundary condition places time ($t \in \R$) and space (the manifold $Y$) on a different footing, and consequently, the estimates we perform on $\R \times Y$ will, in particular, measure regularity in time and space differently.  Function spaces that distinguish among different directions are known as \textit{anisotropic} function spaces (in contrast to the usual isotropic function spaces that measure the regularity of a function equally in every direction.) Specifically, we are required to work with both (isotropic) Besov spaces and also \textit{anisotropic Besov spaces}.  See Section 2 for a definition of these spaces. 

Of course, our Lagrangians are not linear objects, and thus we will have to do a fair amount of nonlinear analysis in conjunction with the pseudodifferential nature of the our Lagrangians in the setting of anisotropic Besov spaces.  This is in contrast to the ASD situation, where since the Lagrangians are finite dimensional modulo gauge, and all norms on a finite dimensional space are equivalent, there are no functional analytic difficulties posed by the nonlinearities of the Lagrangian.  That is, in the ASD case, the nonlinearity of the Lagrangian only becomes a central issue in the bubbling analysis of \cite{We2}, whose importance we described after Theorem B, and not the elliptic regularity analysis.

In fact, it is the necessity of such future bubbling analysis for the Seiberg-Witten case that requires that we work with the $H^{1,p}$ spaces in Theorems A and B.  In \cite{We2}, the bubbling analysis in the instanton case requires that the instanton analogs of Theorems A and B are proven for $H^{1,p}$ for $p > 2$, and not say, for $H^{k,2}$ with $k$ large.  Indeed, to exclude bubbling in the situations relevant for defining Floer homology, \cite{We2} applies mean value inequalities to deduce that the energy density of instantons remain bounded.  Since bounds on the energy density always yield $H^{1,p}$ control (even up to the boundary) of an instanton modulo gauge, we immediately obtain a compactness result from the $H^{1,p}$ compactness theorem.  (If one only had an $H^{k,2}$ compactness result with $k > 1$, bounds on the energy density would not be sufficient to yield compactness.) Thus, it is expected that future analysis of the bubbling phenomenon for the Seiberg-Witten equations will depend on $H^{1,p}$ analysis as well.  Moreover, we should remark that proving our main theorems for $H^{k,2}$ regular solutions for large $k$ does not simplify the analysis in any fundamental or conceptual way, since all the main technical steps, which we outline below, will still need to be performed.\\

\noindent\textsl{Outline:} This paper is organized as follows. In Section \ref{SecPaths}, we define the anisotropic function spaces we will be using. We then apply these function spaces to study the space of paths through a monopole Lagrangian $\fL$.  This is necessary in light of the correspondence (\ref{paths}), which relates configurations on $\R \times Y$ to paths through the configuration space on $Y$.  The boundary condition (\ref{LBC}) specifies that on the boundary, such a path is a path through $\fL$.  Here, we make use of the analysis developed in \cite{N1} in order to show that the space of paths through monopole Lagrangians is again a manifold in (anisotropic) Besov space topologies.  The main theorems of this section, Theorems \ref{ThmPathL1} and \ref{ThmPathL2}  are where we do our main nonlinear analysis on anisotropic Besov spaces. In Section \ref{SecP3Linear}, we show that the linearization of (\ref{BVP}) is Fredholm in the periodic setting and in the appropriate function space topologies (including anisotropic ones).  (That the linearization of (\ref{BVP}) makes sense follows from the results of Section 2.) Here, a key step is to establish a resolvent estimate on anisotropic function spaces. This resolvent estimate also relies crucially on aspects pertaining to self-adjointness, which in the end, amounts to the requirement that our boundary condition in (\ref{BVP}) be given by a Lagrangian.  In Section \ref{SecP3Sec4}, we apply the tools from the previous sections to prove our main theorems, whose proof we summarize as follows.  First, we work locally in time, which means we replace the time interval $\R$ in Theorem A with $S^1$ in Theorem \ref{ThmA'}. After this, the first step is to place the equations in a suitable gauge such that the linear part of the resulting equations is elliptic and falls into the framework of Section \ref{SecP3Linear}.  From this, the second step is to gain regularity in the $\Sigma$ directions for the gauge-fixed $(A,\Phi)$ in a neighborhood of the boundary using the anisotropic estimates of Sections \ref{SecPaths} and \ref{SecP3Linear}.  The third step is to gain regularity in the time direction and normal direction to $\Sigma$ using the theory of Banach space valued Cauchy-Riemann equations due to Wehrheim \cite{We}.  Once we have gained some regularity in all directions, then in our final step, we bootstrap to gain regularity to any desired order, which proves Theorem \ref{ThmA'}.  We then deduce Theorems A and B in Section 4 from Theorem \ref{ThmA'}. (Note if one starts with $H^{k,2}$ regularity in Theorem A for $k$ large, one never has to work with $p \neq 2$ spaces and so on a first reading, one may assume that this is the case for simplicity. One still has to bootstrap in the above anisotropic fashion, however.) As a note, we have deferred the proof of certain technical steps in the above to \cite{N0} since they are rather involved.  These steps include the nonlinear interpolation procedure of Section 2 and establishing a resolvent estimate on anisotropic spaces in Section 3.

Let us finally remark that the analysis we do is of a very general character and is likely to be applicable to other elliptic, semilinear boundary value problems whose linear part is a Dirac type operator and whose slicewise boundary condition (\ref{LBC}) is given by a Lagrangian submanifold $\fL$ which satisfies formally similar properties to those obeyed by the monopole Lagrangians.  Indeed, he proof of Theorem A shows that the essential analysis we do has very little to do with the fact that we are dealing with the Seiberg-Witten equations per se.  In the general situation, for operators of Dirac type on $\R \times Y$, then near the boundary, the part of the operator that differentiates in the $\R$ and normal directions becomes a Cauchy-Riemann operator.  For this operator, we can apply the methods of Section \ref{AppCR} to gain regularity in the $\R$ and normal directions near the boundary so long as we have gained smoothness in the remaining $\Sigma$ directions.  Such regularity may be obtained by the anisotropic linear theory of Section \ref{SecPaths} (where the anisotropy is in the $\Sigma$ directions), so long as the tangent spaces to the Lagrangian $\fL$ fall within the framework of Section \ref{SecPaths} and the nonlinear components of its chart maps for $\fL$ smooth in the $\Sigma$ directions (i.e. the space of paths through $\fL$ satisfies properties similar to those in Theorems \ref{ThmPathL1} and \ref{ThmPathL2}).  It is not unreasonable to expect that any naturally occurring Lagrangians should satisfy such properties. Indeed, if they are locally determined by the zero locus of an analytic function involving multiplication and pseudodifferential--like operators (as monopole Lagrangians are), then as done in Section \ref{SecPaths}, an analysis of the nonlinear compositions of multiplication and pseudodifferential--like operators which enter into the analysis of their chart maps should allow us to recover suitable smoothing properties.  Nevertheless, the analysis we do in Section \ref{SecPaths} is quite involved, and it would have been unwise to obscure the technical exposition of Section \ref{SecPaths} by writing it for abstract Lagrangians, even though after the fact, one is led to believe that the results there should hold more generally.

Of course, constructing a monopole Floer theory with Lagrangian boundary conditions will have to depend on the precise nature of the Seiberg-Witten equations and the chosen Lagrangian. In this regard, the Lagrangians we have chosen, namely the monopole Lagrangians, are the most natural ones to consider.\\

\noindent\textbf{Notation. }All reference numbers with a prefix ``I." refer to reference numbers in the prequel paper \cite{N1}.  The reading of this paper will require familiarity with \cite{N1}.     


\section{Spaces of Paths}\label{SecPaths}

A smooth path $(B(t),\Phi(t)) \in \Maps(\R,\fC(Y))$ satisfying the boundary condition (\ref{LBC}) implies that we get a path
$$r_\Sigma(B(t),\Phi(t)) \in \Maps(\R,\fL)$$
through the smooth Lagrangian $\fL$, which we assume to be a monopole Lagrangian $\fL = \L(Y',\s')$.  Our task in this section is to show that the space of paths through $\fL$ in (anisotropic) Besov space topologies forms a Banach manifold obeying analytic properties suitable for proving Theorems A and B.  These properties are summarized in the main theorems of this section, Theorems \ref{ThmPathL1} and \ref{ThmPathL2}, and among these properties, of particular significance is the fact that the nonlinear component of a Banach manifold chart map for the space of paths through a monopole Lagrangian is smoothing, i.e. it increases regularity, in the $\Sigma$ directions\footnote{This is nothing special to the moduli spaces involved in Seiberg-Witten theory.  The analogous estimates carried out here apply mutatis mutandis to instanton moduli spaces.  This is because the Seiberg-Witten equations and instanton equations have the same shape: they are semi-linear with quadratic nonlinearity.}

In detail, recall from \cite{N1} that given a monopole Lagrangian $\L$, then $\L^{s-1/p,p}$, the $B^{s-1/p,p}(\Sigma)$ closure of $\L,$ is a smooth Banach submanifold of $\fC^{s-1/p,p}(\Sigma)$, the $B^{s-1/p,p}(\Sigma)$ closure of $\fC(\Sigma)$ for $s > \max(3/p,1/2+1/p)$.  Furthermore, the local chart maps of $\L^{s-1/p,p}$ are described by Theorem I.3.4.  In this theorem, the nonlinear part of a particular chart map at a configuration $u \in \L^{s-1/p,p}$, which we denoted by $E^1_u$, is smoothing.  However, if we now consider the space of paths through our Lagrangians, then since the corresponding chart maps on the space of paths will be defined ``slicewise" along the path (see Definition \ref{DefSw}), the smoothing continues to occur but only in the space variables $\Sigma$ and not in the time variable. This naturally leads us to consider spaces which have extra smoothness in some directions and hence anisotropic spaces.  Because our monopole spaces are modeled on Besov spaces, we thus end up with anisotropic Besov spaces.

Before proceeding to define these spaces, we fix some notation.  In this section, we consider Lagrangians $\L(Y,\s)$, which we abbreviate as $\L$, where $Y$ is any $3$-manifold and $\s$ is any $\spinc$ structure on $Y$ such that (\ref{P3:assum}) holds.  Later, we instead use $Y$ to denote the $3$-manifold occurring in the Seiberg-Witten equations on $\R \times Y$. The Lagrangians one considers for $\R \times Y$ are those of the form $\L(Y',\s')$, for some other $3$-manifold $Y'$, with $Y'$ and $\s'$ satisfying Definition \ref{DefML}. Thus, the manifold $Y$ in this section will be replaced with $Y'$ in later sections.

Recall the following notation from \cite{N1}, which we need so that we can analyze how $\L$ is constructed.  With $\s$ fixed, we have the configuration spaces
\begin{align}
  \fC(Y) &= \A(Y) \times \Gamma(\S)\\
  \fC(\Sigma) &= \A(\Sigma) \times \Gamma(\S_\Sigma),
\end{align}
of connections and spinors on $Y$ and $\Sigma$, where $\S$ is the spinor bundle on $Y$ associated to $\s$ and $\S_\Sigma$ is the restriction of $\S$ to $\Sigma$.  Both of these spaces are affine spaces modeled on
\begin{align}
  \T = \T_Y &=  \Omega^1(Y;i\R) \oplus \Gamma(\S) \label{sec2:TY}\\
  \T_\Sigma &= \Omega^1(\Sigma; i\R) \oplus \Gamma(\S_\Sigma), \label{sec2:TS}
\end{align}
respectively, and the tangential restriction map (\ref{P3:rSigma}) on configuration spaces induces one on the tangent spaces:
\begin{align}
  r_\Sigma: \T & \to \T_\Sigma \nonumber \\
  (b,\psi) & \mapsto (b|_\Sigma,\psi|_\Sigma).
\end{align}
The space of monopoles $\fM = \fM(Y,\s)$ is the zero set of the Seiberg-Witten map $SW_3$ given by (\ref{SW3}). Fixing a smooth reference connection $B_\rf$, we have the space
\begin{equation}
  \M = \{(B,\Psi) \in \fM: d^*(B-B_\rf) = 0\} \label{sec2:M}
\end{equation}
of monopoles in Coulomb gauge with respect to $B_\rf$.  The space $\L$ is the space of boundary values of $\fM$, which is equal to the space of boundary values of $\M$, i.e.,
\begin{equation}
  \L = r_\Sigma(\fM) = r_\Sigma(\M),
\end{equation}
where $r_\Sigma$ is the tangential restriction map (\ref{P3:rSigma}).

All these definitions extend to the appropriate Besov completions, as was done in \cite{N1}.  Thus, we have the configuration spaces
$\fC^{s,p}(Y)$ and $\fC^{s,p}(\Sigma)$, the $B^{s,p}(Y)$ and $B^{s,p}(\Sigma)$ closures of $\fC(Y)$ and $\fC(\Sigma)$, respectively. We also have the following Besov monopole spaces
\begin{align}
  \fM^{s,p} &= \{(B,\Psi) \in \fC^{s,p}(Y) : SW_3(B,\Psi) = 0\}, \label{sec2:def-fM}\\
  \M^{s,p} &= \{(B,\Psi) \in \fC^{s,p}(Y) : SW_3(B,\Psi) = 0, d^*(B-B_\rf)=0\}, \label{sec2:def-M}\\
  \L^{s-1/p,p} &= r_\Sigma(\M^{s-1/p,p}). \label{sec2:def-L}
\end{align}
From \cite{N1}, for $p \geq 2$ and $s > \max(3/p,1/2)$, the spaces $\fM^{s,p}$ and $\M^{s,p}$ are Banach submanifolds of $\fC^{s,p}(Y)$.  If in addition, $s > \max(3/p,1/2+1/p)$, then $\L^{s-1/p,p}$ is a Banach submanifold of $\fC^{s-1/p,p}(\Sigma)$ and
\begin{equation}
  r_\Sigma: \M^{s,p} \to \L^{s-1/p,p}
\end{equation}
is a covering.

When $\L^{s-1/p,p}$ is a Banach manifold, the space $C^0(I,\L^{s-1/p,p})$ of continuous paths from an interval $I$ into $\L^{s-1/p,p}$ is naturally a Banach manifold. However, it is far from obvious that the space of paths through $\L^{s-1/p,p}$ in anisotropic Besov topologies is a Banach manifold.  We now define these anisotropic Besov spaces precisely.

\subsection{Anisotropic Function Space Setup}

On Euclidean space, the usual Besov spaces $B^{s,p}(\R^n)$ are well-defined, for any $s \in \R$ and $1<p<\infty$, see \cite{N1} or \cite{Tr}.  Recall that for $p = 2$, the  spaces $B^{s,2}(\R^n)$ coincide with the $L^2$ Sobolev spaces $H^{s,2}(\R^n)$.  Suppose we have a splitting $\R^n = \R^{n_1} \times \R^{n_2}$.  Then for any $s_1 \in \R$ and $s_2 \geq 0$, we define the anisotropic Besov space $B^{(s_1,s_2),p}(\R^{n_1} \times \R^{n_2})$ as follows.  Let $x = (x_{(1)},x_{(2)})$ be the coordinates on $\R^{n_1} \times \R^{n_2}$ and let $\xi = (\xi_{(1)},\xi_{(2)})$ be the corresponding Fourier variables.  Let $\F_2$ denote Fourier transform on $\R^{n_2}$ and define the operator
$$J_{(2)}f = \F_2^{-1}\left(1+\left|\xi_{(2)}\right|^2\right)^{1/2}\F_2f$$
for any smooth compactly supported $f$.  Thus, very roughly speaking, for any $s \geq 0$, $J_{(2)}^s$ takes $s$ derivatives in all the $\R^{n_2}$ directions.

\begin{Def}
  For $s_1 \in \R$ and $s_2 \geq 0$, the \textit{anisotropic Besov space} $B^{(s_1,s_2),p}(\R^{n_1} \times \R^{n_2})$ is defined to be the closure of smooth compactly supported functions on $\R^n$ with respect to the norm
\begin{equation}
  \|f\|_{B^{(s_1,s_2),p}(\R^{n_1}\times\R^{n_2})} := \|J_{(2)}^{s_2}f\|_{B^{s_1,p}(\R^n)}.
\end{equation}
\end{Def}

Equivalently  (see \cite{N3}), the space $B^{(s_1,s_2),p}(\R^{n_1}\times\R^{n_2})$ is the space of all functions lying in the classical (isotropic) Besov space $B^{s_1,p}(\R^n)$ for which $D_{s_2} f \in B^{s_1,p}(\R^n)$, where $D_{s_2}$ is any smooth elliptic operator on $\R^{n_2}$ of order $s_2$.  That is, we have the equivalence of norms
\begin{equation}
  \|f\|_{B^{(s_1,s_2),p}(\R^{n_1}\times\R^{n_2})} \sim \|f\|_{B^{s_1,p}(\R^n)} + \|D_{s_2}f\|_{B^{s_1,p}(\R^n)}. \label{liftAB}
\end{equation}
When $s_2 = 0$, then $B^{(s_1,s_2),p}(\R^{n_1}\times\R^{n_2})$ is just the usual Besov space $B^{s_1,p}(\R^n)$.

The fundamental properties of these spaces are worked out in \cite{N3}.  In the usual way we may then define anisotropic function spaces on products of open sets in Euclidean space and hence on products of manifolds (with boundary) (see \cite{N0}).   Most of the analysis in this section will occur on the anisotropic spaces $B^{(s_1,s_2),p}(I \times Y)$ and $B^{(s_1,s_2),p}(I \times \Sigma)$ where $I$ is an interval.  Indeed, as we have mentioned, we will be considering maps that smooth in the space variables $Y$ and $\Sigma$. 

We apply the anisotropic function space setup as follows.  Let $I$ be a time interval, bounded or infinite.  We have the smooth configuration space
\begin{align}
  \fC(I \times Y) &= \A(I \times Y) \times \Gamma(I \times \S)\nonumber \\
  &= \Maps(I, \A(Y) \times \Gamma(\S) \times \Omega^0(Y; i\R)) \nonumber \\
  &= \Maps(I, \fC(Y) \times \Omega^0(Y; i\R)) \label{sec2:CIY}
\end{align}
via the correspondence  (\ref{A}).  Here, given any space $X$, we write $\Maps(I,X)$ to denote the space of smooth maps from $I$ into $X$.  Likewise, we write $C^0(I,X)$ to denote the space of continuous maps from $I$ into $X$.  We refer to such maps as paths.
Next, replacing $Y$ with $\Sigma$ in the above, we have
\begin{align}
  \fC(I \times \Sigma) &= \A(I \times \Sigma) \times \Gamma(I \times \S_\Sigma) \nonumber \\
  &= \Maps(I,\fC(\Sigma) \times \Omega^0(\Sigma;i\R)). \label{sec2:CIS}
\end{align}
Since we have the anisotropic function spaces $B^{(s_1,s_2),p}(I \times Y)$ and $B^{(s_1,s_2),p}(I \times \Sigma)$, we get induced topologies on the configuration spaces, their subspaces, and their corresponding tangent spaces.

\begin{Notation}\label{Notation}
  Let $\frak{X}$ be a space of configurations over the space $X_1 \times X_2$, where $X_1$ and $X_2$ are manifolds (each of which is either Euclidean space or a compact manifold possibly with boundary).  Then $B^{(s_1,s_2),p}\frak{X}$ denotes the closure of $\frak{X}$ in the $B^{(s_1,s_2),p}(X_1\times X_2)$ topology.  Likewise, if $\frak{X}$ is a space of configurations over a manifold $X$, we let $B^{s,p}\frak{X}$ denote the closure of $\frak{X}$ in the $B^{s,p}(X)$ topology.  If $E$ is a vector bundle over $X = X_1 \times X_2$, we write $B^{(s_1,s_2),p}(E)$ as shorthand for $B^{(s_1,s_2),p}\Gamma(E)$.  Similar definitions apply for prefixes given by other topologies, e.g., $L^p$ and $C^0$. We also write $E_{\partial X} = E|_{\partial X}$ for the bundle $E$ restricted to $\partial X$.  Finally, we will sometimes refer to just the topology of a configuration, e.g., we will say an element of $B^{s,p}\frak{X}$ belongs to $B^{s,p}$.
\end{Notation}

Thus, letting $M = Y$ or $\Sigma$, we can consider the anisotropic configuration spaces
\begin{equation}
  \fC^{(s_1,s_2),p}(I\times M) := B^{(s_1,s_2),p}\fC(I \times M) \label{eq2:BCIM}
\end{equation}
and we can consider their anisotropic tangent spaces
\begin{align}
 T_{\bullet}\fC^{(s_1,s_2),p}(I \times M)  &= B^{(s_1,s_2),p}(\Omega^1(I \times M; i\R) \oplus \Gamma(I \times \S_M)), \label{CIM1} \\
&= B^{(s_1,s_2),p}\left(\Maps(I, \T_M) \times \Maps(I, \Omega^0(M; i\R))\right), \label{CIM2}
\end{align}
where $\bullet$ is any basepoint.  Here, we used (\ref{sec2:CIY}), (\ref{sec2:CIS}), (\ref{sec2:TY}), and (\ref{sec2:TS}).
These anisotropic spaces induce corresponding topologies on their subspaces, in particular, those subspaces given by the space of paths through $\fC(M)$ and $\T_M$, respectively. Thus, the spaces $B^{(s_1,s_2),p}\Maps(I,\fC(M))$ and $B^{(s_1,s_2),p}(I,\T_M)$ are topologized as subspaces of (\ref{eq2:BCIM}) and (\ref{CIM2}), respectively.  Moreover, all these spaces are the completions of spaces of smooth configurations in the $B^{(s_1,s_2),p}(I \times M)$ topology, as is consistent with Notation \ref{Notation}.

The above definitions work out nicely because the spaces we are topologizing are affine spaces.  Suppose we now wish to topologize spaces that are not linear, such as the space of paths through $\M$ and $\L$. For this, we can describe more general path spaces in anisotropic topologies as follows. By the trace theorem, Theorem \ref{ThmATrace}, we have a trace map
\begin{equation}
  r_t: B^{(s_1,s_2),p}(\R \times M) \to B^{s_1+s_2-1/p,p}(\{t\} \times M), \qquad t \in \R \label{2.1-trace}
\end{equation}
for all $s_1 > 1/p$.  Moreover, this trace map is continuous in $t$, in other words, we have the inclusion
\begin{equation}
  B^{(s_1,s_2),p}(\R \times M) \hookrightarrow C^0(\R; B^{s_1+s_2-1/p,p}(M)).
\end{equation}
Thus, we have the following well-defined spaces of paths in anisotropic topologies:

\begin{Def}\label{DefAMaps}
  Let $s_1 > 1/p$, $s_2 \geq 0$, and $p \geq 2$.  Then we can define the following spaces \label{p:Maps}
  \begin{align}
    \Maps^{(s_1,s_2),p}(I, \fC(M)) &= B^{(s_1,s_2),p}\Maps(I,\fC(M)) \label{DefAMapsC}\\
    \Maps^{(s_1,s_2),p}(I, \T_M) &= B^{(s_1,s_2),p}\Maps(I,\T_M) \label{DefAMaps0}\\
    \Maps^{(s_1,s_2),p}(I, \M) &= \{\gamma \in \Maps^{(s_1,s_2),p}(I, \T_Y) : \gamma(t) \in \M^{s_1+s_2-1/p,p}, \textrm{ for all } t \in I\} \label{DefAMaps1} \\
    \Maps^{(s_1,s_2),p}(I, \L) &= \{\gamma \in \Maps^{(s_1,s_2),p}(I, \fC(\Sigma)) : \gamma(t) \in \L^{s_1+s_2-1/p,p}, \textrm{ for all } t \in I\}. \label{DefAMaps2}
  \end{align}
  If $s_2=0$, we write $\Maps^{s_1,p}$ instead of $\Maps^{(s_1,0),p}$. Similar definitions apply if we extend $\fC(M)$ and $\T_M$ in the above by vector bundles over $M$, in particular, $\Omega^0(M;i\R)$.  Finally, if $\X$ is a subspace of $\T_M^{s,p}$ for some $s$, define
  \begin{equation}
    \Maps^{(s_1,s_2),p}(I, \X) = B^{(s_1,s_2),p}\{z \in \Maps(I,\T_M) : z(t) \in \X, \textrm{ for all }t \in I\}, \label{DefAMaps3}\\
  \end{equation}
  the $B^{(s_1,s_2),p}(I, M)$ closure of the space of paths in $\Maps(I,\T_M)$ which take values in $\X$. Note (\ref{DefAMaps3}) generalizes definition (\ref{DefAMaps0}).
\end{Def}

The first two definitions above are just a change of notation since we have already defined the right-hand side.  However, for the next two definitions, there is a subtle point in how we defined these spaces as compared to the previous ones.  Observe that for (\ref{DefAMapsC}) and (\ref{DefAMaps0}), the space $\Maps^{(s_1,s_2),p}(I, \T_M)$, say, is \textit{defined} to be the closure of a space of smooth paths, namely, the closure of $\Maps(I,\T_M)$ in the $B^{(s_1,s_2),p}(I \times M)$ topology.  On the other hand, in (\ref{DefAMaps2}) say, we begin with a path $\gamma$ in the $B^{(s_1,s_2),p}(I \times \Sigma)$ topology, and we impose a restriction on its trace, namely that it always lie in $\L^{s_1+s_2-1/p,p}$.  The resulting space is a priori a larger space than the closure of the space of smooth paths through the space $\L$, i.e., we have
\begin{equation}
  B^{(s_1,s_2),p}\Maps(I,\L) \subseteq \Maps^{(s_1,s_2),p}(I,\L). \label{difftop}
\end{equation}
Indeed, to prove the reverse inequality, one would have to approximate an arbitrary path $\gamma$ in the space $\Maps^{(s_1,s_2),p}(I,\L)$ by a smooth path (both in space \textit{and} time), which because of the nonlinearity of $\L$, is not obvious how to do.  Indeed, even though we showed in  \cite{N1} that the smooth monopole space $\L$ is dense in $\L^{s_1+s_2-1/p}$, this says nothing about the temporal regularity of an approximating sequence to the space of paths. However, by Theorems \ref{ThmPathL1} and \ref{ThmPathL2}, it does turn out to be the case that we can approximate the space of paths by smooth paths in the appropriate range of $s_1$, $s_2$, and $p$, so that we have equality in (\ref{difftop}).\\

We now have all the appropriate definitions and notation in place for our anisotropic path spaces. Our main task in the rest of this section is to prove the analogous results in \cite{N1} for the space of paths through $\L$, i.e., we want to prove that the space $\Maps^{(s_1,s_2),p}(I,\L)$ is a Banach manifold for a large range of parameters and it has chart maps with smoothing properties as described in the discussion at the beginning of this section.  These chart maps are defined from the charts for $\L$ in a ``slicewise fashion" along a path in $\Maps^{(s_1,s_2),p}(I,\L)$.  We thus need to understand the analytic properties of such slicewise operators.

\subsection{Slicewise Operators on Paths}

Unless stated otherwise, from now on, we always assume $p \geq 2$. In the following, we will be estimating the operators studied in \cite{N1} acting slicewise on the space of paths from an interval $I$ into some target configuration space.  More precisely, we have the following definition.

\begin{Def}\label{DefSw}
  Let $O(t)$ be a family of operators acting on a configuration space $\frak{X}$, $t \in I$.  We write $\w{O(t)}$ to denote the slicewise operator associated to the family $O(t)$, that is, $\w{O(t)}$ applied to a path $\g: I \to \frak{X}$ yields the path \label{p:slicewise}
  \begin{align}
    \w{O(t)}\g = \Big(t \mapsto O(t)\g(t)\Big)_{t \in I}.
  \end{align}
\end{Def}

\begin{Notation}
  Note that given a path of configurations $\gamma$, if we write $\g = \g(t)$, it is ambiguous whether we mean the whole path as a function of $t$ or just the single configuration at time $t$.  The above hat notation mitigates this ambiguity.  Moreover, we will mainly be considering the case when $O(t) \equiv O$ is time-independent.  In this case, the $\w{}$\; notation therefore just serves as a notational reminder, although in certain cases, such as when $O$ is a differential operator, we will sometimes just write $O$ to denote $O$ acting slicewise, which is standard notational practice in this case.\\
\end{Notation}

We now proceed to estimate the linear operators of interest to us when they act slicewise.  Recall that we defined the notion of a \textit{local straightening map} in Definition I.B.2 in an abstract framework, and showed how local straightening maps yield natural local charts for Banach submanifolds of a Banach space. In Lemma I.3.9, we defined a local straightening map $F_{\Sigma, \co}$ for $\L^{s-1/p,p}$ at a configuration $r_\Sigma\co \in \L^{s-1/p,p}$, where $\co \in \M^{s,p}$. In this way, we used the local straightening map $F_{\Sigma,\co}$ to obtain chart maps for the Banach submanifold $\L^{s-1/p,p} \subset \fC^{s-1/p,p}(\Sigma)$ in Theorem I.3.4.  We also deduced some important properties of the resulting chart maps, $E_{r_\Sigma\co}$, namely that they are defined on large domains (i.e. $B^{s',p}(\Sigma)$ open balls) and that the nonlinear portion of the chart map $E^1_{r_\Sigma\co}$ smooths by a derivative. When studying the space of paths $\Maps^{s,p}(I,\L)$, our first goal is to show that the slicewise map $\w{F_{\Sigma,\co}}$ gives a local straightening map within a neighborhood of a constant path identically equal to $r_\Sigma\co \in \L^{s,p}$.  From this, the local straightening map yields for us a chart map for $\Maps^{s,p}(I,\L)$ at a constant path.  Later, we will see how to ``glue together" these chart maps for constant paths on small time intervals to obtain a chart map at an arbitrary path in $\Maps^{s,p}(I,\L)$.  

One of the operators that arises in the definition of $F_{\Sigma, \co}$, as seen in (I.3.38) and (I.3.29), is the operator $Q_\co$, defined in (I.3.21). Let us review this operator.  The operator $Q_\co$ is constructed out of the Hessian operator $\H_{\co}$ (more precisely, its inverse on suitable domains), the projection $\Pi_{\K_{\co}^{s,p}}$, and a pointwise quadratic multiplication map $\q$.  Let us briefly review these maps.  Recall that for any configuration $\c \in \fC^{s,p}(Y)$, the operator
\begin{equation}
  \H_\c: \T^{s,p}\to\T^{s-1,p}
\end{equation}
is given by
\begin{equation}
  \H_{(B,\Psi)} = \begin{pmatrix}
    *_Yd & 2i\Im\rho^{-1}(\cdot\Phi^*)_0\\
    \rho(\cdot)\Phi & D_{B}
  \end{pmatrix}.
\end{equation}
It is a first order formally self-adjoint operator and it has kernel $T_\c\fM^{s,p}$ whenever $\c \in \fM^{s,p}$.  By Lemma I.3.3, given any $\co \in \fM^{s,p}$, we can choose a subspace $X_\co^{s,p} \subset \T^{s,p}$ complementary to $T_\c\fM^{s,p}$.  For such $X_\co^{s,p}$, we have by Proposition I.2.2 that
\begin{align}
  \H_\co: X_\co^{s,p} & \to \K^{s-1,p}_\co \label{eq2:HKiso0}\\
  \H_\co: X_\co^{s+1,p} & \to \K^{s,p}_\co, \qquad X_\co^{s+1,p} = X_\co^{s,p} \cap \T^{s+1,p} \label{eq2:HKiso}
\end{align}
are isomorphisms.  Here, the subspace $\K^{s,p}_{\co} \subset \T^{s,p}$ is a complement to the tangent space of the gauge orbit at $\co$ in $\T^{s,p}$, see Lemma I.2.3.  The map
\begin{equation}
  \Pi_{\K_{\co}^{s,p}}: \T^{s,p} \to \K^{s,p}_{\co} \label{eq2:PiK}
\end{equation}
is a bounded projection onto this space, see (I.2.45).  Finally, the map $\q$ arises from the quadratic multiplication that occurs in the map $SW_3$, see (I.3.16).  Thus, both $\q$ and $\w{\q}$ are pointwise multiplication operators, and their mapping properties are controlled by the function space multiplication theorem, Theorem \ref{ThmMult}.

Thus, when we consider the above operators slicewise, the main operators we need to understand are
\begin{equation}
\w{ (\H_{\co} |_{X_\co^{s+1,p}} )^{-1} }, \quad \w{\Pi_{\K^{s,p}_{ \co } }}.  \label{swoperators}
\end{equation}
However, since the operators in (\ref{swoperators}) are time-independent, estimating them on (anisotropic) Besov spaces is not difficult.  Indeed, we have the following general lemma:

\begin{Lemma}\label{LemmaSW}
Let $s_1 > 0$, and $s_2,s_2' \geq 0$, and let $M$ be a compact manifold. Let $T: C^\infty(M) \to C^\infty(M)$ be a linear operator that extends to bounded operators
    \begin{align}
  T: H^{s_2,p}(M) & \to H^{s_2+s_2',p}(M) \label{OpT1}\\
  T: B^{s_1+s_2,p}(M) & \to B^{s_1+s_2+s_2',p}(M). \label{OpT2}
\end{align}
Then the slicewise operator
\begin{equation}
  \widehat{T}: B^{(s_1,s_2),p}(I \times M) \to B^{(s_1,s_2+s_2'),p}(I \times M) \label{whT}
\end{equation}
is bounded and the operator norm of (\ref{whT}) is bounded in terms of the operator norms of (\ref{OpT1}) and (\ref{OpT2}).

\end{Lemma}

\Proof The crucial property we need is the so-called ``Fubini property" of Besov spaces (see \cite{Tr1}).  Namely, for any $s > 0$, the Besov space $B^{s,p}(I \times M)$ can be written as the intersection
\begin{equation}
  B^{s,p}(I \times M) = L^p(I, B^{s,p}(M)) \cap L^p(M, B^{s,p}(I)). \label{fubB}
\end{equation}

In other words, we can separate variables so that a function of regularity of order $s$ on $I \times M$ is a function that has regularity of order $s$ in $I$ and $M$, separately, in the Besov sense.  This Fubini property automatically implies one for anisotropic Besov spaces, and we have from (\ref{liftAB}) that
\begin{equation}
  B^{(s_1,s_2),p}(I \times M) = L^p(I, B^{s_1+s_2,p}(M)) \cap H^{s_2,p}(M, B^{s_1,p}(I)). \label{fubAB}
\end{equation}
Recall that $H^{s_2,p}(M)$ is the fractional Sobolev space of functions whose fractional derivatives up to order $s_2$ belong to $L^p(M)$.\footnote{The space $H^{s_2,p}(M)$ is also known as a \textit{Bessel potential} space. See \cite{N1} for a precise definition in the general case.}  In other words, given an elliptic differential operator $D_{s_2}$ on $M$ of order $s_2$, we can define the norm on $H^{s_2,p}(M)$ by
$$\|f\|_{H^{s_2,p}(M)} = \|f\|_{L^p(M)} + \|D_{s_2}f\|_{L^p(M)}.$$
Thus, the space $H^{s_2,p}(M, B^{s,p}(I))$ appearing in (\ref{fubAB}) is the space of functions $f$ such that both $f$ and $D_{s_2}f$ belong to $L^p(M,B^{s,p}(I))$.

We want to show that (\ref{whT}) is bounded.  We proceed by using the decomposition (\ref{fubAB}).  First, the boundedness of $T: B^{s_1+s_2,p}(M) \to B^{s_1+s_2+s_2',p}(M)$ implies the boundedness of
\begin{equation}
  \w{T}: L^p(I,B^{s_1+s_2,p}(M)) \to L^p(I,B^{s_1+s_2+s_2',p}(M)). \label{Fub-bound1}
\end{equation}
It remains to show that the map
\begin{equation}
  \w{T}: H^{s_2,p}(M,B^{s_1,p}(I)) \to H^{s_2+s_2',p}(M,B^{s_1,p}(I)) \label{Fub-bound2}
\end{equation}
is bounded. To do this, we first show that the space $H^{s_2,p}(M, B^{s,p}(I))$, as defined above, is also equal to the space $B^{s,p}(I, H^{s_2,p}(M))$ which is defined as follows.

Given any Banach space $X$, we can define Besov spaces with values in $X$ as follows.  Namely, for $s > 0$, we can define the space $B^{s,p}(\R,X)$ to be the closure of the space of infinitely Fr\'echet differentiable functions $f: \R \to X$ under the norm \cite[p. 20]{Am} \cite[p. 110]{Tr1}
\begin{equation}
  \|f\|_{B^{s,p}(\R,X)} := \|f\|_{L^p(\R,X)} + \left(\int_0^\infty \left(\int_\R|h^{-s}(\tau_h - \id)^m f(t)|_X^p dt \right)\frac{1}{h}dh\right)^{1/p},
\end{equation}
where $m > s$ is any integer and $\tau_h$ is the translation operator $(\tau_hf)(t) = f(t + h)$. When $X = \R$, this recovers the usual scalar valued Besov spaces.  Observe that for $X = H^{s_2,p}(M)$, we have $B^{s,p}(\R, H^{s_2,p}(M)) = H^{s_2,p}(M, B^{s,p}(\R))$ since the operators $h^{-s}(\tau_h - \id)^m$ and $D_{s_2}$ commute.  From the above, we now define $B^{s,p}(I,X)$ to be the restrictions of elements of $B^{s,p}(\R,X)$ to the domain $I$.

Thus, showing (\ref{Fub-bound2}) is the same thing as showing
\begin{equation}
  \w{T}: B^{s,p}(I,H^{s_2,p}(M)) \to B^{s,p}(I,H^{s_2+s_2',p}(M)). \label{Fub-bound3}
\end{equation}
To show (\ref{Fub-bound2}), we proceed by interpolation \cite{Tr}.  Namely, Besov spaces are interpolation spaces of Sobolev spaces, i.e.,
\begin{equation}
  B^{s,p}(\R, X) = (H^{k_0,p}(\R, X), H^{k_1,p}(\R, X))_{\theta,p}, \label{interp}
\end{equation}
for any nonnegative integers $k_0, k_1$ and $0<\theta<1$ such that $(1 - \theta)k_0 + \theta k_1 = s$.  Here, $(\cdot,\cdot)_{\theta,p}$ is the real interpolation functor and $H^{k,p}(\R,X)$ denotes the Sobolev space of functions\footnote{Technically, we should denote this space $W^{k,p}(\R,X)$, but for all Banach spaces $X$ we encounter, $H^{k,p}(\R,X)$ as it is usually defined (see \cite{N0}) agrees with $W^{k,p}(\R,X)$.} with values in $X$ whose derivatives up to order $k$ belong to $L^p(X)$.

Thus, to establish that (\ref{Fub-bound3}) is bounded, from interpolation, it suffices to establish that
\begin{equation}
  \w{T}: H^{k,p}(I,H^{s_2,p}(M)) \to H^{k,p}(I,H^{s_2+s_2',p}(M)) \label{Fub-bound4}
\end{equation}
is bounded for all integer $k \geq 0$.  For $k = 0$, this follows trivially from (\ref{OpT1}), and for $k \geq 1$, this also follows from (\ref{OpT1}) by commuting derivatives with $T$, since $T$ is linear and time-independent (observe that if $T$ were time-dependent but also smooth, this argument would follow too).  Thus, this proves that (\ref{Fub-bound4}) is bounded, which finishes the proof that (\ref{whT}) is bounded.\End

The above lemma tells us that time-independent slicewise operators on $I \times M$ can be estimated in terms of their mapping properties on $M$.  In particular, we can now easily estimate the slicewise operators in (\ref{swoperators}) on $I \times Y$ because we know how they act on $Y$ from the analysis carried out in \cite{N1}. We have the following corollary:

\begin{Corollary}\label{CorSW}
Let $s > \max(3/p,1/2)$ and $\co \in \M^{s,p}$.
\begin{enumerate}
  \item We have bounded maps
  \begin{align}
    \widehat{\left(\H_{\co}|_{X_{\co}^{s+1,p}}\right)^{-1}}\widehat{\Pi_{\K_{\co}^{s,p}}}: \Maps^{s,p}(I,\T) & \to \Maps^{(s,1),p}(I,\T). \label{hatHmap1}
  \end{align}
  \item Assume in addition that $\co \in \M^{s_1+s_2,2}$, where $s_1 \geq 3/2$ and $s_2 \geq 0$.  Then we have
    \begin{align}
    \widehat{\left(\H_{\co}|_{X_{\co}^{s_1+s_2+1,2}}\right)^{-1}}\widehat{\Pi_{\K_{\co}^{s,p}}}: \Maps^{(s_1,s_2),2}(I,\T) & \to \Maps^{(s_1,s_2+1),2}(I,\T). \label{hatHmap2}
  \end{align}
\end{enumerate}
\end{Corollary}

\Proof By the above lemma, to establish (\ref{hatHmap1}), it suffices to show that the maps
\begin{align*}
  \left(\H_{\co}|_{X_{\co}^{s+1,p}}\right)^{-1} \Pi_{\K_{\co}^{s,p} }: L^p\T & \to H^{1,p}\T \\
  \left(\H_{\co}|_{X_{\co}^{s+1,p}}\right)^{-1} \Pi_{\K_{\co}^{s,p} }: \T^{s,p} & \to \T^{s+1,p}
\end{align*}
are bounded. (By the first line above, we of course mean that the operator in the second line extends to a bounded operator on the stated domain and ranges.  Such notation will always be understood from now on.) The second line follows from (\ref{eq2:PiK}) and (\ref{eq2:HKiso}).  The first line follows from Lemma I.3.5(ii).  Likewise,  for $\co \in \M^{s_1+s_2+1,2}$, (\ref{hatHmap2}) follows from
\begin{align*}
  \left(\H_{\co}|_{X_{\co}^{s+1,p}}\right)^{-1} \Pi_{\K_{\co}^{s,p}} : L^2\T & \to H^{1,2}\T\\
  \left(\H_{\co}|_{X_{\co}^{s+1,p}}\right)^{-1} \Pi_{\K_{\co}^{s,p}} : \T^{s_1+s_2,2} & \to \T^{s_1+s_2+1,2},
\end{align*}
which also follows from I.3.5(ii) and Proposition I.2.2.\End


The above estimates for slicewise operators will now allow us to pass from a local chart map for $\L$ to a local chart map for the space of paths through $\L$ centered at a path through $\L$.  More precisely, we construct such charts for the space of paths through $\L$ in the Besov topology $B^{s,p}(I \times \Sigma)$ topology, and this shows that the space $\Maps^{s,p}(I,\L)$ is a Banach submanifold of $\Maps^{s,p}(I,\fC(\Sigma))$. To do this, we prove the following fundamental lemma below, whose main idea we can describe as follows:

Recall from earlier discussion that we construct a chart map for $\L^{s-1/p,p}$ from the data of a \textit{local straightening map}, see Definition I.B.2.  Specifically, to a configuration $r_\Sigma(B_0,\Psi_0) \in \L^{s-1/p,p}$, the local straightening map involved is the map $F_{\Sigma,\co}$ of Lemma I.3.9.  It readily follows from this construction that to a constant path $\g(t) \equiv (B_0(t), \Psi_0(t))$ through $\Maps^{s,p}(I,\L)$, the slicewise operator $\w{F_{\Sigma,\co}}$ is a local straightening map for $\g(t)$ for the Banach manifold $C^0(I, \L^{s-1/p,p})$, which induces for us a time-independent slicewise chart map for $C^0(I, \L^{s-1/p,p})$ in the natural way.  Our first main goal is to show that this construction is compatible with the $B^{s,p}(I \times \Sigma)$ topology in the sense that all maps involved preserve the subspace $B^{s,p}(I \times \Sigma) \subset C^0(I, B^{s-1/p,p}(\Sigma))$, which we show in Lemma \ref{LemmaFSigmaHat} below.  We will then use this result in Theorem \ref{ThmPathL1} to construct chart maps for $\Maps^{s,p}(I,\L)$.  However, proving the aforementioned preservation of $B^{s,p}(I \times \Sigma)$ regularity is a nontrivial matter for two reasons.  First, the space $B^{s,p}(I \times \Sigma)$ has fractional regularity (in time), which makes it difficult to estimate nonlinear operators on this space.  Second, we want to show that the $B^{s,p}(I \times \Sigma)$ preservation holds on ``large" open sets, i.e. on an open set containing a ball defined with respect to the $C^0(I, B^{s',p}(\Sigma))$ topology, where $s' < s - 1/p$.  This is so that the local Banach manifold chart maps we obtain for $\Maps^{s,p}(I,\L)$ will contain a $C^0(I, B^{s',p}(\Sigma))$ neighborhood of an element of  $\Maps^{s,p}(I,\L)$, which is key in the proof of compactness in Theorem B. (In a few words, this is so that sequences of bounded configurations in $\Maps^{s,p}(I,\L)$, which we know will converge strongly along a subsequence in the weaker topology $C^0(I, B^{s',p}(\Sigma))$, converges within the range of a fixed chart map.)  Moreover, the fact that chart maps are defined on large neighborhoods is what allows us to construct chart maps for $\Maps^{s,p}(I,\L)$ at nonconstant paths via approximation by constant paths in Theorem \ref{ThmPathL1}.

We now prove our fundamental lemma, which is the slicewise analog of Lemma I.3.9.

\begin{Lemma}\label{LemmaFSigmaHat}
  Let $s > \max(3/p,1/2+1/p)$ and let $r_\Sigma\co \in \L^{s-1/p,p}$, where $\co \in \M^{s,p}$.  We have
  $$C^0(I, \T_\Sigma) = C^0(I, T_{r_\Sigma\co}\L^{s-1/p,p}) \oplus C^0(I, J_\Sigma T_{r_\Sigma\co}\L^{s-1/p,p})$$
  and we can define the map
  \begin{align}
  \widehat{F_{\Sigma,\co}} : \mathcal{V}_\Sigma & \to  C^0(I, T_{r_\Sigma\co}\L^{s-1/p,p}) \oplus C^0(I, J_\Sigma T_{r_\Sigma\co}\L^{s-1/p,p}) \nonumber \\
  z = (z_0,z_1) & \mapsto (z_0, z_1 - \widehat{ r_\Sigma }\widehat{ E^1_{\co} }(\widehat{ P_\co } z_0)), \label{eq4:FSigma}
\end{align}
where $\mathcal{V}_\Sigma \subset C^0(I, \T_\Sigma)$ is an open subset containing $0$, $E^1_\co$ is defined as in Theorem I.3.2, and $P_\co$ is the Poisson operator defined as in Theorem I.2.1. For any $\max(1/2,2/p) < s' \leq s-1/p$, we can choose $\mathcal{V}_\Sigma$ to contain a $C^0(I,B^{s',p}(\Sigma))$ ball, i.e., there exists a $\delta > 0$, depending on $r_\Sigma\co$, $s'$, and $p$, such that
$$\V_\Sigma \supseteq \{z \in C^0(I, \T_\Sigma) : \|z\|_{C^0(I,B^{s',p}(\Sigma))} < \delta\}.$$
Moreover, we can choose $\delta$ such that the following hold:
\begin{enumerate}
  \item We have $\w{F_{\Sigma,\co}}(0)=0$ and $D_0\w{F_{\Sigma,\co}} = \id$.  For $V_\Sigma$ sufficiently small, $\w{F_{\Sigma,\co}}$ is a local straightening map for $C^0(I,\L^{s',p})$ within $V_\Sigma$.
  \item The maps $\w{F_{\Sigma,\co}}$ and $\w{F_{\Sigma,\co}}^{-1}$ preserve $B^{s,p}(I\times\Sigma)$ regularity.  Moreover, the term $\widehat{r_\Sigma }\widehat{ E^1_{\co} }(\widehat{ P_\co })$ maps $B^{s,p}(I \times \Sigma)$ configurations to $B^{(s,1-1/p-\eps),p}(I\times\Sigma)$, where $\eps > 0$ is arbitrary.
  \item We can choose $\delta$ uniformly for $r_\Sigma\co$ in a sufficiently small $B^{s',p}(\Sigma)$ neighborhood of any configuration in $\L^{s-1/p,p}$.
\end{enumerate}
\end{Lemma}

\Proof By Lemma I.3.9, we know that $F_{\Sigma,\co}$ is a local straightening map for $\L^{s-1/p,p}$.  Based on that lemma, it easily follows that $\w{F_{\Sigma,\co}}$ is a local straightening map for $C^0(I,\L^{s-1/p,p})$, thus establishing (i).  It is (ii) that mainly needs verification, which we do by estimating the operators occurring in $\w{F_{\Sigma,\co}}$ one by one.

First, we show the boundedness of the slicewise Calderon projection $\widehat{P^+_{\co}}$ in the $B^{s,p}(I \times \Sigma)$ topology, so that the decomposition $x = (x_0,x_1)$ above in the $C^0(I, B^{s-1/p,p}(\Sigma))$ topology
induces one in the $B^{s,p}(I \times \Sigma)$ topology.  Here, $P^+_{\co}$ is the Calderon projection as defined in Theorem I.2.1, and it yields a bounded map
\begin{align}
  P^+_\co: \T^{t,p}_\Sigma \to B^{t,p}(T_{r_\Sigma\co}\L^{s-1/p,p}), \qquad 0 \leq t \leq s+1-1/p, \label{eq:caldproj}
\end{align}
since $\co \in \M^{s,p}$.  Here, by a slight abuse of notation with regard to Notation \ref{Notation}, we define $B^{t,p}(T_{r_\Sigma\co}\L^{s-1/p,p})$ to be the intersection of $T_{r_\Sigma\co}\L^{s-1/p,p}$ with $B^{t,p}(\Sigma)$ configurations in case $t \geq s - 1/p$; otherwise, we take the $B^{t,p}(\Sigma)$ closure. The map $P^+_{\co}$ also extends to a bounded operator on Sobolev spaces
\begin{align}
  P^+_\co: H^{t,p}\T_\Sigma \to H^{t,p}(T_{r_\Sigma\co}\L^{s-1/p,p}), \qquad 0 \leq t < s+1-1/p, \label{eq:caldproj2}
\end{align}
where $H^{t,p}(T_{r_\Sigma\co}\L^{s-1/p,p})$ is defined as above, see Remark I.3.2. From this, since $P^+_\co$ is time-independent, one can apply Lemma \ref{LemmaSW} to conclude, in particular, that
\begin{equation}
  \widehat{P^+_\co}: \Maps^{s,p}(I,\T_\Sigma) \to \Maps^{s,p}(I,T_{r_\Sigma\co}\L^{s-1/p,p}),  \label{eq:caldproj3}
\end{equation}
is bounded with kernel $\Maps^{s,p}(I, J_\Sigma T_{r_\Sigma\co}\L^{s-1/p,p})$.

Proceeding as above, we also obtain a bounded Poisson operator 
\begin{align}
  \widehat{P_\co}: \Maps^{s,p}(I,\T_\Sigma) & \to \Maps^{(s,1/p),p}(I, T_{\co}\M^{s,p}) \nonumber \\
  & \qquad = \{z \in \Maps^{(s,1/p),p}(I, \T) : z(t) \in T_{\co}\M^{s,p},\quad \forall t \in I\}
\end{align}
and its range is contained in $C^0(I,L^\infty\T)$, due to the embedding
$$\Maps^{(s,1/p),p}(I,\T) \hookrightarrow C^0(I, B^{s,p}(Y)) \hookrightarrow C^0(I,L^\infty(Y)).$$
(Here, when we apply Lemma \ref{LemmaSW} to $\widehat{P_\co}$, it makes no difference that we map configurations on $I \times \Sigma$ to configurations on $I \times Y$, i.e., Lemma \ref{LemmaSW} is unchanged if the domain and range manifolds are different.)

Next, we estimate the mapping properties of $\widehat{E^1_\co}$.  Namely, we show that
\begin{align}
  \widehat{E^1_\co}: \Maps^{(s,1/p),p}(I, T_{\co}\M^{s,p}) & \dashrightarrow \Maps^{(s,1),p}(I,\T),  \label{E1Lhat}
\end{align}
is bounded, i.e. $\widehat{E^1_\co}$ smooths by $1-1/p$ derivatives in the $Y$ directions. Here, we restrict the domain of (\ref{E1Lhat}) so that it lies inside the set
\begin{equation}
  \Maps^{(s,1/p),p}(I, U) = \{z \in \Maps^{(s,1/p),p}(I, T_{\co}\M^{s,p}): z(t) \in U\} \label{eq:mapsU}
\end{equation}
where $U \subset T_\co\M^{s,p}$, as defined in Theorem I.3.2, is a domain on which $E^1_\co$ is defined. This is where our slicewise estimates made in Lemma \ref{LemmaSW} and Corollary \ref{CorSW} come into play.  By (I.3.29), we have
\begin{align}
  \widehat{E^1_\co}(z) = -\w{Q_\co}\big(\w{F_\co}^{-1}(z), \w{F_\co}^{-1}(z) \big), \label{eq:E1hatQ}
\end{align}
where $F_\co$ is a local straightening map for $\co \in \M^{s,p}$, see Lemma I.3.4 and Theorem I.3.2.

To estimate $\widehat{E^1_\co}$, we first estimate $\w{F_\co}^{-1}$.  This is the most difficult step of all.  We want to use the ideas from Lemma \ref{LemmaSW} to conclude that the smooth time-independent map $\w{F_\co}^{-1}$ preserves the $B^{(s,1/p),p}(I\times Y)$ topology.  Namely, using the Fubini property, we can write
\begin{equation}
  B^{(s,1/p),p}(I\times Y) = L^p(I,B^{s+1/p,p}(Y)) \cap B^{s,p}(I,H^{1/p,p}(Y)). \label{sec2:fubSigma}
\end{equation}
Since $\co \in \M^{s,p}$, the proofs of Lemma I.3.5, Theorem I.3.2, and Corollary I.3.1 show that $F_\co^{-1}$ preserves $B^{t,p}(Y)$ regularity for $1/p \leq t \leq s+1$, hence for $t = s + 1/p$ in particular.  Thus, we have trivially that the time-independent operator $\w{F_\co}^{-1}$ preserves $L^p(I,B^{s+1/p,p}(Y))$ regularity.  The nontrivial step is to show that $\w{F_\co}^{-1}$ preserves $B^{s,p}(I,H^{1/p,p}(Y))$ regularity.  To show this, we want to interpolate between the estimates
\begin{equation}
  \w{F_\co}^{-1}: H^{k,p}(I,H^{1/p,p}(Y)) \dashrightarrow H^{k,p}(I,H^{1/p,p}(Y)), \qquad k \geq 0, \label{sec2:NLinterp}
\end{equation}
where we abuse notation by letting $H^{k,p}(I,H^{1/p,p}(Y))$ denote the closure of $\Maps(I,\T)$ in said topology. Here the domain of (\ref{sec2:NLinterp}) will be specified in a moment.

We want to use interpolation because fractional regularity in time is difficult to grasp; on the other hand, integer Sobolev spaces, such as those appearing in (\ref{sec2:NLinterp}) are amenable to estimates via the Leibnitz rule.  The difficulty of course is that $\w{F_\co}^{-1}$ is nonlinear so that Lemma \ref{LemmaSW} does not apply.  Additionally, $\w{F_\co}^{-1}$ is a highly nonexplicit map, being the slicewise inverse of the operator $F_\co$ and hence defined by the inverse function theorem.  In \cite{N1}, we were able to estimate $F_\co^{-1}$ (see the proof of Theorem I.3.2) because it differed from the identity map by a smoothing operator, but in this case, since we only have smoothing in the space directions, we cannot carry over such a perturbative argument. What saves us in our situation is that there is a theory of interpolation of Lipschitz operators due to Peetre \cite{Pe}.  Moreover, such an interpolation theory extends to operators which are only locally defined (as is the case for $\w{F_\co}^{-1}$) under the appropriate conditions. Here, the relevant theorem is Theorem 14.8 of \cite{N0}.  We can apply this theorem
to $ \w{F_\co}^{-1}$ by verifying that the following properties hold:
\begin{enumerate}\renewcommand{\theenumi}{(\Roman{enumi})}
  \item The map $F_{\co}^{-1}$ is a bounded operator from a subset $\tilde V$ of $H^{1/p,p}\T$ into $H^{1/p,p}\T$.  (This implies that the map $\w{F_{\co}}^{-1}$ maps $L^\infty(I, \tilde V)$ into $L^\infty(I,H^{1/p,p}\T)$.  We then have (\ref{sec2:NLinterp}) on this domain since $F_\co^{-1}$ is smooth.)
  \item On the same subset $\tilde V$ above, the map $F_{\co}^{-1}$ and all its Fr\'echet derivatives are Lipschitz. (Hence the corresponding statement is true for $\w{F_{\co}}^{-1}$ on $L^\infty(I, \tilde V)$ by the Leibnitz rule).
  \item The vector valued Gagliardo-Nirenberg inequality holds: for all integers $0 < m < n$, we have
$$\|\partial_t^m f\|_{L^r(I, \X)} \leq \|f\|_{L^q(I,\X)}^{1-\theta} \|\partial_t^n f\|_{L^p(I,\X)}^{\theta},$$
where $\X = H^{1,p}(Y)$, and $r$, $p$, $q$ are any numbers satisfying
$$1 < r,p,q \leq \infty, \quad \frac{1}{r} = \frac{\theta}{q} + \frac{1-\theta}{p}, \quad 0 < \theta = m/n < 1.$$
\end{enumerate}\renewcommand{\theenumi}{(\roman{enumi})}

Let us explain why these properties hold. For (I), the statement for $F_{\co}^{-1}$ holds due to Corollary I.3.1.  For (II), we use the fact that $F_{\co}^{-1}$, being the inverse of an analytic (in fact quadratic) map $F_\co$, is itself analytic (i.e. it has a local power series expansion) and hence so is $\w{F_{\co}}^{-1}$.  This uses the fact that the inverse function theorem holds in the analytic category of functions, see e.g. \cite{De}.  Consequently, the fact that $\w{F_{\co}}^{-1}$ is analytic means we have estimates for all the derivatives of $F_{\co}^{-1}$ on some fixed small neighborhood (just as one would have for a holomorphic function of one complex variable, using Cauchy's integral formula or its local power series expansion), which allows us to establish (II). These statements follow from Proposition 21.3 of \cite{N0}.\footnote{In fact, since $\w{F_\co}$ is just a quadratic map, one can understand the form of power series expansion of $\w{F_\co}^{-1}$ sufficiently well so that one can use avoid the nonlinear interpolation method we have used to establish the boundedness of $\w{F_\co}$ on the $B^{(s,1/p),p}(I \times Y)$.  Using this method, however, has the unfortunate consequence that the radius of convergence of the power series of $\w{F_\co}^{-1}$ (about $0$) depends on $s$ and $p$, and so the set $V_\Sigma$ of the lemma would depend on $s$ and $p$ (i.e. it may shrink with $s$ and $p$).  This would result in a somewhat awkward proof of Theorem \ref{ThmPathL1}(iv) and Theorem \ref{ThmA'} later on. Hence, we adhere to the interpolation method since in some sense it is the ``optimal" method, even though for the purposes of this paper, we can get away with other methods that use less machinery.} For (III), the case when $\X = \R$ is the classical scalar Gagliardo-Nirenberg inequality.  As it turns out, by the work of \cite{SS}, the general vector valued case holds for \textit{arbitrary} Banach spaces $\X$ .  Thus (III) holds in particular for $\X = H^{1/p}\T$.

Altogether, from the above, we can apply Theorem 14.8 in \cite{N0} to our situation, where in that theorem, $\cZ = H^{1,p}\T$ and the set $U_r$ is the set $L^\infty(I,\tilde V)$ for sufficiently small $\tilde V$.
Here the Lipschitz hypothesis of Theorem 14.8 is satisfied by the same argument as in
\cite[p. 330]{Pe}, because of (III).  Thus, interpolation between the estimates (\ref{sec2:NLinterp}) shows that
\begin{equation}
  \w{F_\co}^{-1}: B^{s,p}(I,H^{1/p,p}(Y)) \dashrightarrow B^{s,p}(I,H^{1/p,p}(Y)) \label{sec2:NLinterp3}
\end{equation}
is bounded, where the domain of (\ref{sec2:NLinterp3}) is the open subset
\begin{equation}
  \{z \in \Maps^{(s,1/p),p}(I,\T): z(t) \in \tilde V, \textrm{ for all }t \in I\} \label{eq:mapstV}
\end{equation}
of $\Maps^{(s,1/p),p}(I,\T)$ with $\tilde V$ given by (I) above.  Thus, the estimate (\ref{sec2:NLinterp3}) and the trivial estimate that $\w{F_{\co}}^{-1}$ preserves $L^p(I,B^{s+1/p,p}(Y))$ regularity on (\ref{eq:mapsU}) implies that $\w{F_{\co}}^{-1}$ preserves $B^{(s,1/p),p}(I \times Y)$ regularity when restricted to
\begin{equation}
  \Maps^{(s,1/p),p}(I,U \cap \tilde V). \label{eq:mapsUtV}
\end{equation}

Finally, to estimate (\ref{eq:E1hatQ}), we estimate $\w{Q_\co}$ by applying function space multiplication and Corollary \ref{CorSW}.  Since
$$B^{(s,1/p),p}(I\times\Sigma) \hookrightarrow B^{s,p}(I\times\Sigma) \cap L^\infty(I \times \Sigma)$$
by Corollary \ref{CorEmbed}, and since the latter space is an algebra by Theorem I.C.7, we have the multiplication
$$B^{(s,1/p),p}(I\times\Sigma) \times B^{(s,1/p),p}(I\times\Sigma) \to B^{s,p}(I\times\Sigma) \cap L^\infty(I \times \Sigma).$$
Since
\begin{equation}
\w{Q_\co} = \w{\left(\H_{\co}|_{X_{\co}^{s+1,p}}\right)^{-1}}\w{\Pi_{\K_{\co}^{s,p}}} \w{\q} \label{eq:swQ}
\end{equation}
is quadratic multiplication followed by the operator (\ref{hatHmap1}), then
$$\w{Q_\co}: \Maps^{(s,1/p),p}(I , \T) \times \Maps^{(s,1/p),p}(I , \T) \to \Maps^{(s,1),p}(I , \T).$$
Finally, by Theorem \ref{ThmATrace}, we have a bounded restriction map
$$\w{r_\Sigma}: \Maps^{(s,1),p}(I , \T) \to \Maps^{(s,1-1/p-\eps),p}(I , \T_\Sigma),$$
where $\eps > 0$ is arbitrary.

Altogether, this shows that
\begin{equation}
\widehat{ r_\Sigma }\widehat{ E^1_{\co} }(\widehat{ P_\co(z) } ) \in \Maps^{(s,1-1/p-\eps),p}(I , \T_\Sigma), \qquad z \in \mathcal{V}_\Sigma
\end{equation}
if we choose $\V_\Sigma$ small enough so that $\widehat{ P_\co }\V_\Sigma$ is contained in the space (\ref{eq:mapsUtV}). In particular, since $B^{(s,1-1/p-\eps),p}(I \times \Sigma) \subset B^{s,p}(I \times \Sigma)$ for $\eps$ small,
altogether, the above estimates and the formula (\ref{eq4:FSigma}) show that $\widehat{F_{\Sigma,\co}}$ preserves the $B^{s,p}(I \times \Sigma)$ topology.  Since $\w{F_{\Sigma,\co}}^{-1}$ is just given by
$$\w{F_{\Sigma,\co}}^{-1}(z) = (z_0, z_1 - \widehat{ r_\Sigma }\widehat{ E^1_{\co} }(\widehat{ P_\co } z_0)),  \qquad z \in \V_\Sigma,$$
then $\w{F_{\Sigma,\co}}^{-1}$ also preserves $B^{s,p}(I \times \Sigma)$ regularity on $\V_\Sigma$.  Since the open set $U$ in (\ref{eq:mapsU}) contains an $L^\infty(Y)$ ball, by Theorem I.3.2, then that $\V_\Sigma$ can be chosen to contain a $C^0(I,B^{s',p}(\Sigma))$ ball, since
$$\w{P_\co}: C^0(I,\T^{s',p}_\Sigma) \to C^0(I,\T^{s'+1/p,p}) \hookrightarrow C^0(I,L^\infty\T) \cap C^0(I,H^{1/p,p}\T).$$
This implies that $\w{P_\co}$ maps a small $C^0(I,B^{s',p}(\Sigma))$ ball in $\Maps^{s,p}(I,\T_\Sigma)$ into (\ref{eq:mapsUtV}).

For (iii), note that the operators
\begin{align}
  \w{P_\co}: C^0(I,\T^{s',p}_\Sigma) & \to C^0(I,\T^{s'+1/p,p}) \label{eq:uniform1} \\
  \w{F_\co}^{-1}: C^0(I,H^{1/p,p}\T) & \dashrightarrow C^0(I,H^{1/p,p}\T) \label{eq:uniform2}
\end{align}
vary continuously with $r_\Sigma\co \in \L^{s-1/p,p}$ in the $B^{s',p}(\Sigma)$ topology, i.e. the above maps vary continuously with $\co \in \M^{s,p}$ in the $B^{s'+1/p,p}(Y)$ topology.  This follows from the work in \cite{N1}, which establishes the continuous dependence of all the operators involved in the construction of these operators.  From this continuous dependence, we can now deduce that the sets $U$ and $\tilde V$ in (\ref{eq:mapsUtV}) can be constructed locally uniformly in $\co$.  Namely, the set $U$ contains a uniform $L^\infty(Y)$ ball by Theorem I.3.2, and the set $V$ contains a uniform $H^{1/p,p}(Y)$ ball by Corollary I.3.1.  Hence $\V_\Sigma$ can be chosen to contain a $\delta$-ball in the $C^0(I,B^{s',p}(\Sigma))$ topology, with $\delta$ locally uniform in $\co$.\End

With the above lemmas, we have made most of the important steps in proving our first main theorem in this section.  Namely, the above lemma constructs for us chart maps for $\Maps^{s,p}(I,\L)$ at constant paths $\g_0$ via the local straightening map $\w{F_{\Sigma,\co}}$, where $\g_0$ is identically $r_\Sigma\co$.  Now we want to construct chart maps for $\Maps^{s,p}(I,\fC(\Sigma))$ at an \textit{arbitrary} path $\g \in \Maps^{s,p}(I,\L)$.  In addition, just as we did for the chart maps at constant paths above, we want to show that the nonlinear portion of the chart map for $\g$ smooths in the $\Sigma$ directions, and moreover, we want to show the chart map contains ``large" domains, i.e., is defined on open neighborhoods of our manifolds with respect to a topology weaker than $B^{s,p}(I \times \Sigma)$.  These results will be important for our proofs of Theorems A and B in Section 4.  Namely, by having chart maps that smooth, we will be able to gain regularity in the $\Sigma$ directions for a gauge fixed connection in the proof of Theorem A.  Furthermore, as explained prior to Lemma \ref{LemmaFSigmaHat}, being able to define chart maps for $C^0(I,B^{s',p}(\Sigma))$ neighborhoods of $\Maps^{s,p}(I,\L)$ will be instrumental in the proof of compactness in Theorem B.

We have the following theorem:

\begin{Theorem}\label{ThmPathL1} (Besov Regularity for Paths through $\L$)
  Let $I$ be a bounded interval and let $s > \max(3/p,1/2+1/p)$.
  \begin{enumerate}
    \item Then $\Maps^{s,p}(I,\L)$ is a closed submanifold of $\Maps^{s,p}(I,\fC(\Sigma))$.
  \item For any $\g \in \Maps^{s,p}(I,\L)$, there exists a neighborhood $\cU$ of $0 \in T_\g\Maps^{s,p}(I,\L)$ and a smoothing map
      $$\cE^1_\g: \cU \to \Maps^{(s,1-1/p-\eps),p}(I,\T_\Sigma),$$
      where $\eps > 0$ is arbitrary, such that the map
      \begin{align*}
        \cE_\g: \cU & \to \Maps^{s,p}(I,\L)\\
        z & \mapsto \g + z + \cE^1_\g(z)
      \end{align*}
  is a diffeomorphism onto a neighborhood of $\g$ in $\Maps^{s,p}(I,\L)$.  
  \item For any $\max(2/p,1/2) < s' \leq s-1/p$, we can choose both $\cU$ and $\cE_\g(\cU)$ to contain open $C^0(I, B^{s',p}(\Sigma))$ neighborhoods of $0 \in T_\g\Maps^{s,p}(I,\L)$ and $\g \in \Maps^{s,p}(I,\L)$ respectively, i.e., there exists a $\delta > 0$, depending on $\g$, $s'$, and $p$, such that
      \begin{align*}
        \cU & \supseteq \{z \in T_\g\Maps^{s,p}(I,\L) : \|z\|_{C^0(I,B^{s',p}(\Sigma))} < \delta\}\\
        \cE^1_\g(\cU) & \supseteq \{\g' \in \Maps^{s,p}(I,\L) : \|\g'-\g\|_{C^0(I,B^{s',p}(\Sigma))} < \delta\}.
      \end{align*}
      We can choose $\delta$ uniformly for all $\g$ in a sufficiently small $C^0(I,B^{s',p}(\Sigma))$ neighborhood of any configuration in $\Maps^{s,p}(I,\L)$.
  \item Smooth paths are dense in $\Maps^{s,p}(I,\L)$.
  \end{enumerate}
\end{Theorem}

\Proof The local straightening map $\w{F_\Sigma,\co}$ in Lemma \ref{LemmaFSigmaHat} yields for us induced chart maps in the case when $\g$ is a constant path (see Definition I.B.3).  From this, to obtain a chart map centered at a general configuration $\g$, we divide up $I$ as union $I = \cup_{j=0}^n I_j$ of finitely many smaller overlapping subintervals $I_j = [a_j, b_j]$, $a_j < a_{j+1} < b_j$ for $0 \leq j \leq n - 1$, and take a constant path $\g_j$ on $I_j$ near $\g|_{I_j}$ for which to use as a chart map for $\g|_{I_j}$ in $\Maps^{s,p}(I_j,\L)$. (We will see why we want overlapping intervals in a bit.) This is possible because the chart maps at constant paths contain $C^0(I,B^{s',p}(\Sigma))$ neighborhoods, and on the small time intervals $I_j$, an arbitrary path in $\Maps^{s,p}(I_j,\L)$ can be approximated in $C^0(I_j,B^{s-1/p,p}(\Sigma)) \hookrightarrow C^0(I_j, B^{s',p}(\Sigma))$ by constant paths.  Here, it is important that these neighborhoods depend locally uniformly on the configuration in $C^0(I_j,B^{s',p}(\Sigma))$ by Lemma \ref{LemmaFSigmaHat}, and that we have the embedding $\Maps^{s,p}(I,\L) \to C^0(I,\L^{s-1/p,p})$.

With the above considerations, let $\g$ be a constant smooth path identically equal to $u_0 \in \L$.  Then the map $\cE_\g$ is the chart map associated to $\w{F_{\Sigma,\co}}$ via Lemma \ref{LemmaFSigmaHat}, where $\co \in \M$ is any configuration satisfying $r_\Sigma\co = u_0$.  Namely, define
$$\cU := \cV_\Sigma \cap T_\g\Maps^{s,p}(I,\L),$$
and then define
\begin{align}
  \cE_\g: \cU &\to \Maps^{s,p}(I,\L) \nonumber \\
  z & \mapsto \g + \widehat{F_{\Sigma,\co}}^{-1}(z). \label{cEconst}
\end{align}
The proof of Lemma \ref{LemmaFSigmaHat} shows that $\cU$ contains a $C^0(I,B^{s',p}(\Sigma))$ ball for $\max(2/p,1/2)<s'\leq s-1/p$.  Furthermore, the map
\begin{equation}
  \cE^1_\g(z) := \widehat{F_{\Sigma,\co}}^{-1}(z) - z, \label{cE1}
\end{equation}
which is therefore the nonlinear part of $\widehat{F_{\Sigma,\co}}^{-1}(z)$ has the desired mapping properties, since it is simply the map
\begin{equation}
    \cE^1_\g(z) = \widehat{ r_\Sigma }\widehat{ E^1_{\co} }(\widehat{ P_\co }(\cdot)):  \Maps^{s,p}(I,\T_\Sigma) \dashrightarrow  \Maps^{(s,1-1/p-\eps),p}(I,\T_\Sigma). \label{cE1-2}
\end{equation}
Altogether, the map (\ref{cEconst}) yields the desired chart map for $\Maps^{s,p}(I,\L)$ for the constant path $\g$, since the map $\widehat{F_{\Sigma,\co}}$ is a local straightening map for $\g$. 

We now consider a general nonconstant, nonsmooth path $\g \in \Maps^{s,p}(I,\L)$.  As we explained, Lemma \ref{LemmaFSigmaHat} implies that if we choose $a_1 > a_0$ small enough, then on $I_0 = [a_0,a_1]$, the path $\g|_{I_0} \in \Maps^{s,p}(I_0,\L)$ lies within the range of a chart map for $C^0(I,\L^{s-1/p,p})$ at the constant path identically equal to $\g(a_0)$.  By the continuous dependence of the size of the chart map in the $C^0(I,B^{s-1/p,p}(\Sigma))$ topology and since smooth configurations are dense in $\L^{s-1/p,p}$ by Theorem I.3.4, we can choose a constant smooth path $\g_0$ that remains $C^0(I,B^{s-1/p,p}(\Sigma))$ near $\g(a_0)$, and such that the chart map $\cE_{\g_0}$ for $\g_0$, constructed as above, contains $\g$ in its image.  Thus, we have
\begin{align*}
  \g|_{I_0} &= \cE_{\g_0}(z_0^*)\\
   &= \g_0 + z_0^* + \cE^1_{\g_0}(z_0^*), \qquad z_0^* \in T_{\g_0}\Maps^{s,p}(I_0,\L)
\end{align*}
for some $z_0^*$.  Here, to conclude that $z_0^* \in T_{\g_0}\Maps^{s,p}(I_0,\L)$ has the same regularity as $\g|_{I_0}$, we used that $\g_0$ is smooth so that $\g|_{I_0} - \g_0 \in \Maps^{s,p}(I_0,\T_\Sigma)$, and then we used that $\w{F_{\Sigma,\co}}^{-1}$ preserves the $B^{s,p}(I\times\Sigma)$ topology by Lemma \ref{LemmaFSigmaHat}.

Of course, we can continue the above process, whereby we have constant smooth paths $\g_j \in \Maps(I_j,\L)$, and $\g|_{I_j}$ is in the image of the chart map $\cE_{\g_j}$ for $\g_j$, i.e.,
\begin{align*}
  \g|_{I_j} &= \cE_{\g_j}(z_j^*)\\
   &= \g_j + z_j^* + \cE^1_{\g_j}(z_j^*), \qquad z_j^* \in T_{\g_j}\Maps^{s,p}(I_j,\L)
\end{align*}
for some $z_j^*$. In this way, concatenating all the local chart maps $\cE_{\g_j}$ for the constant smooth paths $\g_j$, we can define a map
\begin{align}
  \cE_\g: \oplus_{j=0}^n T_{\g_j}\Maps^{s,p}(I_j,\L) & \dashrightarrow \times_{j=0}^n \Maps^{s,p}(I_j,\fL) \nonumber \\
  (z_j)_{j=0}^n & \mapsto \Big( t \mapsto \cE_{\g_j}(z_j^*+z_j) \Big)_{j=0}^n, \quad t \in I_j. \label{chartcE0}
\end{align}
Here, the domain of $\cE_\g$ is of course restricted to the direct sum of the domains of the individual $\cE_{\g_j}$.

To get an actual chart map, we must restrict the domain of $\cE_\g$ above so that its image under $\cE_\g$ gives an honest path in $\Maps^{s,p}(I,\fL)$ when we concatenate all the local paths on the $I_j$.  Thus, define $\cU$ by
\begin{align}
  \cU &\subset \Big\{(z_j) \in \oplus_{j=0}^n T_{\g_j}\Maps^{s,p}(I_j,\L) \; : \; \cE_{\g_j }( z_j )\big|_{I_j \cap I_{j+1}} = \cE_{ \g_{j+1} }( z_{j+1} )\big|_{I_j \cap I_{j+1}}, \quad 0 \leq j \leq n - 1\Big\},
\end{align}
where $\cU$ is any sufficiently small open subset containing $0$ on which $\cE_\g$ is defined.  In this way, the map (\ref{chartcE0}) induces a well-defined map
\begin{align}
  \cE_\g: \cU & \to \Maps^{s,p}(I_j,\L) \nonumber \\
  (z_j)_{j=0}^n & \mapsto \Big( t \mapsto \cE_{\g_j}(z_j^*+z_j) \Big), \quad t \in I_j. \label{chartcE}
\end{align}
which maps $0$ to $\g$ and the open set $\cU$ diffeomorphically onto a neighborhood of $\g \in \Maps^{s,p}(I,\L)$.  Moreover, we can choose $\cU$ so that it contains a $C^0(I,B^{s',p}(\Sigma)$ open ball, in which case, $\cE_\g(\cU)$ contains a $C^0(I,B^{s',p}(\Sigma))$ neighborhood of $\g$ in $\Maps^{s,p}(I,\fL)$.

Note that in defining $\cE_\g$ as above, while the tangent space $T_\g\Maps^{s,p}(I,\L)$ naturally sits inside $T_\g\Maps^{s,p}(I,\L)$ as the subspace\footnote{Note that the tangent space $T_\g\Maps^{s,p}(I,\L)$ only makes sense since we can prove that $\Maps^{s,p}(I,\L)$ is a submanifold of $\Maps^{s,p}(I,\fC(\Sigma))$.  Otherwise, (\ref{TMapsL0}) would just be a formal equality instead of an actual equality.}
\begin{align}
  T_\g\Maps^{s,p}(I,\L) = \{z \in \Maps^{s,p}(I,\T_\Sigma): z(t) \in T_{\g(t)}\L^{s-1/p,p}, \quad\textrm{ for all } t \in I\}, \label{TMapsL0}
\end{align}
we really study this space under the identification
\begin{align}
  T_\g\Maps^{s,p}(I,\L) \cong & \Big\{(z_j)_{j=0}^n \in \Big(\oplus_{j=0}^n T_{\g_j}\Maps^{s,p}(I_j,\L)\Big) : \nonumber \\
   & \qquad (D_{z_j^*}\cE_{\g_j})( z_j ) \big|_{I_j \cap I_{j+1} } = (D_{z_{j+1}^*}\cE_{\g_{j+1}}) ( z_{j+1} ) \big|_{I_j \cap I_{j+1} }, \quad 0 \leq j \leq n-1 \Big\}. \label{TMapsL}
\end{align}
The space (\ref{TMapsL0}) is difficult to study because it is described as a family of varying subspaces, which makes it hard to understand, for example, the mapping properties of projections onto this space.  Indeed, any resulting projection would be time-dependent, but not in a smooth way, since $\g \in \Maps^{s,p}(I,\L)$ is in general not smooth.  Consequently, it would be hard to prove an analog of Lemma \ref{LemmaSW} for such a nonsmooth time-dependent slicewise operator.  On the other hand, with (\ref{TMapsL}), we understand each factor $T_{\g_j}\Maps^{s,p}(I,\L)$ from Lemma \ref{LemmaFSigmaHat}.  These spaces do have bounded projections onto them, given by a slicewise Calderon projection as in (\ref{eq:caldproj3}), which we understand because it is time independent.  Thus, in (\ref{TMapsL}), we have simplified matters by constructing ``local trivializations" of the space $T_\g\Maps^{s,p}(I,\L)$.  Indeed, (\ref{TMapsL}) tells us that on a small time interval $I_j$, we can identify $T_{\g|_{I_j}}\Maps^{s,p}(I_j,\L)$ with $T_{\g_j}\Maps^{s,p}(I_j,\L)$, and the total space $T_\g\Maps^{s,p}(I,\L)$ is obtained by gluing together these local spaces. (This is why we chose the intervals $I_j$ in $I = \cup_{j=0}^n I_j$ to overlap.) When performing estimates for chart maps, it is thus convenient to work with the identification (\ref{TMapsL}), whereas when we wish to regard $T_\g\Maps^{s,p}(I,\L)$ as a natural subspace of $\Maps^{s,p}(I,\fC(\Sigma))$, then we have the equality (\ref{TMapsL0}).

Thus, with the identification (\ref{TMapsL}), the map $\cE_\g^1$ is just the concatenation of all the local $\cE_{\g_j}^1$ and thus has the requisite smoothing properties.  One could also compute $\cE_\g^1$ in the case where we regard $T_\g\Maps^{s,p}(I,\L)$ as (\ref{TMapsL0}), but this will not be necessary.  Altogether, we have proven (i) and (ii).  Statement (iii) now follows directly from the corresponding property for chart maps at constant paths.

For (iv), the last statement follows from the fact that smooth configurations are dense in $T_{\g_j}\Maps^{s,p}(I_j,\L)$.  Indeed, since $\g_j \in \Maps(I,\L)$ is constant and smooth, it has a lift to a constant smooth path in $\Maps(I,\M)$, which is identically $\cx[j] \in \M$ for some $\cx[j]$.  This follows from Theorem I.4.3.  Since $\cx[j]$ is smooth, the time-independent slicewise Calderon projection $\w{P^+_{\cx[j]}}$ gives a projection of the space of smooth paths $\Maps(I_j,\T_\Sigma)$ onto $T_{\g_j}\Maps(I_j,\L)$, see Theorem I.2.1.  If we mollify $z_j^* \in T_{\g_j}\Maps^{s,p}(I_j,\L)$ on $I_j \times \Sigma$, we get elements $z_j^{\eps} \in \Maps(I_j,\T_\Sigma)$ such that $z_j^{\eps} \to z_j^*$ in $\Maps^{s,p}(I_j,\T_\Sigma)$ as $\eps \to 0$.  Hence, we now obtain smooth elements $\w{P^+_\co}(z_j^\eps)$ belonging to $T_{\g_j}\Maps(I_j,\L)$ that converge to $z_j^*$ in $T_{\g_j}\Maps^{s,p}(I_j,\L)$ as $\eps \to 0$.  We then have that the $\cE_{\g_j}(z_j^\eps)$ are smooth paths which approach $\g|_{I_j}$ as $\eps \to 0$ (indeed, since $\g_j$ is smooth, one can see from the proof of Lemma \ref{LemmaFSigmaHat} that $\w{F_{\Sigma,\g_j(a_j)}}$ preserves $C^\infty$ smoothness, and hence so does $\cE_{\g_j}$).  Gluing together all these paths on the $I_j$ yields a smooth path in $\Maps^{s,p}(I,L)$ approximating $\g_j$ in the $B^{s,p}(I \times \Sigma)$ topology.\End

Because of the regularity preservation property in Lemma \ref{LemmaFSigmaHat}(ii) and because we have the mixed regularity estimates in Corollary \ref{CorSW}(ii), whereby if the base configuration $\co \in \M^{s_1+s_2,2}$ is more regular we obtain the additional mapping property (\ref{hatHmap2}), we also get the following theorem for mixed regularity paths through $\L$.  As explained earlier, the reason we consider these mixed topologies is because they will arise in the bootstrapping procedure occurring in the proof of Theorem A.  We only need the below result for small $s_2$, but we state it for general $s_2 \geq 0$.

\begin{Theorem}\label{ThmPathL2} (Mixed Regularity Paths through $\L$)
Assume the hypotheses of Theorem \ref{ThmPathL1}.  In addition, assume $s_1 \geq 3/2$ and $s_2 \geq 0$.
  \begin{enumerate}
    \item Then $\Maps^{s,p}(I,\L) \cap \Maps^{(s_1,s_2),2}(I,\L)$ is a closed submanifold of the space\\
          $\Maps^{s,p}(I,\fC(\Sigma)) \cap \Maps^{(s_1,s_2),2}(I,\fC(\Sigma))$.
  \item For any $\g \in \Maps^{s,p}(I,\L) \cap \Maps^{(s_1,s_2),2}(I,\L)$, there exists a neighborhood $\cU$ of\\
      $0 \in T_\g\Big(\Maps^{s,p}(I,\L) \cap \Maps^{(s_1,s_2),2}(I,\L)\Big)$ and a smoothing map
      $$\cE^1_\g: \cU \to \Maps^{(s,1-1/p-\eps),p}(I,\T_\Sigma) \cap \Maps^{(s_1, s_2 + 1 - \eps'),2}(I,\T_\Sigma),$$
      where $\eps,\eps' > 0$ are arbitrary, such that the map
      \begin{align*}
        \cE_\g: \cU & \to \BMaps^{s,p}(I,\L) \cap \HMaps^{(s_1,s_2),2}(I,\L)\\
        z & \mapsto \g + z + \cE^1_\g(z)
      \end{align*}
      is a diffeomorphism onto a neighborhood of $\g$ in $\BMaps^{s,p}(I,\L) \cap \HMaps^{(s_1,s_2),2}(I,\L)$.  If $s_1 > 3/2$ or $s_2 > 0$, we can take $\eps' = 0$ above.    
  \item For any $\max(s_2,2/p,1/2) < s' \leq s-1/p$, we can choose both $\cU$ and $\cE_\g(\cU)$ to contain open $C^0(I, B^{s',p}(\Sigma))$ neighborhoods of $0 \in T_\g(\Maps^{s,p}(I,\L)\cap \Maps^{(s_1,s_2),2}(I,\L))$ and $\g \in \Maps^{s,p}(I,\L)\cap \HMaps^{(s_1,s_2),2}(I,\L)$ respectively, i.e., there exists a $\delta > 0$, depending on $\g$, $s'$, $s_2$, and $p$ such that
      \begin{align*}
        \cU & \supseteq \{z \in T_\g(\Maps^{s,p}(I,\L) \cap \HMaps^{(s_1,s_2),2}(I,\L)) : \|z\|_{C^0(I,B^{s',p}(\Sigma))} < \delta\}\\
        \cE^1_\g(\cU) & \supseteq \{\g' \in \Maps^{s,p}(I,\L) \cap \HMaps^{(s_1,s_2),2}(I,\L) : \|\g'-\g\|_{C^0(I,B^{s',p}(\Sigma))} < \delta\},
      \end{align*}
      We can choose $\delta$ uniformly for all $\g$ in a sufficiently small $C^0(I,B^{s',p}(\Sigma))$ neighborhood of any configuration in $\Maps^{s,p}(I,\L) \cap \Maps^{(s_1,s_2),2}(I,\L)$.
  \item Smooth paths are dense in $\Maps^{s,p}(I,\L) \cap \Maps^{(s_1,s_2),2}(I,\L)$.
  \end{enumerate}
\end{Theorem}

\Proof The proof is exactly as the same as in Theorem \ref{ThmPathL1} and the steps made in Lemma \ref{LemmaFSigmaHat}, only we have to check that the relevant operators have the right mapping properties when we take into account the new topology we have introduced.  First, generalizing Lemma \ref{LemmaFSigmaHat}, we show that given $\co \in \M^{s,p} \cap \M^{s_1+s_2,2}$, we obtain a
local straightening map $\w{F_{\Sigma,\co}}$ for $C^0(I,\L^{s-1/p,p})$ that preserves the $B^{s,p}(I \times \Sigma) \cap B^{(s_1,s_2),2}(I \times \Sigma)$ topology.  The same proof of Lemma \ref{LemmaFSigmaHat}(ii), redone taking into account the $B^{(s_1,s_2),2}(I \times \Sigma)$ topology, shows that this is indeed the case.  Here, the dependence of $\delta$ on $s_2$ reflects the fact that we now have to interpolate the estimate
\begin{equation}
  \w{F_\co}^{-1}: H^{k,2}(I,H^{s_2+1/2,2}(Y)) \dashrightarrow H^{k,2}(I,H^{s_2+1/2,2}(Y)), \qquad k \geq 0, \label{sec2:NLinterp2}
\end{equation}
in addition to (\ref{sec2:NLinterp2}), since we want to show that $\w{F_\co}^{-1}$ preserves $\Maps^{(s_1,s_2+1/2),2}(I,  Y)$ regularity.  (To minimize notation, in the above and in the rest of this proof, we identify configuration spaces with their function space topologies). The analogous set $\tilde V$ we obtain in the proof of Lemma \ref{LemmaFSigmaHat} is thus a subset of $H^{1/p,p}\T \cap H^{s_2+1/2,2}\T$, and hence depends on $s_2$.  Thus, if we define the set
\begin{equation}
  \{z \in \Maps^{(s,1/p),p}(I,\T) \cap \Maps^{(s_1,s_2+1/2),2}(I,\T) : z(t) \in \tilde V, \textrm{ for all }t \in I\} \label{eq.2.2.domV}
\end{equation}
analogous to (\ref{eq:mapstV}), then on the domain (\ref{eq.2.2.domV}), the map $\w{F_\co}^{-1}$ preserves $\Maps^{(s_1,s_2+1/2),2}(I , Y)$ regularity.  Using Corollary \ref{CorSW}, the analogous interpolation argument as before shows that $\w{F_{\Sigma,\co}}$ and $\w{F_{\Sigma,\co}}^{-1}$ preserve $\Maps^{(s_1,s_2),2}(I , \Sigma)$ regularity.

Next, we want $s' > s_2$, since this implies $C^0(I,B^{s',p}(\Sigma)) \hookrightarrow C^0(I,H^{s_2,2}(\Sigma))$ and so that
$$\w{P_\co}: C^0(I,B^{s',p}(\Sigma)) \to C^0(I,H^{s_2+1/2,2}(Y))$$
is bounded. This implies that $\w{P_\co}$ maps a $C^0(I,B^{s',p}(\Sigma))$ small ball into $C^0(I, \tilde V)$, and hence into a domain on which $\w{F_\co}^{-1}$ is well-defined and preserves regularity. Moreover, we also get the requisite local uniformity of $\delta$ with respect to $\g$, by doing the analogous continuous dependence analysis of Lemma \ref{LemmaFSigmaHat}(iii).

To establish the mapping property of $\cE^1_\g$, it remains to show that
\begin{equation}
  \widehat{ r_\Sigma }\w{ E^1_{\co} }(\w{ P_\co }(\cdot) ) : B^{s,p}(I \times \Sigma) \cap B^{(s_1,s_2),2}(I \times \Sigma) \to B^{(s_1,s_2+1 - \eps'),2}(I \times \Sigma). \label{eq:E1-A}
\end{equation}
for $\co$ smooth.  (Here, $\co$ is a smooth configuration that is nearby $\g$ on a small interval, which we may take to be $I$, as in the analysis of the previous theorem.) By the exact same argument as in Theorem \ref{ThmPathL1}, we have bounded maps.
\begin{align*}
  \w{P_\co}: B^{s,p}(I \times \Sigma) & \to B^{(s,1/p),p}(I \times Y),\\
  \w{P_\co}: B^{(s_1,s_2),2}(I \times \Sigma) & \to B^{(s_1,s_2+1/2),2}(I \times Y),
\end{align*}
for $\co \in \M$.  Next, we show that
\begin{equation}
  \w{Q_\co}: B^{(s,1/p),p}(I \times Y) \cap B^{(s_1,s_2+1/2),2}(I \times Y) \to B^{(s_1,s_2+3/2-\eps'),2}(I \times Y). \label{hatQ-A}
\end{equation}
Thus, we need to estimate the operators appearing in equation (\ref{eq:swQ}).  First, using the embedding $B^{(s,1/p),p}(I \times Y) \hookrightarrow L^\infty(I\times Y)$ by Corollary \ref{CorEmbed}, we have a quadratic multiplication map
$$\w{\q}: \Big(B^{(s_1,s_2+1/2),2}(I \times Y) \cap L^\infty(I\times Y) \Big)^2  \to B^{(s_1,s_2+1/2-\eps'),2}(I \times Y)$$
where $\eps' > 0$ is arbitrary (we can take $\eps' = 0$ if $s_1 > 3/2$ or $s_2 > 0$).  This follows from Theorem \ref{ThmMult}.
From this, Corollary \ref{CorSW}(ii) implies (\ref{hatQ-A}).  Finally, we apply Theorem \ref{ThmATrace}, which gives us a bounded restriction map
$$\w{r_\Sigma}: B^{(s_1,s_2+3/2-\eps'),2}(I \times Y) \to B^{(s_1,s_2+1-\eps'),2}(I \times \Sigma).$$
Altogether, this completes the proof of (\ref{eq:E1-A}).  The proof of the theorem now follows as in Theorem \ref{ThmPathL1}.\End

\begin{Rem}\label{RemStraight}${}\;$
  \begin{enumerate}
    \item From now on, for any $\g \in \Maps^{s,p}(I,\L)$, we will not be concerned with precisely which model of $T_\g\Maps^{s,p}$ we need, i.e., the subspace model (\ref{TMapsL0}) or the locally ``straightened" model (\ref{TMapsL}), since both are equivalent.  All that matters is that we have chart maps $\cE_\g$ and $\cE^1_\g$ as in Theorems \ref{ThmPathL1} and \ref{ThmPathL2} which obey the analytic properties stated.  This will be the case for the proofs of Theorems A and B in Section 4.
    \item Since we have just shown that $\Maps^{s,p}(I,\L)$ is indeed a manifold, for $s > \max(3/p,1/2+1/p)$, then the family of spaces $T_{\g(t)}\L^{s-1/p,p}$, $t \in I$, does indeed comprise the tangent space $T_\g\Maps^{s,p}(I,\L)$ via (\ref{TMapsL0}).  By the density of smooth configurations, we have $\Maps^{s,p}(I,\L)$ is the $B^{s,p}(I \times \Sigma)$ closure of the space of smooth paths $\Maps(I,\L) = \{z \in C^\infty(I \times \Sigma): z(t) \in \L\}$ through the Lagrangian.  Thus, (\ref{TMapsL0}) is the same space as
  \begin{equation}
    B^{s,p}(I \times \Sigma) \textrm{ closure of } \{z \in \Maps(I,\T_\Sigma): z(t) \in T_{\g(t)}\L, \quad\textrm{ for all } t \in I\}. \label{formal-lin}
  \end{equation}
  In general, if we replace the submanifold $\L \subset \fC(\Sigma)$, with another submanifold $\fL \subset \fC(\Sigma)$, which we suppose, like $\L$, is a Fr\'echet submanifold of the Frech\'et affine space $\fC(\Sigma)$, then in general, it may not be the case that (\ref{formal-lin}) with $\L$ replaced by $\fL$ gives the true tangent space $T_{\g}\Maps^{s,p}(I,\fL)$.  Indeed, $\Maps^{s,p}(I,\fL)$ may not even be a manifold.  Of course, if $\fL$ satisfies very reasonable properties (i.e. it is defined by local straightening maps obeying the same formal analytic properties as those of $\L$), then $\Maps^{s,p}(\fL)$ will be a manifold and $T_{\g}\Maps^{s,p}(I,\fL)$ will coincide with the appropriate modification of \ref{formal-lin}).  In other words, if we define the space (\ref{formal-lin}) to be the \textit{formal tangent space} of $\Maps^{s,p}(I,\fL)$ at $\g$, then under reasonable hypotheses on $\fL$, this space will coincide with the honest tangent space $T_{\g}\Maps^{s,p}(I,\fL)$, in the appropriate range of $s$ and $p$.

  \hspace{.53cm}In the next section, we will be considering abstract Lagrangian submanifolds $\fL \subset \fC(\Sigma)$.  All tangent spaces, therefore, will be constructed formally, in the sense above.  Of course, when we specialize to $\fL = \L$ a monopole Lagrangian, there is no distinction.
  \end{enumerate}
\end{Rem}



\section{Linear Estimates}\label{SecP3Linear}

Based on the results of the previous section, we know that the linearization of (\ref{BVP}) is well-defined.  Indeed, by Theorem \ref{ThmPathL1}, $\Maps^{1-1/p,p}(\R, \L)$ is a manifold for $p > 4$.  Thus, in this section, as a preliminary step towards proving our main theorems in the next section, we will study the linearized Seiberg-Witten equations, which for simplicity, we consider about a smooth configuration.  Since smooth configurations are dense in $\Maps^{1-1/p,p}(\R, \L)$, we will see that there is no harm in doing so.

In short, our goal is to show that the linearization of (\ref{BVP}) about a smooth configuration in a suitable gauge makes the problem elliptic.  The corresponding elliptic estimates we obtain for the linearized equations will be important in studying the nonlinear equations (\ref{BVP}) in the next section.  For now, we work abstractly and take $\fL$ to be an \textit{arbitrary} Lagrangian submanifold of $\fC(\Sigma)$.  Along the way, we will see what kinds of properties such a Lagrangian should possess in order for the linearized equations to be well-behaved, i.e., the associated linearized operator of the equations is Fredholm when acting between suitable function spaces (including anisotropic spaces) and satisfies an elliptic estimate.  At the end of this section, we show in Theorem \ref{ThmMonLag} that our monopole Lagrangians obey all such properties.  This shows that the Seiberg-Witten equations with monopole Lagrangians are well-behaved at the linear level.  Nevertheless, by working with abstract Lagrangians, not only do we isolate the essential properties of monopole Lagrangians, but we also leave room for the possibility of generalizing our results to other Lagrangians that obey suitable properties but which are not monopole Lagrangians.  (Note that in working with abstract Lagrangians, when we consider their tangent spaces, we do so formally, in the sense of Remark \ref{RemStraight}(ii).)

Altogether, our main results in this section can be roughly described as follows.  Here, we replace the time interval $\R$ for our equations with $S^1$, so that we do not have to worry about issues dealing with asymptotic behavior at infinity.  Since our main results, Theorems A and B, are of a local in time nature, there is no harm in working in a compact setting as we will see in their proofs.  Moreover, it is clear that all the results of Section 2 carry over verbatim to the periodic setting. Our first main result is Theorem \ref{ThmSWlin}, which tells us that if we consider a path $\g(t)$ through our abstract Lagrangian $\fL$, then if the family of tangent spaces $L(t) = T_{\g(t)}\fL$ satisfy the hypotheses of Definition \ref{DefStraight}, then the operator (\ref{SWlin2}) induced from the family of spaces $L(t)$ is a Fredholm operator between the appropriate spaces and obeys an elliptic estimate.  Indeed, this is the relevant operator to consider, since if $\fL$ is a monopole Lagrangian and $\g(t) = r_\Sigma(B(t),\Phi(t))$, the operator considered in Theorem \ref{ThmSWlin} is precisely the linearized operator associated to (\ref{BVP}).  (Note however, that in Theorem \ref{ThmSWlin}, we only linearize about a smooth configuration and we only consider $p = 2$ Besov spaces.)  In light of Theorem \ref{ThmMonLag}, which tells us in particular that monopole Lagrangians satisfy the hypotheses of Theorem \ref{ThmSWlin}, we have Theorem \ref{ThmFredProp}, which tells us that the operator associated to the linearized Seiberg-Witten equations with monopole Lagrangian boundary conditions is a Fredholm operator.

Our second result concerns the analog of Theorem \ref{ThmSWlin} in the anisotropic setting.  Whereas Theorem \ref{ThmSWlin} is global, in the sense that it holds on all of $S^1 \times Y$, for the anisotropic setting, we work only in a collar neighborhood of the boundary of $S^1 \times Y$, namely $S^1 \times [0,1] \times \Sigma$.  This is because the anisotropic spaces we consider will be those that have extra regularity in the $\Sigma$ directions, and so we must restrict ourselves near the boundary where there is a spltting of the underlying space into the $\Sigma$ directions and the remaining $S^1 \times [0,1]$ directions.  In fact, when we prove Theorems A and B, we will only need to worry about what happens near the boundary, since the Seiberg-Witten equations are automatically elliptic in the interior modulo gauge.  The anisotropic estimates we establish for the linearized Seiberg-Witten equations in the neighborhood $S^1 \times [0,1]\times\Sigma$ will allow us to gain regularity for the \textit{nonlinear} Seiberg-Witten equations in the $\Sigma$ direction.  In short, this is because the nonlinear part of the Lagrangian boundary condition smooths in the $\Sigma$ directions, thanks to Theorems \ref{ThmPathL1} and \ref{ThmPathL2}, and hence the nonlinearity arising from the boundary condition appears only as a lower order term. Using the linear anisotropic estimates of Theorem \ref{ThmSWlinA} and Corollary \ref{CorSubspaceFred}, this extra smoothness in the $\Sigma$ directions at the boundary allows us to gain regularity in the full neighborhood $S^1 \times [0,1] \times \Sigma$ of the boundary (again, only in the $\Sigma$ directions).  This step (which is Step Two in Theorem \ref{ThmA'}) will be key in the next section.

Having described our main results, we give a brief roadmap of this section.  In the first part of this section, we describe the appropriate gauge fixing for our linearized equations.  We end up with an operator of the form $\frac{d}{dt} + D(t)$, where $D(t)$ is a time-dependent self-adjoint operator.  The Fredholm properties of such operators are well-understood, and we want to adapt these known methods to our situation.  In the second part, using the same ansatz as before, we then generalize our results to the anisotropic setting.  Here, some nontrivial work must be done since the presence of anisotropy is a rather nonstandard situation. In particular, a key result we need to establish is that the resolvent of a certain self-adjoint operator satisfies a decay estimate  (\ref{Aresolventest}) on anisotropic spaces.\\

There are two natural choices of gauge for the equations (\ref{BVP}) that will make them elliptic.  Recall that a gauge transformation $g \in \Maps(S^1\times Y, S^1)$ acts on a configuration via
\begin{equation}
  g^*(A,\Phi) := (A - g^{-1}dg, g\Phi).
\end{equation}
The first choice of gauge is to find a nearby smooth configuration $(A_0,\Phi_0)$ and find a gauge transformation $g$ such that $g^*(A,\Phi) - (A_0,\Phi_0)$ lies in the subspace\footnote{For more details on gauge fixing, see Section 2 of \cite{N1}.}
\begin{equation}
  \K_{(A_0,\Phi_0),n} := \{(a,\phi) \in T_{(A_0,\Phi_0)}\fC(S^1 \times Y): -d^*a + i\Re(i\Phi_0,\phi) = 0, *a|_{S^1 \times \Sigma} = 0\}
\end{equation}
orthogonal to the tangent space of the gauge orbit through $(A_0,\Phi_0)$. While this is the most geometric choice, it is not the most convenient, since such gauge-fixing can only be done locally, i.e., for $(A,\Phi)$ near $(A_0,\Phi_0)$.  The second choice of gauge fixing is to pick a smooth connection $A_0$ and place $(A,\Phi)$ is the Coulomb-Neumann slice through $A_0$, i.e., pick a gauge in which $(A,\Phi)$ satisfies
\begin{equation}
  d^*(A-A_0) =0, \qquad *(A-A_0)|_{S^1\times\Sigma}=0. \label{CNgauge}
\end{equation}
For any $A \in \A(S^1 \times Y)$, one can find a gauge transformation $g \in \G_{\id}(S^1\times Y)$, unique up to a constant, such that $g^*A$ satisfies (\ref{CNgauge}).  (Indeed, if we write $g = e^\eta$, with $\eta \in \Omega^0(S^1 \times Y; i\R)$, this involves solving an inhomogeneous Neumann problem for $\eta$.)

From now on, we assume our smooth solution $(A,\Phi)$ to $SW_4(A,\Phi) = 0$ is such that $A$ satisfies (\ref{CNgauge}) with respect to some $A_0$ (to be determined later).  Picking any smooth spinor $\Phi_0$ on $S^1 \times Y$, then we have the equation
\begin{equation}
  SW_4(A,\Phi) - SW_4(A_0,\Phi_0) = -SW_4(A_0,\Phi_0), \label{SW4diff}
\end{equation}
which, via (\ref{SW4}), is a semilinear partial differential equation in $(A-A_0,\Phi-\Phi_0)$ with a quadratic linearity.

So write $A_0 = B_0(t) + \alpha_0(t)dt$ as given by (\ref{paths}), whereby a connection on $S^1 \times Y$ is expressed as a path of connections on $Y$ plus its temporal component, and define
\begin{align}
  b(t) + \xi(t)dt &= A - A_0\\
  \phi(t) &= \Phi(t)-\Phi_0(t).
\end{align}
Hence, we can express the left-hand-side of (\ref{SW4diff}) as a quadratic function of $(b(t),\phi(t),\xi(t)) \in \Maps(S^1, \Omega^1(Y;i\R) \oplus \Gamma(\S) \oplus \Omega(Y;i\R))$ depending on the reference configuration $(A_0,\Phi_0)$.  Then (\ref{SW4diff}) becomes
\begin{equation}
  \left(\frac{d}{dt} + \H_{(B_0(t),\Phi_0(t))}\right)(b,\phi) + (\rho^{-1}(\phi\phi^*)_0 - d\xi, \rho(b)\phi + \xi\phi + \xi\Phi_0 + \alpha_0\phi) = -SW_4(A_0,\Phi_0). \label{eq3:linSW}
\end{equation}
Observe that the first term arises from linearizing $\frac{d}{dt}(B,\Phi) - SW_3(B,\Phi)$ and the rest are just remaining terms, which are $-d\xi$ and a quadratic function of $(b,\phi,\xi)$.  This is now a semilinear equation in $(b,\phi,\xi)$ but it is not elliptic.  We now add in the Coulomb-Neumann gauge fixing condition to remedy this.  The condition (\ref{CNgauge}) becomes
\begin{equation}
  \frac{d}{dt}\xi - d^*b = 0, \qquad *b|_{S^1 \times \Sigma} = 0. \label{CNgauge2}
\end{equation}
If we add this equation to (\ref{eq3:linSW}), then we obtain the system of equations
\begin{align}
  \left(\frac{d}{dt} + \tH_{(B_0(t),\Phi_0(t))}\right)(b,\phi,\xi) &= -\left(\rho^{-1}(\phi\phi^*)_0, \rho(b)\phi + \xi\phi +  \xi\Phi_0 + \alpha_0\phi, 0\right) - SW_4(A_0,\Phi_0) \label{eq3:dtHeq} \\
  *b|_{S^1\times\Sigma} &= 0, \label{eq3:NBC}
\end{align}
where for any configuration $\co \in \fC(Y)$, the operator $\tH_{(B_0,\Phi_0)}$ is the augmented Hessian given by \begin{align}
  \tH_\co=
  \begin{pmatrix}
    \H_\co & -d \\ -d^* & 0
  \end{pmatrix}:\T \oplus \Omega^0(Y;i\R) & \to  \T \oplus \Omega^0(Y;i\R).  \label{tHmatrix}
\end{align}
Write
$$\tT = \T \oplus \Omega^0(Y;i\R)$$
for the domain and range of $\tH_\co$. Thus, for every $t$, the operator $\tH_{(B_0(t),\Phi_0(t))}$ augments the original Hessian $\H_{(B_0(t),\Phi_0(t))}$ by taking into account the Coulomb gauge fixing and the additional $-d\xi$ term that appears in (\ref{eq3:linSW}). To simplify the form of the equations (\ref{dtHeq}) even a bit more, we can fix a smooth reference connection $B_\rf \in \A(Y)$ and consider the time-independent augmented Hessian
\begin{equation}
  \tH_0 := \tH_{(B_\rf,0)}. \label{eq3:tH0}
\end{equation}
Then we can write (\ref{eq3:dtHeq})-(\ref{eq3:NBC}) as the system
\begin{align}
  \left(\frac{d}{dt} + \tH_0\right)(b,\phi,\xi) &= N_{(A_0,\Phi_0)}(b,\phi,\xi) - SW_4(A_0,\Phi_0) \label{dtHeq} \\
  *b|_{S^1\times\Sigma} &= 0, \label{NBC}
\end{align}
where $N_{(A_0,\Phi_0)}(b,\phi,\xi)$ is the quadratic multiplication map
\begin{equation}
  N_{(A_0,\Phi_0)}(b,\phi,\xi) = -\left(\rho^{-1}(\phi\phi^*)_0, \rho(b)\phi + \xi\phi + \xi\Phi_0 + \alpha_0\phi, 0\right) - (B_0(t)-B_\rf,\Phi_0(t))\#(b,\phi). \label{sec3:Nmap}
\end{equation}
Here and elsewhere, $\#$ denotes a bilinear pointwise multiplication map whose exact form is immaterial.

Thus, the equations (\ref{dtHeq}) and (\ref{NBC}) are altogether the Seiberg-Witten equations in Coulomb-Neumann gauge.  Observe that (\ref{dtHeq}) is a semilinear elliptic equation.  Indeed, the left-hand side is a smooth constant coefficient (chiral) Dirac operator\footnote{Indeed, one can check that $\left(\frac{d}{dt} - \tH_0\right)\left(\frac{d}{dt} + \tH_0\right)$ is a Laplace-type operator.} $\frac{d}{dt} + \tH_0$ while the right-hand side is a quadratic nonlinearity.

The boundary condition we impose on our configuration $(A,\Phi)$ (aside from the Neumann boundary condition arising from gauge-fixing) is that
\begin{equation}
  r_\Sigma(B(t),\Phi(t)) \in \fL, \qquad t \in S^1, \label{sec4:LBC}
\end{equation}
where $\frak{L} \subset \fC(\Sigma)$ is a Lagrangian submanifold.  Recall that the symplectic structure on $\fC(\Sigma)$ is given by the constant symplectic form (\ref{P3:omega}) on each tangent spaces to $\fC(\Sigma)$.  We will see shortly why the Lagrangian property is important. Altogether then, it is the linearization of the full system (\ref{dtHeq}), (\ref{NBC}), and (\ref{sec4:LBC}) that we want to study.

If we linearize the equations (\ref{dtHeq}) at a smooth configuration $(A,\Phi)$, then we obtain the linear operator
\begin{equation}
  \frac{d}{dt} + \tH_0 - D_{(A,\Phi)}N_{(A_0,\Phi_0)} \label{sec3:dtHDN}
\end{equation}acting on the space
\begin{equation}
  \Maps(S^1, \Omega^1(Y; i\R) \oplus \Gamma(\S) \oplus \Omega^0(Y;i\R)) = \Maps(S^1, \tT).
\end{equation}
To linearize the boundary condition\footnote{Unlike Section 2, here we work with smooth configurations, so if $\fL$ is a smooth manifold, then $\Maps(S^1,\fL)$ is automatically a smooth manifold and it can be linearized in the expected way.  In subsequent steps, we will be taking various closures of the tangent spaces to $\Maps(S^1,\fL)$ and in this regard, see Remark \ref{RemStraight}(ii).} (\ref{sec4:LBC}), we introduce the following setup.  Consider the full restriction map
\begin{align}
  r: \tT \to \tT_\Sigma & = \Omega^1(\Sigma; i\R) \oplus \Gamma(\S_\Sigma) \oplus \Omega^0(\Sigma;i\R)\oplus \Omega^0(\Sigma;i\R)\nonumber \\
  (b,\phi,\xi) & \mapsto (b|_\Sigma, \phi|_\Sigma, -b(\nu),\xi|_\Sigma), \label{restr}
\end{align}
mapping $\tT$ to its boundary data $\tT_\Sigma$ on $\Sigma$.  Here, $\nu$ is the unit outward normal to $\Sigma$, and so $b(\nu)$ is the normal component of $b$ at the boundary; the rest of the components of $r$ represent the tangential components of $(b,\phi,\xi)$ restricted to the boundary.  On the space of paths, the map $r$ induces a slicewise restriction map, which by abuse of notation we again denote by $r$ (instead of $\w{r}$):
\begin{equation}
  r: \Maps(S^1,\tT) \to \Maps(S^1,\tT_\Sigma).
\end{equation}
Thus, specifying boundary conditions on the space $\Maps(S^1, \tT)$ for the linearized operator (\ref{sec3:dtHDN}) is equivalent to specifying a subspace of $\Maps(S^1,\tT_\Sigma)$ which determines the admissible boundary values.

The linearization of $\Maps(S^1,\fL)$ along a path $\g$ through $\fL$ along with the Neumann boundary condition (\ref{NBC}) will determine for us a subspace of $\Maps(S^1,\tT_\Sigma)$.  In order to express this, let $L \subset \T_\Sigma$ be any subspace and define the augmented space
\begin{align}
\wt{L} & := L \oplus 0 \oplus \Omega^0(\Sigma;i\R)\\
 &\subseteq \T_\Sigma \oplus \Omega^0(\Sigma; i\R) \oplus \Omega^0(\Sigma; i\R) = \tT_\Sigma. \nonumber
\end{align}
This subspace determines the subspace
\begin{equation}
  \tT_L = \{(b,\phi,\xi) \in \tT : r(b,\phi,\xi) \in \wt{L}\}
\end{equation}
of $\tT$ whose boundary values lie in $\wt{L} \subset \tT_\Sigma$. Given a family of subspaces $L(t) \subset \T_\Sigma$, $t \in S^1$, we thus get a corresponding family of spaces $\wt{L}(t) \subset \tT_\Sigma$ and $\tT_{L(t)} \subset \tT$.  The spaces $\tT_{L(t)}$ can be regarded as a family of domains for the operator $\tH_0$, and thus, we can regard the space
\begin{align}
  \Maps(S^1, \tT_{L(t)}) & := \{(b(t),\phi(t),\xi(t)) \in \Maps(S^1, \tT) : r(b(t),\phi(t),\xi(t)) \in \tT_{L(t)}, \textrm{ for all }t \in S^1\}\\
  &  \subset \Maps(S^1,\tT) \nonumber
\end{align}
as a domain for (\ref{sec3:dtHDN}).

Altogether then, using this setup, the linearization of the Seiberg-Witten equations with boundary conditions (\ref{sec4:LBC}) at a smooth configuration $(A,\Phi)$, where $A$ is in Coulomb-Neumann gauge (\ref{CNgauge}), yields the operator
\begin{align}
  \frac{d}{dt} + \tH_0 - D_{(A,\Phi)}N_{(A_0,\Phi_0)}: \Maps(S^1,\tT_{L(t)}) & \to \Maps(S^1,\tT) \label{SWlin0}\\
  L(t) &= T_{r_\Sigma(B(t),\Phi(t))}\fL, \qquad t\in S^1. \label{SWlin1}
\end{align}
Indeed (\ref{SWlin1}) is precisely the linearization of (\ref{sec4:LBC}), and this linearized boundary condition along with the Neumann boundary condition (\ref{NBC}), is precisely what defines the domain $\tT_{L(t)}$.  We want to obtain estimates for the operator (\ref{SWlin0}) on the appropriate function space completions.  For this, it suffices to consider the constant-coefficient operator
\begin{equation}
    \frac{d}{dt} + \tH_0: \Maps(S^1,\tT_{L(t)}) \to \Maps(S^1,\tT), \label{SWlin2}
\end{equation}
since (\ref{SWlin2}) differs from (\ref{SWlin0}) by a smooth multiplication operator for $(A,\Phi)$ and $(A_0,\Phi_0)$ smooth.

\begin{align}
\begin{split}
    \xy
(0,0)*+{\tT_\Sigma}="1";
(35,0)*+{\tT}="2";
(0,-20)*+{\tiL}="3";
(35,-20)*+{\tT_L}="4";
(55,-10)*+{\tT}="4.5";
{\ar_{r} "2";"1"};
{\ar^{r} "4";"3"};
{\ar^{\tH_0} "2";"4.5"};
{\ar_{\tH_0} "4";"4.5"};
{\ar@{^{(}->} (0,-16)*{};"1"};    
{\ar@{^{(}->} (35,-16)*{};"2"};
(-5,-40)*+{\Maps(S^1,\tT_\Sigma)}="b1";
(40,-40)*+{\Maps(S^1,\tT)}="b2";
(-5,-60)*+{\Maps(S^1,\tiL(t))}="b3";
(40,-60)*+{\Maps(S^1,\tT_{L(t)})}="b4";
(75,-50)*+{\Maps(S^1,\tT)}="b5";
{\ar_{r} "b2";"b1"};
{\ar^{r} "b4";"b3"};
{\ar^{\frac{d}{dt}+\tH_0} "b2";"b5"};
{\ar_{\frac{d}{dt}+\tH_0} "b4";"b5"};
{\ar@{^{(}->} (-5,-56)*{};"b1"};    
{\ar@{^{(}->} (40,-56)*{};"b2"};
\endxy
\end{split} \label{MapsCD} \\ \nonumber
\end{align}

Let us now make use the requirement that the manifold $\fL$ is Lagrangian submanifold of $\fC(\Sigma)$.  What this implies is that each tangent space $L(t) = T_{\g(t)}\fL$ to $\fL$ is a Lagrangian subspace of $\T_\Sigma$.  Consequently, each augmented space $\wt{L(t)}$ is product Lagrangian in the symplectic space $\tT_\Sigma$, where the symplectic form on $\tT_\Sigma$ is given by the product symplectic form
\begin{align}
  \tilde\omega: \tT_\Sigma \oplus \tT_\Sigma & \to \R \nonumber\\
    \tilde\omega((a,\phi,\alpha_1,\alpha_0),(b,\psi,\beta_1,\beta_0)) &= \omega((a,\phi),(b,\psi)) + \int_\Sigma (-\alpha_0\beta_1+\alpha_1\beta_0). \label{eq3:tomega}
\end{align}
Recall from \cite{N1} that the symplectic forms $\omega$ and $\tilde\omega$ naturally arise from Green's formula (I.2.65) for the Hessian and augmented Hessian operators $\H_\c$ and $\tH_\c$, respectively, for any $\c \in \fC(Y)$.  In this way, each Lagrangian subspace $L(t) \subset \T_\Sigma$ yields for us a domain $\tT_{L(t)}$ on which $\tH_0$ is symmetric, since for all $x,y \in \tT_{L(t)}$, we have
\begin{equation}
  \tomega(r(x),r(y)) = -(\tH_0 x, y)_{L^2(Y)} + (x,\tH_0 y)_{L^2(Y)} = 0,
\end{equation}
since $r(x),r(y) \in \wt{L(t)}$.

Thus, when we impose Lagrangian boundary conditions for the Seiberg-Witten equations, the associated linear operator (\ref{SWlin2}) is of the form $\frac{d}{dt} + D(t)$, where $D(t)$ is a formally self-adjoint operator with time-varying domain.  In this abstract situation, we can say the following.  When the domain of $D(t)$ is constant and furthermore, $D(t)$ is a self-adjoint, Fredholm operator (with respect to the appropriate topologies), then there is a vast literature concerning the corresponding operator $\frac{d}{dt} + D(t)$, since the Fredholm and spectral flow properties of such operators, for example, constitute a rich subject.  When the domains of $D(t)$ are varying, the results of the constant domain case can be carried over as long as the domains of $D(t)$ satisfy appropriate ``trivialization" conditions (e.g., see \cite[Appendix A]{SaWe}).  We will study (\ref{SWlin2}) from this point of view.  There are thus two things we wish to impose on our Lagrangian $\fL$ so that its tangent spaces $L(t) = T_{\g(t)}\fL$ all obey the following loosely formulated conditions:
\begin{description}
  \item{(I)} for each domain $\tT_{L(t)}$, the operator $\tH_0: \tT_{L(t)} \to \tT$ is self-adjoint (as opposed to formally self-adjoint) and Fredholm with respect to the appropriate topologies;
  \item{(II)} the time-varying domains $\tT_{L(t)}$ satisfy the appropriate trivialization conditions.
\end{description}

To make (I) and (II) more precise, we explain the function spaces we will be using.  Namely, we will work with the Besov spaces $B^{s,2}$ with exponent $p = 2$.  Recall these spaces are also denoted by $H^s$, though since we have been working primarily with Besov spaces in this paper, we will stick to the notation $B^{s,2}$ to be consistent.  We want to work with $p = 2$, because on $L^2$ spaces, one can employ Hilbert space methods, in particular, one has the spectral theorem and unitarity of the Fourier transform.  The $p \neq 2$ analysis developed in the previous section will come into play for the nonlinear analysis of Seiberg-Witten equations, which we take up in the next section.     We may thus consider the operator
\begin{equation}
    \frac{d}{dt} + \tH_0: \Maps^{k,2}(S^1,\tT_{L(t)}) \to \Maps^{k-1,2}(S^1,\tT) \label{P3:SWlin} \\
\end{equation}
for all integers $k \geq 1$.  The spaces $\Maps^{k,2}(S^1,\tT)$ is defined as in Definition \ref{DefAMaps}.  To define the space $\Maps^{k,2}(S^1,\tT_{L(t)})$ with varying domain $\tT_{L(t)}$, one proceeds in a similar way for $k \geq 2$. In this case, one can take a trace twice for functions in $B^{k,2}(S^1 \times Y)$, and so
$\Maps^{k,2}(S^1,\tT_{L(t)})$ is the subspace of $\Maps^{k,2}(S^1,\tT)$ whose paths $(b(t),\phi(t),\xi(t))$ satisfy
$$r(b(t),\phi(t),\xi(t)) \in B^{k-1,2}\wt{L(t)} \subset \tT^{k-1,2}_\Sigma, \qquad t \in S^1.$$
Unfortunately, this definition does not work for $k=1$.  Thus, we have the following definition which works for all $k \geq 1$ and which coincides with the above definition.  Namely, we define
\begin{equation}
  \Maps^{k,2}(S^1,\tT_{L(t)}) := B^{k,2}(S^1, \tT^{0,2}) \cap L^2(S^1, \tT^k_{L(t)}). \label{sepvariables}
\end{equation}
In other words, $\Maps^{k,2}(S^1,\tT_{L(t)})$ consists of those paths that have $k$ time derivatives in the space $L^2(S^1, \tT^{0,2})$ and which belong to $\tT^k_{L(t)}$ in the $L^2(S^1)$ sense, where
$$\tT^{k,2}_{L(t)} = \{(b,\psi,\xi) \in \tT^{k,2}: r(b,\psi,\xi) \in B^{k-1/2,2}\wt{L(t)} \subset \tT^{k-1/2,2}_\Sigma\}$$
makes sense for $k\geq 1$.  Thus, we have merely separated variables in  (\ref{sepvariables}) and ask that all derivatives of order $k$ exist in $L^2$ and that the appropriate boundary conditions hold.

With the above definitions, we can thus consider the operator (\ref{P3:SWlin}) and the family of operators
\begin{equation}
  \tH_0: \tT^{k+1,2}_{L(t)} \to \tT^{k,2}, \qquad t\in S^1. \label{tH0-map}
\end{equation}
It is these operators and their domains that have to satisfy the appropriate assumptions, by the above discussion, in order for us to obtain suitable estimates for (\ref{P3:SWlin}).  According to the first condition (I) above, we want that the family of operators (\ref{tH0-map}) to be Fredholm for all $k \geq 0$, and furthermore, for $k =0$, that they are all self-adjoint.  To obtain (II), we make the following definition:

\begin{Def}\label{DefStraight}
  Let $L(t)$ be a smoothly varying family of subspaces\footnote{See Definition I.A.3.} of $\T_\Sigma$, $t \in S^1$ or $\R$.  We say that the $L(t)$ are \textit{regular} if for every $t_0 \in \R$, there exists an open interval $I \ni t_0$ such that for all $t \in I$, there exist isomorphisms $S(t): \tT \to \tT$, satisfying the following properties for all nonnegative integers $k$:
  \begin{enumerate}
    \item The map $S(t)$ extends to an isomorphism $S(t): \tT^{k,2} \to \tT^{k,2}$.
    \item The map $S(t)$ straightens the family $L(t)$ in the sense that $S(t): \tT_{L(t_0)}^{k,2} \to \tT_{L(t)}^{k,2}$ is an isomorphism ($k \geq 1$).
    \item The commutator $[D,S(t)]$, where $D: \tT \to \tT$ is any first order differential operator, is an operator bounded on $\tT^{k,2}$.
  \end{enumerate}
\end{Def}

The reason for this definition is that then, on the interval $I$ in the above, where say $t_0 = 0$, the conjugate operator
\begin{equation}
  S(t)^{-1}\left(\frac{d}{dt} + \tH_0\right)S(t): \Maps^{k,2}(I, \tT_{L(0)}) \to \Maps^{k-1,2}(I,\tT) \label{conj-op}
\end{equation}
has constant domain, for all $k \geq 1$.  Conditions (i) and (ii) ensure that (\ref{conj-op}) is well-defined.  Condition (iii) ensures that the conjugate operator (\ref{conj-op}) gives us a lower order perturbation of the original operator, since
\begin{equation}
  S(t)^{-1}\left(\frac{d}{dt} + \tH_0\right)S(t) = \frac{d}{dt} + \tH_0 + \left(S(t)^{-1}\frac{d}{dt}S(t) + S^{-1}(t)[\tH_0,S(t)]\right).
\end{equation}
Thus, (iii) and the fact that $S(t)$ depends smoothly on $t$ implies that the right-most term of the above is a bounded operator.  One can now understand the time-varying domain case in terms of the constant domain case via this conjugation (see Theorem \ref{ThmSWlin}). Thus, the map $S(t)$ trivializes, or ``straightens", the family of subspaces $\tT_{L(t)}$, and it is the existence of such an $S(t)$ that expresses precisely what we mean by condition (II) above.

We will show later that when $L(t)$ are the subspaces arising from linearizing a monopole Lagrangian $\fL$ along a smooth path, then the operators (\ref{tH0-map}) are Fredholm for all $k$, self-adjoint for $k=0$, and the $L(t)$ are regular (see Theorem \ref{ThmMonLag}).  This uses the fundamental analysis concerning monopole Lagrangians in \cite{N1}.  Assuming these properties hold for some general Lagrangian submanifold $\fL \subset \fC(\Sigma)$, we have our first result concerning the linearization of the Seiberg-Witten equations with (general) Lagrangian boundary conditions:

\begin{Theorem}\label{ThmSWlin}
\begin{enumerate}
  \item Suppose the operators (\ref{tH0-map}) are Fredholm for all $k \geq 0$, and that furthermore, they are self-adjoint for $k=0$.  Also, suppose the family of spaces $L(t) = T_{\g(t)}\fL$, $t \in S^1$, are regular. Then
\begin{equation}
  \frac{d}{dt} + \tH_0 : \Maps^{k+1,2}(S^1, \tT_{L(t)}) \to \Maps^{k,2}(S^1,\tT) \label{dtHk}
\end{equation}
is Fredholm for every $k \geq 0$, and we have the elliptic estimate
\begin{equation}
  \|(b,\phi,\xi)\|_{B^{k+1,2}(S^1 \times Y)}^2 \leq C\left(\left\|\left(\frac{d}{dt}+ \tH_0\right)(b,\phi,\xi)\right\|_{B^{k,2}(S^1\times Y)}^2 + \|(b,\phi,\xi)\|_{B^{k,2}(S^1\times Y)}^2\right). \label{est_dtHk}
\end{equation}
\item (Elliptic regularity) If $(b,\phi,\xi) \in \Maps^{1,2}( S^1, \tT_{L(t)} )$ satisfies $(\frac{d}{dt} + \tH_0)(b,\phi,\xi) \in \Maps^{k,2}(S^1,\tT)$, then $(b,\phi,\xi) \in \Maps^{k+1,2}( S^1, \tT_{L(t)} )$ and it satisfies the elliptic estimate (\ref{est_dtHk}).
\end{enumerate}
\end{Theorem}

\Proof (i) First, let $k=0$ and suppose $L(t) \equiv L(0)$ is independent of $t$.  Let us make the abbreviations $x = (b,\phi,\xi)$ and $\partial_t = \frac{d}{dt}$. There are two methods to obtain (\ref{est_dtHk}).  The first proceeds as follows.  We prove the identity
\begin{equation}
  \|(\partial_t + \tH_0)x\|_{L^2(S^1\times Y)}^2 = \|\partial_t x\|_{L^2(S^1 \times Y)}^2 + \|\tH_0 x\|_{L^2(S^1\times Y)}^2. \label{dtHIBP}
\end{equation}
via an integration by parts.  Here, we can write the cross term of (\ref{dtHIBP}) as
\begin{equation}
(\partial_tx, \tH_0x)_{L^2(Y)} + (\tH_0x,\partial_tx)_{L^2(Y)} = \int_{S^1}\frac{d}{dt}(x,\tH_0 x)_{L^2(Y)}dt, \label{intS^1}
\end{equation}
because of the self-adjointness of
\begin{equation}
  \tH_0: \tT^{1,2}_{L(t)} \to \tT^{0,2}. \label{tH1to0}
\end{equation}
and the time-independence of $L(t)$. The term (\ref{intS^1}) vanishes since we are integrating an exact form over $S^1$.  Next, since (\ref{tH1to0}) is Fredholm by hypothesis, we also have
\begin{equation}
  \|x(t)\|_{B^{1,2}(Y)}^2 \leq C(\|\tH_0 x(t)\|_{L^2(Y)}^2 + \|x(t)\|_{L^2(Y)}^2). \label{eq4:ellest12}
\end{equation}
for every $t \in S^1$.
Integrating this estimate over $S^1$, using this in (\ref{dtHIBP}), and using the fact that $\|x\|_{B^{1,2}(S^1 \times Y)}$ is equivalent to $\|\partial_tx\|_{L^2(S^1\times Y)} + \|x\|_{L^2(S^1, B^{1,2}(Y)}$, we have the elliptic estimate
\begin{equation}
  \|x\|_{B^{1,2}(S^1 \times Y)}^2 \leq C(\|(\partial_t + \tH_0)x\|_{L^2(S^1\times Y)}^2 + \|x\|_{L^2(S^1\times Y)}^2). \label{est_dtH1}
\end{equation}
This shows that the map (\ref{dtHk}) has closed range and finite dimensional kernel.  To show that the cokernel is finite dimensional, we use the fact that we have the following weak regularity estimate (for time-varying domains):
\begin{align}
  y \in \Maps^{0,2}(S^1, \tT) \textrm{ and }\left(\left(\partial_t + \tH_0\right)x, y\right) & = 0 \textrm{ for all }x \in \Maps(S^1, \tT_{L(t)})\nonumber \\ & \Rightarrow y \in \Maps^{1,2}(S^1, \tT_{L(t)}). \label{wreg}
\end{align}
This is proven in \cite[Appendix A]{SaWe}).  In light of this, an integration by parts shows that the cokernel of (\ref{dtHk}) is finite dimensional for $k = 0$ (in fact, any $k$), since the adjoint operator $-\partial_t + \tH_0$ obeys the same estimate (\ref{est_dtH1}).

There is a second approach to proving (\ref{eq4:ellest12}) which generalizes to a more general setting that we will need later.  Together with weak regularity (\ref{wreg}), this proves the Fredholm property of (\ref{dtHk}) for $k=0$.  The method we use is to apply the Fourier transform (in $t \in S^1 = [0,2\pi]/\sim$) to the time-independent operator $\partial_t + \tH_0$, which means we analyze the operators $i\tau + \tH_0$, for $\tau \in \Z$.

Without loss of generality, we can suppose the Fredholm operator (\ref{tH1to0}) is an invertible operator, which we can always do by perturbing $\tH_0$ by a bounded operator.  Indeed, the operator $\tH_0: \tT^{1,2}_{L(0)} \to \tT^{0,2}$, being a self-adjoint Fredholm operator, has discrete spectrum, and so we can perturb $\tH_0$ by some multiple of the identity to achieve invertibility.  By self-adjointness, $i\tau + \tH_0$ is invertible for all $\tau \in \R$.  Thus, if we have $\left(\partial_t + \tH_0\right)x = y$ and we want to solve for $x$, we just solve for the Fourier modes.  In other words, we have
$$\hat x(\tau) = (i\tau + \tH_0)^{-1}\hat y(\tau), \qquad \tau \in \Z,$$
where, if $z \in \HMaps^{0,2}(S^1,\tT)$, we have
$$\hat z(\tau) = \int_0^{2\pi} e^{-i\tau t}z(t) dt \in \tT^{0,2}.$$
Thus, by Plancherel's theorem,
\begin{align}
  \|\partial_tx\|_{L^2(S^1 \times Y)}^2 &= \frac{1}{2\pi}\sum_{\tau \in\Z}\|i\tau \hat x(\tau)\|_{L^2(Y)}^2 \nonumber \\
  &= \frac{1}{2\pi}\sum_{\tau \in \Z}\left\|\left(\frac{i\tau}{i\tau + \tH_0}\right)\hat y(\tau)\right\|_{L^2(Y)}^2. \label{planch}
\end{align}
From the spectral theorem, we have
\begin{equation}
  \|(i\tau - \tH_0)^{-1}\|_{Op(L^2(Y))} = O(\tau^{-1}). \label{resolvent_est}
\end{equation}
From (\ref{planch}), this implies
\begin{equation}
  \|\partial_tx\|_{L^2(S^1 \times Y)} \leq C\|y\|_{L^2(S^1 \times Y)} = C\|(\partial_t+\tH_0)x\|_{L^2(S^1 \times Y)}, \label{controldt}
\end{equation}
which implies
\begin{align}
  \|x\|_{B^{1,2}(S^1\times Y)} & \leq C(\|\partial_tx\|_{L^2(S^1 \times Y)} + \|\tH_0 x\|_{L^2(S^1 \times Y)} + \|x\|_{L^2(S^1\times Y)})\nonumber \\
  & \leq C(2\|\partial_tx\|_{L^2(S^1 \times Y)} + \|(\partial_t+\tH_0)x\|_{L^2(S^1 \times Y)} + \|x\|_{L^2(S^1\times Y)}) \nonumber \\
  & \leq C'(\|(\partial_t+\tH_0)x\|_{L^2(S^1 \times Y)} + \|x\|_{L^2(S^1\times Y)}). \label{est_dtH1''}
\end{align}
This establishes the desired elliptic estimate.  Thus, this establishes the Fredholm property of (\ref{dtHk}) by previous remarks for $k=0$.

It remains to consider the general case where $L(t)$ is time-dependent and $k \geq 0$.  To control the varying domains, we employ the straightening maps $S(t)$ associated to $L(t)$, as given by Definition \ref{DefStraight}.  We reduce to the constant domain case via conjugation by $S(t)$ as in (\ref{conj-op}).  For notational convenience, let us suppose the straightening maps can be made periodic, i.e., for all $t \in S^1$, so that we can replace $I$ with $S^1$ in (\ref{conj-op}).  (In general, we have to apply a partition of unity in time and apply local straightening maps on local intervals.  We then apply the analysis that is to follow on small time intervals and then sum up the estimates.) By Definition (\ref{DefStraight})(iii), the terms $S(t)\frac{d}{dt}S(t)^{-1}$ and $[\tHt,S(t)]S(t)^{-1}$ in (\ref{conj-op}) are bounded on $\HMaps^{1,2}(S^1,\tT_{L(0)})$.  Thus, we also have the elliptic estimate (\ref{est_dtH1}) for the conjugate operator (\ref{conj-op}) and hence the original operator on varying domains by isomorphism property, property (ii) of Definition \ref{DefStraight}.   This establishes (\ref{dtHk}) for $k=0$ and for varying domains.

For $k \geq 1$, we proceed inductively in $k$.  Suppose we have established the result for all nonnegative integers up to some $k \geq 0$ and we want to prove it for $k + 1$.  As in the previous part, we can assume we are in the time-independent case since the time-dependent case reduces to this case via conjugation by $S(t)$.  Moreover, by the weak regularity result (\ref{wreg}), we only need to establish the elliptic estimate (\ref{est_dtHk}).

So suppose $x \in \HMaps^{k+1,2}(S^1,\tT_{L(0)})$.  Then $\partial_tx \in \HMaps^{k,2}(S^1,\tT_{L(0)})$.  By the inductive hypothesis, we have the estimate (\ref{est_dtHk}) for $\partial_tx$.  Now, $\|x\|_{B^{k+1,2}(S^1\times Y)}$ is equivalent to $\|\partial_tx\|_{B^{k,2}(S^1\times Y)} + \|x\|_{L^2(S^1,B^{k+1}(Y))}$.  The inductive hypothesis gives us control of the first of these two terms from (\ref{est_dtHk}); it remains to control the second term $\|x\|_{L^2(S^1,B^{k+1,2}(Y))}$ in order to prove (\ref{est_dtHk}) for $k+1$.  Since $\tH_0: \tT^{k+1,2}_{L(0)} \to \tT^{k,2}$ is Fredholm, then
\begin{equation}
 \|x\|_{L^2(S^1,B^{k+1,2}(Y))} \leq C(\|\tH_0 x\|_{L^2(S^1,B^{k,2}(Y)} + \|x\|_{L^2(S^1 \times Y)}).
\end{equation}
Altogether, simple rearrangement yields
\begin{align*}
  \|x\|_{B^{k+1,2}(S^1\times Y)} & \sim (\|\partial_tx\|_{B^{k,2}(S^1\times Y)} + \|x\|_{L^2(S^1,B^{k+1}(Y))})\\
  & \leq C(\|\partial_tx\|_{B^{k,2}(S^1\times Y)} + \|\tH_0 x\|_{L^2(S^1,B^{k,2}(Y))} + \|x\|_{L^2(S^1 \times Y)})\\
  & \leq C(2\|\partial_tx\|_{B^{k,2}(S^1\times Y)} + \|(\partial_t + \tH_0)x\|_{L^2(S^1,B^{k,2}(Y))} + \|x\|_{B^{k,2}(S^1 \times Y)})\\
  & \leq 2C(\|(\partial_t + \tH_0)\partial_tx\|_{B^{k-1,2}(S^1\times Y)} + \|\partial_tx\|_{B^{k-1,2}(S^1\times Y)}\\ & \qquad + \|(\partial_t + \tH_0)x\|_{L^2(S^1,B^{k,2}(Y))} + \|x\|_{B^{k,2}(S^1 \times Y)})\\
  & \leq 2C(\|\partial_t(\partial_t + \tH_0)x\|_{B^{{k-1},2}(S^1 \times Y)}  + \|(\partial_t + \tH_0)x\|_{L^2(S^1,B^{k,2}(Y))}\\
  & \qquad +  \|x\|_{B^{k,2}(S^1 \times Y)}) \\
  & \sim \left\|(\partial_t + \tH_0)x\right\|_{B^{k,2}(S^1 \times Y)} +  \|x\|_{B^{k,2}(S^1 \times Y)}.
\end{align*}
In the fourth line, we applied the inductive hypothesis to $\partial_t x$.  The above computation completes the induction and establishes the elliptic estimate (\ref{est_dtHk}) for all $k$.

For (ii), we note that this easily follows from (i), since in proving (i), we implicitly constructed a (left) parametrix for $\frac{d}{dt}+\tH_0$ by the above steps.  Thus, the a priori elliptic estimate in (i) implies the elliptic regularity statement (ii).\En

\begin{Rem}\label{RemER}
  As is usual for elliptic equations, an a priori elliptic estimate implies elliptic regularity as in the above, since there is always a smoothing parametrix implicit in problem.  Henceforth, we will always prove a priori elliptic estimates and then state the corresponding elliptic regularity result without additional proof.
\end{Rem}

\subsection{Anisotropic Estimates}\label{SecLinA}

While Theorem \ref{ThmSWlin} tells us that the linearized operator associated to the Seiberg-Witten equations with suitable Lagrangian boundary conditions is Fredholm, as mentioned in the introduction to this section, we will also need an anisotropic analog.  Indeed, we did a great deal of analysis in Section \ref{SecPaths} on anisotropic Besov spaces and we will need to generalize the above theorem to such spaces.  The reason for this is that in our proof of Theorem A in the next section, we will be bootstrapping the regularity of a configuration in the $\Sigma$ directions in a neighborhood of the boundary of $\R \times Y$.  This boostrapping requires that we gain regularity in certain individual directions from the linearized Seiberg-Witten equations, which means that we want the operator (\ref{dtHk}) to be Fredholm on anisotropic spaces and to have a corresponding elliptic estimate (\ref{est_dtHk}) on anisotropic spaces.\footnote{The analysis developed in this section is completely absent in the ASD case, as seen in \cite{We1}.  Indeed, in the ASD case, there are no spinors, and the structure of the ASD equations alone allows one to easily gain $\Sigma$ regularity in a collar neighborhood of the boundary of $\R \times Y$ without even using the Lagrangian boundary condition.  As a consequence, the hard work we do in Section \ref{SecPaths} and in this section, which is to gain $\Sigma$ regularity near the boundary, is completely unnecessary and absent in the ASD case.}

In detail, on $Y$ we now work in a collar neighborhood $[0,1] \times \Sigma$ of the boundary, and consequently, on $S^1 \times Y$, we work in the collar neighborhood $S^1 \times [0,1] \times \Sigma$. In either case, we let $v \in [0,1]$ denote the inward normal coordinate.  We have the corresponding restricted configuration space
\begin{equation}
  \tT_{[0,1]\times\Sigma} := \tT|_{[0,1]\times\Sigma} = \Omega^1([0,1]\times\Sigma;i\R) \oplus \Gamma(\S|_{[0,1]\times\Sigma}) \oplus\Omega^0([0,1]\times\Sigma;i\R).
\end{equation}
The restriction map $r$ induces two separate restriction maps
\begin{equation}
  r_j: \tT_{[0,1]\times\Sigma} \to \tT_{\Sigma_j}
\end{equation}
corresponding to the two boundary components
$$\Sigma_j := \{j\} \times \Sigma, j=0,1$$
of $[0,1]\times\Sigma$. In this case, we write
\begin{equation}
  r = (r_0,r_1): \tT_{[0,1]\times\Sigma} \to \tT_{\Sigma_0} \oplus \tT_{\Sigma_1}
\end{equation}
for the total restriction map.  If we just write $\Sigma$, we will always mean $\Sigma_0$.

The space $[0,1]\times\Sigma$ is a product manifold and so we can define anisotropic Besov spaces on it.  We have the space $B^{(s_1,s_2),2}([0,1] \times \Sigma)$, the space of functions whose derivatives up to order $s_2$ in the $\Sigma$ directions belong to $B^{s_1,2}([0,1] \times \Sigma)$.  We define the spaces
$$\tT^{(k,s),2}_{[0,1]\times\Sigma} = B^{(k,s),2}\tT_{[0,1]\times\Sigma}$$
of configurations in the $B^{(k,s),2}([0,1] \times \Sigma)$ topology, where $k \geq 0$ is a nonnegative integer and $s \geq 0$.  By the anisotropic trace theorem, Theorem \ref{ThmATrace}, the restriction maps extend to bounded operators
\begin{equation}
  r_j: \tT^{(k,s),2}_{[0,1]\times\Sigma} \to \tT_{\Sigma_j}^{k-1/2+s,2}, \qquad k \geq 1, \quad j =0,1.
\end{equation}

We want to generalize Theorem \ref{ThmSWlin} to the anisotropic setting, which means we need to adapt (I) and (II) of the previous subsection to the anisotropic setting.  For (I), in order to get Fredholm operators mapping between the $\tT^{(k,s),2}_{[0,1]\times\Sigma}$ spaces, we need to impose boundary conditions as before, only now we have to impose them on the two boundary data spaces $\tT_{\Sigma_j}^{k-1/2+s,2}$.  For $\tT_{\Sigma_0}^{k-1/2+s,2}$ we impose the same type of boundary condition as before, namely by specifying a subspace $L_0$ of $\T_{\Sigma_0}$, considering the augmented space $\wt{L_0} \subseteq \tT_{\Sigma_0}$, and then taking the $B^{k-1/2+s,2}(\Sigma)$ completion.  On $\tT_{\Sigma_1}^{k-1/2+s,2}$ we also choose a subspace as a boundary condition, and for this, we choose any suitable subspace
$$L_1 \subseteq \tT_{\Sigma_1}$$
as follows.

Define the restricted configuration space
$$\tT^{(k,s),2}_{[0,1]\times\Sigma, L_0, L_1} = \{x \in \tT^{(k,s),2}_{[0,1]\times\Sigma} :\; r_0(x) \in B^{k-1/2+s,2}\wt{L_0}, \quad r_1(x) = B^{k-1/2+s,2}L_1\},$$
of configurations in $\tT^{(k,s),2}_{[0,1]\times\Sigma}$ whose boundary data on $\tT^{k-1/2+s,2}_{\Sigma_0}$ and $\tT^{k-1/2+s}_{\Sigma_1}$ lie in $\wt{L_0}$ and $L_1$, respectively.  Ultimately, we want the operator
\begin{equation}
  \tH_0: \tT^{(k+1,s),2}_{[0,1]\times\Sigma, L_0,L_1} \to \tT^{(k,s),2}_{[0,1]\times\Sigma} \label{tH0-LL}
\end{equation}
to be Fredholm, for $k \geq 0$ and $s \geq 0$.  The determining of which subspaces $L_0$ and $L_1$ determine a Fredholm operator for $\tH_0$ falls within the study of elliptic boundary value problems which we described in \cite{N1}.  From Theorem \ref{ThmEBP2}, we have the following result:  If the subspace on the boundary which determines the boundary condition for $\tH_0$, in this case the subspace
\begin{equation}
  B^{k+1/2+s,2}(\wt{L_0} \oplus L_1) \subseteq \tT^{k+1/2+s,2}_{\Sigma_0} \oplus \tT^{k+1/2+s,2}_{\Sigma_1}, \label{L0L1}
\end{equation}
is such that it is Fredholm with $r(\ker \tH_0)$ (see Definition \ref{DefFred}), then the associated operator (\ref{tH0-LL}) is Fredholm.  Indeed, when $\wt{L_0} \oplus L_1$ is the range of a pseudodifferential operator, this follows from the standard theory of pseudodifferential elliptic boundary conditions.  For us however, the space $\wt{L_0}$, being a tangent space to $\frak{L}$, a monopole Lagrangian, is only ``nearly" pseudodifferential (see Theorem I.2.1(ii)), and so we use the more general framework of Theorem \ref{ThmEBP2}.

In fact, we can say more.  Recall from \cite{N1} that $r(\ker \tH_0) = r(\ker \tH_{(B_\rf,0)})$ is given by the range of a zeroth order pseudodifferential operator, the Calderon projection
\begin{equation}
  \tP^+_0 := \tP^+_{(B_\rf,0)}: \tT_{\Sigma_0}^{k+1/2+s,2} \oplus \tT_{\Sigma_1}^{k+1/2+s,2} \circlearrowleft. \label{eq:P^+}
\end{equation}
Because the symbol of $\tP^+_0$ is determined locally by the symbol of $\tH_0$, on $\Sigma_0$, the principal symbol of $\tP^+_0$ coincides with that of $\Pi^+: \tT_{\Sigma_0}^{k+1/2+s,2} \to \tT_{\Sigma_0}^{k+1/2+s,2}$, the positive spectral projection of the tangential boundary operator (see Definition I.2.1) associated to $\tH_0$ on $\Sigma_0$.  Identifying $\tT_{\Sigma_0}$ with $\tT_{\Sigma_1}$, since the restriction map $r_1: \tT_{[0,1]\times\Sigma} \to \tT_{\Sigma_1}$ is defined as in (\ref{restr}), with $-\nu$ now the \textit{outward} normal to $\Sigma_1$ (as opposed to $-\nu$ being the inward normal at $\Sigma_0$), the principal symbol of $\tP^+_0$ on $\Sigma_1$ coincides with $\Pi^-: \tT_{\Sigma_1}^{k+1/2+s,2} \to \tT_{\Sigma_1}^{k+1/2+s,2}$, the negative spectral projection of the tangential boundary operator.  Indeed, choosing the opposite choice of normal at $\Sigma_1$ reverses the sign of the tangential boundary operator and so changes the associated positive spectral projection to a negative spectral projection.  Altogether then, the range of (\ref{eq:P^+}) is a compact perturbation of the range of
$$\im \Pi^+ \oplus \im \Pi^-: \tT_{\Sigma_0}^{k+1/2+s,2} \oplus \tT_{\Sigma_1}^{k+1/2+s,2} \circlearrowleft.$$
This is convenient because $\im \tP^+_0$ certainly not a direct sum of a subspace of $\tT_{\Sigma_0}^{k+1/2+s,2}$ with a subspace of $\tT_{\Sigma_1}^{k+1/2+s,2}$, but the above analysis tells us it is a compact perturbation of this. From this, we easily deduce
\begin{Lemma}\label{LemmaF0F1}
  Suppose we have
\begin{description}
  \item{(F0)} $B^{k+1/2+s,2}{\wt{L_0}}$ is Fredholm with $\im\Pi^+$ in $\tT^{k+1/2+s,2}_{\Sigma_0}$;
  \item{(F1)} $B^{k+1/2+s,2}{L_1}$ is Fredholm with $\im\Pi^-$ in $\tT^{k+1/2+s,2}_{\Sigma_1}$.
\end{description}
Then $B^{k+1/2+s,2}(\wt{L_0} \oplus L_1)$ Fredholm with $r(\ker \tH_0)$ in $\tT^{k+1/2+s}_{\Sigma_0}\oplus \tT^{k+1/2+s}_{\Sigma_1}$, where $\tH_0$ is the operator in (\ref{tH0-LL}).  This implies (\ref{tH0-LL}) is Fredholm.
\end{Lemma}

If $\fL$ is a monopole Lagrangian, Theorem \ref{ThmMonLag} tells us that (F0) is satisfied for $L_0$ a tangent space to $\fL$.  Thus, (F1) is the only condition that needs to be satisfied.  This latter condition is a generic open condition, and so there is great freedom in our choice of $L_1$.  Moreover, if we let $\chi: [0,1] \to \R^+$ be a smooth cutoff function, $\chi(v) = 1$ on $v \leq 1/2$ and $\chi(v) = 0$ for $v \geq 3/4$, then for any $(b,\psi,\xi) \in \tT^{(k+1,s),2}$, we have
$$\chi(b,\psi,\xi) \in \tT^{(k+1,s),2}_{[0,1]\times\Sigma,L_0,L_1}$$
for any choice of $L_1$.  Furthermore, it is these truncated configurations\footnote{More properly, the truncations of the time-dependent configurations on $S^1 \times [0,1] \times \Sigma$.} for which we will bootstrapping regularity in the proof of Theorem A.  Thus, in this sense, the choice of $L_1$ is just a ``dummy" boundary condition to make the operator (\ref{tH0-LL}) Fredholm. Thus, assuming (F0), we have described a sufficient condition (F1) that makes (\ref{tH0-LL}) Fredholm.

To generalize (II) to the anisotropic setting, we study the operator $\frac{d}{dt}+\tH_0$ with time-varying domains on anisotropic function spaces.  We are thus led to consider the space
\begin{equation}
  \Maps^{(k,s),2}(S^1, \tT_{[0,1]\times\Sigma}),
\end{equation}
the closure of the space of smooth paths $\Maps(S^1,\tT_{[0,1]\times\Sigma})$ in the topology $B^{(k,s),2}((S^1\times[0,1])\times\Sigma)$.  Thus, the $s$ measures anisotropy in the $\Sigma$ directions only.  Likewise, we can define $ \Maps^{(k,s)}(S^1, \tT_{[0,1]\times\Sigma,L(t),L_1})$ as the anisotropic analogue of (\ref{sepvariables}), namely,
\begin{equation}
  \Maps^{(k,s),2}(S^1, \tT_{[0,1]\times\Sigma,L(t),L_1}) :=  B^{k,2}(S^1, \tT^{(0,s),2}_{[0,1]\times\Sigma}) \cap L^2(S^1, \tT^{(k,s),2}_{[0,1]\times\Sigma,L(t),L_1}). \label{Asepvariables}
\end{equation}
The notation involved, while systematic, is unfortunately a nightmare. The below diagram summarizes all the spaces involved with their appropriate topologies:


\begin{align}
\begin{split}
    \xy
(-15,1)*+{\tT_{\Sigma_0}^{k+1/2+s,2} \oplus \tT_{\Sigma_1}^{k+1/2+s,2}}="1";
(40,1)*+{\tT_{[0,1]\times\Sigma}^{(k+1,s),2}}="2";
(-15,-21)*+{B^{(k+1/2+s),2}(\tiL \oplus L_1)}="3";
(40,-21)*+{\tT_{[0,1]\times\Sigma,L,L_1}^{(k+1,s),2}}="4";
(75,-10)*+{\tT_{[0,1]\times\Sigma}^{(k,s),2}}="4.5";
{\ar_(0.35){r=(r_0,r_1)} "2";"1"};
{\ar^(0.4){r} "4";"3"};
{\ar^{\tH_0} "2";"4.5"};
{\ar_{\tH_0} "4";"4.5"};
{\ar@{^{(}->} (-15,-16)*{};"1"};    
{\ar@{^{(}->} (40,-16)*{};"2"};
(-17,-40)*+{\Maps^{(k+1/2,s),2}(S^1,\tT_{\Sigma_0} \oplus \tT_{\Sigma_1})}="b1";
(45,-40)*+{\Maps^{(k+1,s),2}(S^1,\tT_{[0,1]\times\Sigma})}="b2";
(-17,-66)*+{\Maps^{(k+1/2,s),2}(S^1,\tiL(t) \oplus L_1)}="b3";
(45,-66)*+{\Maps^{(k+1,s),2}(S^1,\tT_{[0,1]\times\Sigma,L(t),L_1})}="b4";
(100,-53)*+{\Maps^{(k,s),2}(S^1,\tT_{[0,1]\times\Sigma})\hspace{1in}}="b5";
{\ar_(0.45){r} "b2";"b1"};
{\ar^{r} "b4";"b3"};
{\ar^{\frac{d}{dt}+\tH_0} "b2";"b5"};
{\ar_{\frac{d}{dt}+\tH_0} "b4";"b5"};
{\ar@{^{(}->} (-17,-62)*{};"b1"};    
{\ar@{^{(}->} (45,-62)*{};"b2"};
\endxy
\end{split} \label{MapsCD2} \\ \nonumber
\end{align}

We want to study the operator
\begin{equation}
  \frac{d}{dt}+\tH_0: \Maps^{(k+1,s),2}(S^1, \tT_{[0,1]\times\Sigma,L(t),L_1}) \to \Maps^{(k,s),2}(S^1, \tT_{[0,1]\times\Sigma})
\end{equation}
where $L(t)$ is a family of tangent spaces to $\fL$ along a smooth path.  A priori, it is not at all obvious why this operator should be Fredholm and satisfy a corresponding estimate as in Theorem \ref{ThmSWlin}.  Indeed, we can no longer use simple self-adjointness techniques, since the Hilbert space $\tT^{(0,s),2}_{[0,1]\times\Sigma}$ no longer admits $\tH_0$ as an (unbounded) symmetric operator when $s > 0$. (Indeed, the inner product on $\tT^{(0,s),2}_{[0,1]\times\Sigma}$ is no longer defined in terms of the bundle metrics implicit in the definition of $\tT$ but contains operators in the $\Sigma$ direction which capture the anisotropy.) This means we no longer have the integration by parts formula (\ref{intS^1}), nor can we apply the spectral theorem as in (\ref{resolvent_est}) to understand the resolvent of $\tH_0$ on anisotropic function spaces.  However, all is not lost, since we can still prove an analogous estimate to (\ref{resolvent_est}) in the anisotropic setting.  This is because the resolvent of $\tH_0$, which is a resolvent associated to an elliptic boundary value problem, is a pseudodifferential type operator, and such an operator lends itself to estimates on a variety of function spaces, including anisotropic spaces.  Indeed, this is the reason we proved Theorem \ref{ThmSWlin} using resolvents as an alternative method to the integration by parts method, since the robust methods of pseudodifferential operator theory will carry over to anisotropic spaces.

On anisotropic spaces, the resolvent we wish to understand is the resolvent of
\begin{equation}
  \tH_0: \tT^{(1,s),2}_{[0,1]\times\Sigma, L_0,L_1} \to \tT^{(0,s),2}_{[0,1]\times\Sigma} \label{tH0-LL10},
\end{equation}
and we want the estimate
\begin{equation}
\|(i\tau - \tH_0)^{-1}\|_{B^{(0,s),2}([0,1]\times\Sigma)} \leq O(\tau^{-1}). \label{Aresolventest}
\end{equation}
As in the proof of Theorem \ref{ThmSWlin}, we assume here that (\ref{tH0-LL10}) is invertible and self-adjoint for $s = 0$, which we can always do by perturbing $\tH_0$ by a bounded operator.  By doing so, the resolvent in (\ref{Aresolventest}) makes sense for all $\tau \in \R$.

There is a well-developed theory for understanding the resolvent of elliptic boundary value problems, dating back to the work of Seeley in \cite{Se3}.  There the boundary conditions considered were differential and later extensions were made to pseudodifferential boundary conditions satisfying certain hypotheses (see e.g. \cite{GS} and \cite{G99}).  For us, the tangent space to a monopole Lagrangian is only ``nearly" pseudodifferential, in the sense that its tangent spaces are given by the range of projections which differ from a pseudodifferential operator by a smoothing operator (see Theorem I.2.1(ii)).  However, after a detailed analysis, one can adapt the methods of \cite{GS} and \cite{G99} to carry over to our present situation.  In carrying out this analysis, we should remark here that it is key that the operator (\ref{tH0-LL10}) is self-adjoint for $s=0$. We develop a sufficiently general framework for the construction of the resolvent of an elliptic boundary value problem in \cite{N0}, and via Theorem 15.32 and Corollary 15.34 there, we prove the resolvent estimates we need on anisotropic function spaces.

Having made the above remarks, let us finally state the generalization of Theorem \ref{ThmSWlin} to the anisotropic situation.  In order to do this, we have to introduce the anisotropic version of Definition \ref{DefStraight}, so that the maps which straighten the domains are well-behaved on anisotropic spaces.  For this purpose, define the subspace
$$\tT_{[0,1]\times\Sigma, L} = \{x \in \tT_{[0,1]\times\Sigma}, \; r_0(x) \in \wt{L}\} \subset \tT_{[0,1]\times\Sigma},$$
where we only impose boundary conditions on $\Sigma=\Sigma_0$.

\begin{Def}\label{DefAStraight}
  Let $L(t)$ be a smoothly varying family of subspaces\footnote{See Definition I.A.3.} of $\T_\Sigma$, $t \in S^1$ or $\R$.  We say that the $L(t)$ are \textit{anisotropic regular} if for every $t_0 \in \R$, there exists an open interval $I \ni t_0$ such that for every $t \in I$, there exist isomorphisms $S(t): \tT_{[0,1]\times\Sigma} \to \tT_{[0,1]\times\Sigma}$ satisfying the following properties for all nonnegative integers $k$ and every $s \in [0,1]$:
  \begin{enumerate}
    \item The map $S(t)$ extends to an isomorphism $S(t): \tT^{(k,s),2}_{[0,1]\times\Sigma} \to \tT^{(k,s),2}_{[0,1]\times\Sigma}$.  Furthermore, $S(t)$ acts as the identity on $[1/2,1]\times\Sigma$, i.e., for every $(b,\phi,\xi) \in \tT^{(k,s),2}_{[0,1]\times\Sigma}$, we have
        $$S(t)\big((b,\phi,\xi)\big)\Big|_{[1/2,1]\times\Sigma} = (b,\phi,\xi)\Big|_{[1/2,1]\times\Sigma}.$$
    \item The map $S(t)$ straightens the family $L(t)$ in the sense that $S(t): \tT_{[0,1]\times\Sigma, L(t_0)}^{(k,s),2} \to \tT_{[0,1]\times\Sigma,L(t)}^{(k,s),2}$ is an isomorphism ($k \geq 1$).
    \item The commutator of $[D,S(t)]$ where $D: \tT \to \tT$ is any first order differential operator, is an operator bounded on $\tT^{(k,s),2}_{[0,1]\times\Sigma}$.
  \end{enumerate}
\end{Def}

Observe that if the family $L(t)$ is anisotropic regular then it is also regular, since we can take $s = 0$ in the above and extend $S(t)$ by the identity to the rest of $Y$.

\begin{Theorem}\label{ThmSWlinA}
  Assume the family of spaces $L(t) \subset \T_{\Sigma_0}$, $t \in S^1$, and the space $L_1 \subseteq \tT_{\Sigma_1}$ satisfy the following:
  \begin{enumerate}
    \item The $L(t)$ are anisotropic regular.
    \item The operators
    \begin{equation}
  \tH_0: \tT^{(k+1,s),2}_{[0,1]\times\Sigma, L_0,L_1} \to \tT^{(k,s),2}_{[0,1]\times\Sigma}, \qquad t \in S^1, \label{AFredholm}
\end{equation}
are Fredholm for all $k \geq 0$ and $s \in [0,1]$, and moreover, for $k=s=0$, are self-adjoint.
\item The resolvent estimate (\ref{Aresolventest}) holds for all $s \in [0,1]$.
  \end{enumerate}
  Then
\begin{equation}
  \frac{d}{dt} + \tH_0 : \Maps^{(k+1,s),2}(S^1, \tT_{[0,1]\times\Sigma,L(t),L_1}) \to \HMaps^{(k,s),2}(S^1,\tT_{[0,1]\times\Sigma}) \label{AdtHs}
\end{equation}
is Fredholm for all $k \geq 0$ and $s \in[0,1]$, and we have the elliptic estimate
\begin{equation}
  \|x\|_{B^{(k+1,s),2}((S^1 \times [0,1])\times\Sigma)} \leq C\left(\left\|\left(\frac{d}{dt} + \tH_0\right)x\right\|_{B^{(k,s),2}((S^1 \times [0,1])\times\Sigma)} + \|x\|_{B^{(k,s),2}((S^1 \times [0,1])\times\Sigma)}\right). \label{Aest-dtHK}
\end{equation}
Furthermore, if $x \in \Maps^{(1,s),2}(S^1, \tT_{[0,1]\times\Sigma,L(t),L_1})$ satisfies $(\frac{d}{dt} + \tH_0)x \in \Maps^{(k,s),2}(S^1,\tT_{[0,1]\times\Sigma})$, then $x \in \Maps^{(k+1,s),2}(S^1, \tT_{[0,1]\times\Sigma,L(t),L_1})$ and it satisfies (\ref{Aest-dtHK}).
\end{Theorem}

\Proof The same proof of (\ref{controldt}) using the Fourier transform shows that (\ref{Aresolventest}) implies
\begin{equation}
  \|\partial_tx\|_{B^{(0,s),2}((S^1 \times [0,1])\times\Sigma)} \leq C\|(\partial_t+\tH_0)x\|_{B^{(0,s),2}((S^1 \times [0,1])\times\Sigma)}, \label{controldtA}
\end{equation}
Here, we used that $B^{s,2}([0,1]\times\Sigma)$ is a Hilbert space, so that we have a Plancherel theorem on $B^{(0,s),2}((S^1 \times [0,1])\times\Sigma) = L^2(S^1, B^{s,2}([0,1]\times\Sigma))$.  From (\ref{controldtA}), the same reasoning we used to derive (\ref{est_dtH1''}) shows that we now the anisotropic estimate
\begin{equation}
  \|x\|_{B^{(1,s),2}((S^1 \times [0,1])\times\Sigma)} \leq C(\|(\partial_t + \tH_0)x\|_{B^{(0,s),2}((S^1 \times [0,1])\times\Sigma)} + \|x\|_{B^{(1,s),2}((S^1 \times [0,1])\times\Sigma)}).
\end{equation}
This proves (\ref{Aest-dtHK}) for $k = 0$. From here, commuting time derivatives as in the proof of (\ref{est_dtHk}) shows that (\ref{Aest-dtHK}) holds for all $k \geq 0$. For the last statement, see Remark \ref{RemER}.\End

Summarizing, we have studied the operator $\frac{d}{dt} + \tH_0$ on the manifold with boundary $S^1 \times Y$, both on the usual $L^2$ spaces and on anisotropic spaces.  We have listed general properties that families of subspaces $L(t)$, serving as boundary conditions for $\frac{d}{dt}+\tH_0$, should satisfy if the spaces
\begin{equation}
  \Maps(S^1, \tT_{L(t)}), \qquad \Maps(S^1, \tT_{[0,1]\times\Sigma,L(t),L_1}),
\end{equation}
in suitable function space completions, are to yield a domain for which the operator $\frac{d}{dt} + \tH_0$ is Fredholm.  We have phrased matters in this generality, because this is the general model for the linearized Seiberg-Witten equations with Lagrangian boundary conditions, where the $L(t)$ come from linearizing a path along the Lagrangian $\fL$.  This allows us to understand the general framework for these equations, in particular, which Lagrangians are suitable for a well-posed boundary value problem.  We make further remarks on this and related issues in Section 5.

Of course, all of this discussion would be fruitless if we could not produce any examples of Lagrangians which satisfy the properties we have imposed.  Fortunately, for $\fL$ a monopole Lagrangian, all the properties we have used in Theorem \ref{ThmSWlinA} hold.  We have the following theorem.

\begin{Theorem}\label{ThmMonLag}
  Let $\fL$ be a monopole Lagrangian and let $\g \in \Maps(S^1,\fL)$ be a smooth path.  Define $L(t) = T_{\g(t)}\fL$, $t \in S^1$.  Then there exists $L_1 \subset \tT_{\Sigma_1}$ such that all the hypotheses (i)-(iii) of Theorem \ref{ThmSWlinA} hold.
\end{Theorem}

\Proof (i) Since $\fL \subset \fC(\Sigma)$ is a submanifold, then the tangent spaces to any smooth path $\g \in \Maps(S^1,\fL)$ automatically form a smoothly varying family of subpaces of $\T_\Sigma$.  By Theorem I.2.1(ii), we know that each tangent space $L(t) = T_{\g(t)}\fL$ is the range of a ``Calderon projection"
\begin{equation}
  P^+(t) := P^+_{\tilde \g(t)}: \T_\Sigma \to \T_\Sigma.
\end{equation}
Here $\tilde \g \in \Maps(I,\M)$ is a smooth path that lifts $\g \in \Maps(I,\L)$, i.e. $r_\Sigma(\tg(t)) = \g(t)$ for all $t$.  That such a smooth lift exists follows from the techniques used in the proof of Theorem \ref{ThmPathL1}(iv).   For each $t$, the resulting projection $P^+(t)$ extends to a bounded map, in particular, on $\T^{s,2}_\Sigma$ for all $s \geq 0$. Furthermore, it differs from a pseudodifferential projection $\pi^+$ by an operator $T(t) := P^+(t) - \pi^+$ that smooths by one derivative, i.e.,  $T(t): T^{s,2}_\Sigma \to T^{s+1,2}_\Sigma$ for all $s \geq 0$.  Indeed, because $\tilde\g$ is smooth, one can check from the arithmetic of Theorem I.2.1 that the maps $P^+(t)$ and $T(t)$ have the mapping properties on all the function spaces just stated.

To construct the straightening maps $S(t)$ in Definition \ref{DefAStraight}, we construct straightening maps on the boundary using Lemma I.A.5, and then extend these to maps in a collar neighborhood of the boundary in a slicewise fashion.  In detail, given any $t_0 \in S^1$, say $t_0 = 0$, Lemma I.A.5 tells us that there exists a time interval $I \ni 0$ such that we have straightening maps
\begin{align*}
  S_\Sigma(t): \T_\Sigma & \to \T_\Sigma, \qquad t \in I,
\end{align*}
with each $S_\Sigma(t)$ is an isomorphism and $S_\Sigma(t)(L(0)) = L(t)$.  Moreover, since the closures $B^{s,2}L(t)$ are also smoothly varying and complemented in $\T_\Sigma^{s,2}$ for all $s \geq 0$ (again by the Theorem I.2.1), the straightening maps extend to the Besov closures as well, and we have
\begin{align}
  S_\Sigma(t): \T^{s,2}_\Sigma & \to \T^{s,2}_\Sigma, \qquad t \in I,\; s \geq 0, \label{eq:SSigma}
\end{align}
with $S_\Sigma(t)(B^{s,2}L(0)) = B^{s,2}L(t)$.

We now use these boundary straightening maps to construct straightening maps on $\tT_{[0,1]\times\Sigma}$.  Let $(b,\phi,\xi) \in \tT_{[0,1]\times\Sigma}$.  In the collar neighborhood $[0,1]\times\Sigma$, let $v \in [0,1]$ be the inward normal coordinate and write $b = b_1 + \beta dv$ in terms of its tangential and normal components, respectively, where $b_1 \in \Gamma([0,1], \Omega^1(\Sigma;i\R))$ and $\beta \in \Gamma([0,1], \Omega^0(\Sigma;i\R))$.  We will use the $S_\Sigma(t)$ maps to straighten out the tangential components $(b_1,\phi)|_\Sigma$ and we need not do anything to the $\beta_1$ and $\xi$.  More precisely, let $h = h(v)$ be a smooth cutoff function, $0 \leq h(v) \leq 1$, where $h(v) = 1$ for $v \leq 1/4$ and $h = 0$ for $v \geq 1/2$.  We define
\begin{align}
  S(t): \tT_{[0,1]\times\Sigma}& \to \tT_{[0,1]\times\Sigma} \nonumber\\
  (b,\phi,\xi) &\mapsto  h(v)S_\Sigma(t)(b_1(v),\phi(v)) + (1-h(v))(b_1(v),\phi(v)) + \beta dv + \xi, \label{S(t)},
\end{align}
where $S(t)$ is defined to be the identity outside of $[0,1]\times\Sigma$. Thus, at $v = 0$, $S(t)$ acting on the tangential components $(b_1(0),\phi(0))$ is just the map $S_\Sigma(t)$, on $v \geq 1/2$, the map $S(t)$ is the identity, and in between, we linearly interpolate. For all $v$, we do nothing to $\beta dv$ and $\xi$. We now have to check that all the properties of Definition \ref{DefAStraight} hold.  For (i), we first have that (\ref{S(t)}) is an isomorphism in the smooth setting for all $t \in I$, where $I$ is sufficiently small.  Indeed, at $t = 0$, the map $S_\Sigma(0)$ is just the identity and hence so is $S(0)$.  For small enough $t$, the map $S_\Sigma(t)$ is sufficiently close to the identity that all its linear interpolants with the identity map on $\T_\Sigma$ are still isomorphisms.  Shrinking $I$ if necessary, it follows that (\ref{S(t)}) is an isomorphism for all $t\in I$.  It remains to show that $S(t)$ is bounded on anisotropic Besov spaces.  However, this follows from a similar analysis as was done in Lemma \ref{LemmaSW}, since although $S(t)$ is not time-independent, it is smoothly so.  Indeed, the mapping properties of $S(t)$ are determined from $S_\Sigma(t)$ acting slicewise in the $v$ direction.  This latter slicewise map clearly acts on integer Sobolev spaces because of the Leibnitz rule, the smoothness of $S_\Sigma(t)$, and (\ref{eq:SSigma}).  The Fubini property and interpolation property of Besov spaces, as explained in Lemma \ref{LemmaSW}, now show that the $v$-slicewise $S_\Sigma(t)$ is bounded on anisotropic Besov spaces and hence so is $S(t)$. This proves (i) in Definition \ref{DefAStraight}.  Next, for (ii), it follows from $S_\Sigma(t)(B^{k+s-1/2,2}L(0)) = B^{k+s-1/2,2}L(t))$ and the boundedness of $S(t)$ on $\tT^{(k,s),2}_{[0,1],\Sigma}$ that
\begin{equation}
    S(t): \tT^{(k,s),2}_{[0,1]\times\Sigma,L(0)} \to \tT^{(k,s),2}_{[0,1]\times\Sigma,L(t)}, \qquad k \geq 1.
\end{equation}
 Finally, for the commutator property (iii) in Definition \ref{DefAStraight}, we need only work in the collar neighborhood $[0,1]\times\Sigma$, since $S(t)$ is the identity outside of it.  There are two cases.  If $D$ is a differential operator in the $\Sigma$ directions, then $[D,S(t)]$ is essentially given by $[D,P^+(t)]$ and $[D,T(t)]$, both of which yield bounded operators; the former because the commutator of a first and zeroth order pseudodifferential operator is a zeroth order pseudodifferential operator (hence bounded) and the latter because $T(t)$ is already smoothing by one derivative in the $\Sigma$ directions.  If $D$ is a differential operator in the $v$ direction, then $[S(t),D]$ is still a bounded operator, since we just differentiate the smooth cutoff function $h(v)$ in the commutator.  Altogether, this proves $S(t)$ satisfies all properties of Definition \ref{DefAStraight} and hence, the family $L(t)$ is anisotropic regular.

(ii) Using Lemma \ref{LemmaF0F1}, if we choose any $L_1$ such that (F1) is satisfied, then it suffices to prove that condition (F0) holds for $L_0 = L(t)$ for every $t$.  It is here where we use that the monopole Lagrangian $\fL = \mathcal{L}(Y')$ comes from a manifold $Y'$ such that $\partial Y' = -\partial Y$.  Following the discussion in Section 2 of \cite{N1}, on $Y$, the operator $\tH_0$ is a Dirac operator that decomposes as the sum of two Dirac operators,
\begin{equation}
  \tH_0 = D_{\dgc} \oplus D_{B_\rf}, \label{tH0-DD}
\end{equation}
the div-grad-curl operator
$$D_{\dgc} = \begin{pmatrix}
  *d & -d\\ d^* & 0
\end{pmatrix}: \Omega^1(Y;i\R) \oplus \Omega^0(Y;i\R) \circlearrowleft,$$
on differential forms, and the Dirac operator
$$D_{B_\rf}: \Gamma(\S)\to\Gamma(\S)$$
on spinors. Thus, the tangential boundary operator $\mathsf{B}: \tT_\Sigma \to \tT_\Sigma$ associated to $\tH_0$ on $\Sigma = \Sigma_0$ splits as a sum
\begin{equation}
  \mathsf{B} = \mathsf{B}_{\dgc} \oplus \mathsf{B}_{\S} \label{B=BB2}
\end{equation}
of the tangential boundary operators associated to $D_{\dgc}$ and $\D_{B_\rf}$, respectively.  We get an associated spectral decomposition of $\tT_\Sigma$ via
\begin{equation}
  \tT_\Sigma = \cZ^+ \oplus \cZ^- \oplus \cZ^0
\end{equation}
given by the positive, negative, and zero eigenspace decomposition of $\mathsf{B}$.  Furthermore, by (\ref{B=BB2}), we have
\begin{equation}
  \cZ^{\pm} = (\cZ^{\pm}_e \oplus \cZ^{\pm}_c) \oplus \cZ^{\pm}_\S \label{Zspaces2}
\end{equation}
where
$$\cZ^{\pm}_e \oplus \cZ^{\pm}_c \subset \Omega^1(\Sigma;i\R) \oplus \Omega^0(\Sigma;i\R)\oplus\Omega^0(\Sigma;i\R)$$
and
$$\cZ^{\pm}_\S\subseteq \Gamma(\S_\Sigma)$$
are the positive and negative eigenspaces associated to $\mathsf{B}_{\dgc}$ and $\mathsf{B}_{\S}$, respectively, see Lemma I.2.5 and (I.2.77). (Note that $\cZ^+$ is the same space as $\im \Pi^+$ appearing in Lemma \ref{LemmaF0F1}).

Thus, to show (F0) in Lemma \ref{LemmaF0F1}, we have to show that $B^{k+s+1/2,2}\wt{L(t)}$ is Fredholm with $B^{k+s+1/2,2}\cZ^+$.  By Theorem I.2.1(i), we have that $B^{k+s+1/2,2}L(t)$ is a compact perturbation of $B^{k+s+1/2,2}(\im d \oplus \cZ^-_\S)$ in $\T^{k+s+1/2,2}_\Sigma$.  Note the important minus sign in the last factor.  This minus sign arises because when we apply Theorem I.2.1, we apply it to the manifold $Y'$, and since we have the opposite orientation $\partial Y' = -\partial Y$, the tangential boundary operators for the operators on $Y'$ differ by a minus sign from the corresponding ones on $Y$.  Altogether then, $B^{k+s+1/2,2}\wt{L(t)}$ is a compact perturbation of
\begin{equation}
  B^{k+s+1/2,2}(\im d \oplus \cZ^-_\S \oplus 0 \oplus \Omega^0(\Sigma;i\R)) \subset \tT^{k+s+1/2,2}_\Sigma. \label{Y'space}
\end{equation}
From the definition of $\cZ^+_e$ and $\cZ^+_c$ in Lemma I.2.5, one can now easily see that (\ref{Y'space}) is Fredholm with $B^{k+s+1/2,2}\cZ^+$ via (\ref{Zspaces2}).  Thus, (\ref{Y'space}) and hence $B^{k+s+1/2,2}\wt{L(t)}$ is Fredholm with $B^{k+s+1/2,2}\cZ^+$.  So (F0) is satisfied, and this proves the Fredholm property of (\ref{AFredholm}) by Lemma \ref{LemmaF0F1}.

For $k=s=0$, the operator (\ref{AFredholm}) is symmetric since $B^{1/2,2}\wt{L(t)} \subseteq \tT^{1/2,2}_\Sigma$ is an isotropic subspace.  Here, we suppose the subspace $L_1$ chosen in (ii), which satisfies (F1), is such that  $B^{1/2,2}L_1 \subset \tT^{1/2,2}_{\Sigma_1}$ is a Lagrangian subspace. From this, it turns out that since $B^{1/2,2}\wt{L(t)} \subseteq \tT^{1/2,2}_\Sigma$ is a Lagrangian subspace which furthermore satisfies (F0), then (\ref{AFredholm}) is self-adjoint\footnote{Not every Lagrangian subspace of $\tT^{1/2,2}_\Sigma \oplus \tT^{1/2,2}_{\Sigma_1}$ will yield a self-adjoint operator.  For instance, $r(\ker \tH_0) \subset \tT^{1/2,2}_{\Sigma} \oplus \tT^{1/2,2}_{\Sigma_1}$ is a Lagrangian subspace, but this boundary condition is not a self-adjoint boundary condition for $\tH_0$.}.  This follows from the general and abstract framework of finding self-adjoint extensions of closed-symmetric operators by finding Lagrangian subspaces of a suitable quotient Hilbert space, dating back to classical work of von-Neumann.  Here, the relevant theorem is Theorem 22.4 of \cite{N0}.  This shows (\ref{AFredholm}) is self-adjoint for $k=0$.

(iii) This follows from Theorem I.2.1(ii) and Corollary 15.34 of \cite{N0}. \End

Thus, both Theorems \ref{ThmSWlin} and \ref{ThmSWlinA} hold for the linearized operator arising from the (gauge-fixed) Seiberg-Witten equations with Lagrangian boundary condition determined by a monopole Lagrangian.  In particular, from Theorem \ref{ThmSWlin}, the linearized operator of our boundary value problem is a Fredholm operator.  This proves

\begin{Theorem}\label{ThmFredProp} (Fredholm Property)
  Let $\fL$ be a monopole Lagrangian.  Consider the equations (\ref{BVP}) on $S^1 \times Y$, where we impose the gauge fixing condition
 \begin{equation}
   d^*(A - A_0) = 0, \qquad *(A-A_0)|_{S^1\times\Sigma}=0, \label{GFCN}
 \end{equation}
where $A_0$ is a smooth connection. Then the linearization of the system (\ref{BVP}) and (\ref{GFCN}) at a smooth configuration $(A,\Phi) = (B(t)+\alpha(t)dt,\Phi(t))$ determines an operator
  \begin{equation}
    \frac{d}{dt} + \tH_0  -D_{(A,\Phi)}N_{(A_0,\Phi_0)}: \Maps^{k+1,2}(S^1,\tT_{L(t)}) \to \Maps^{k,2}(S^1,\tT_{L(t)}), \label{linSWFred}
  \end{equation}
where $\tH_0$ and $N_{(A_0,\Phi_0)}$ are given by (\ref{eq3:tH0}) and (\ref{sec3:Nmap}), respectively, $L(t) = T_{r_\Sigma(B(t),\Phi(t))}\L$, and where $k \geq 0$.  The operator (\ref{linSWFred}) is a Fredholm operator for all $k \geq 0$.
\end{Theorem}

In particular, this means that if we have transversality for our Seiberg-Witten system (say by perturbing the equations in the interior in a mild way), then the moduli space of solutions modulo gauge to our boundary value problem is finite dimensional.

We have one more corollary to the above analysis which we will need in the next section, which yields for us the inhomogeneous version of Theorem \ref{ThmSWlinA}.

\begin{Corollary}\label{CorSubspaceFred}
  Let $\g \in \Maps(S^1,\fL)$ be a smooth path, let $L(t) = T_{\g(t)}\fL$, and let $L_1$ be as in Theorem \ref{ThmMonLag}.  Then we have the following:
  \begin{enumerate}
    \item The space
  \begin{align}
    \Maps^{(k+1/2,s),2}(S^1,\wt{L(t)}\oplus L_1) & := \{z \in \Maps^{(k+1/2,s),2}(S^1, \tT_{\Sigma_0}\oplus\tT_{\Sigma_1}) : \nonumber\\
    & \qquad\qquad\qquad z(t) \in \wt{L(t)}\oplus L_1, \textrm{ for all }t \in S^1\} \label{Fredholmspace}
  \end{align}
  is Fredholm with $r\left(\ker \left(\partial_t+\tH_0\right)\right) \subset \Maps^{(k+1/2,s),2}(S^1,\tT_{\Sigma_0}\oplus\tT_{\Sigma_1})$, for all $k\geq 0$ and $0 \leq s \leq 1$.
  \item The intersection of (\ref{Fredholmspace}) and $r\left(\ker \left(\partial_t+\tH_0\right)\right)$ is spanned by finitely many smooth elements, and the span of these two spaces is complemented by a space spanned by finitely many smooth elements.
  \item There exists a projection
  \begin{equation}
    \Pi: \Maps^{(k+1/2,s),2}(S^1,\tT_{\Sigma_0}\oplus\tT_{\Sigma_1}) \circlearrowleft \label{eq3:Pi}
  \end{equation}
  such that $\ker \Pi$ is (\ref{Fredholmspace}) and $\im \Pi \cap r\left(\ker \left(\partial_t+\tH_0\right)\right)$ is of finite codimension in \\ $r\left(\ker \left(\partial_t+\tH_0\right)\right)$.  The map $\Pi$ is independent of $k \geq 0$ and $s \in [0,1]$.  Moreover, we have the inhomogeneous elliptic estimate
  \begin{multline}
     \|x\|_{B^{(k+1,s),2}((S^1 \times [0,1])\times\Sigma)} \leq C\bigg(\left\|\left(\partial_t + \tH_0\right)x\right\|_{B^{(k,s),2}((S^1 \times [0,1])\times\Sigma)} + \\\|\Pi r(x)\|_{B^{(k+1/2,s),2}(S^1\times(\Sigma_0 \cup \Sigma_1))}
 + \|x\|_{B^{(k,s),2}((S^1 \times [0,1])\times\Sigma)}\bigg). \label{inhom-est}
  \end{multline}
  \end{enumerate}
\end{Corollary}

\Proof (i) The space (\ref{Fredholmspace}) is precisely the space of boundary values of the domain of (\ref{AdtHs}).  Since the operator (\ref{AdtHs}) is Fredholm, (\ref{Fredholmspace}) is Fredholm with $r\left(\ker \left(\frac{d}{dt}+\tH_0\right)\right)$, the boundary values of the kernel of
$$\frac{d}{dt}+\tH_0: \Maps^{(k+1,s),2}(S^1,\tT_{[0,1]\times\Sigma}) \to \Maps^{k,2}(S^1,\tT_{[0,1]\times\Sigma})$$
where no boundary conditions are imposed.  This follows from Theorem \ref{ThmEBP2}.

(ii) The intersection of the two spaces consists of smooth elements because we have the elliptic estimate (\ref{Aest-dtHK}) for all $k$, which tells us that all elements in the kernel of (\ref{AdtHs}) are smooth.  The same analysis applies to the adjoint problem, and so the cokernel of (\ref{AdtHs}) (that is, the orthogonal complement of its range) is also spanned by smooth configurations.  We now apply Theorem \ref{ThmEBP2}.

(iii) This follows from Theorem \ref{ThmEBP3}.  By (ii), the projection $\Pi$ differs from the projection $(1 - \Pi_\U)$ in (\ref{EBP-inhom}) by a smooth error (where $\U$ is taken to be (\ref{Fredholmspace})), and so (\ref{inhom-est}) follows from (\ref{EBP-inhom}). The map $\Pi$ is independent of $k$ and $s$, since the kernel and cokernel of \ref{AdtHs} are independent of $k$ and $s$ and the Fredholm property in (i) holds for every $k$ and $s$.\End


\section{Proofs of the Main Theorems}\label{SecP3Sec4}

In the previous section, we studied the linearized Seiberg-Witten equations with abstract Lagrangian boundary conditions.  In the course of doing so, we found that Lagrangians satisfying certain analytic properties yield elliptic estimates for the linearized Seiberg-Witten equations.  Furthermore, we showed that monopole Lagrangians satisfy all such properties.  Thus, with $\fL$ a monopole Lagrangian, we can now prove our main theorems, using both the linear analysis in the previous section, and the nonlinear analysis in Section 2 concerning the space of paths through $\fL$.  It is convenient to prove the results on $S^1 \times Y$.  A careful patching argument (the boundary condition is nonlinear) then proves the main results on $\R \times Y$.

Recall that when $p = 2$, the Besov spaces $B^{s,2}$ are the usual fractional Sobolev spaces $H^s = H^{s,2}$, for all $s \in \R$.  When $s$ is a nonnegative integer,  the spaces $B^{s,2} = H^{s,2}$ are also denoted by $W^{s,2}$.  We will use either notation $B^{s,2}$ or $H^{s,2}$ wherever convenient.

\begin{Theorem}\label{ThmA'}
  Let $p > 4$, and let $A = B(t) + \alpha(t)dt \in H^{1,p}\A(S^1 \times Y)$ and $\Phi \in H^{1,p}\Gamma(\S^+)$ solve the boundary value problem
\begin{align}
SW_4(A,\Phi) & = 0, \label{SWeq1}\\
  r_\Sigma(B(t),\Phi(t)) & \in \fL^{1-2/p,p}, \quad \forall t \in S^1, \label{SWeq3}
\end{align}
where $\fL$ is a monopole Lagrangian.  Then there exists a gauge transformation $g \in H^{2,p}\G_{\id}(S^1 \times Y)$ such that $g^*(A,\Phi)$ is smooth.  In particular, if $A$ is in Coulomb-Neumann gauge with respect to any smooth connection, i.e. (\ref{GFCN}) holds, then $(A,\Phi)$ is smooth.
\end{Theorem}

\Proof  From the previous section, we know we can find a gauge in which the equations are a semilinear elliptic equation (in the interior) as in (\ref{dtHeq}).  Since $H^{k,p}(S^1 \times Y)$ is an algebra for $k \geq 1$ and $p>4$, it follows that we can elliptic bootstrap in the interior to any desired regularity.  Thus, we need only prove regularity near the boundary.

From here, the proof proceeds in four main steps.  The first step is to rewrite the equations in a suitable gauge so that the linear portion of the equations satisfy all the hypotheses of the previous section.  In particular, the linearized equations now satisfy elliptic estimates on anisotropic function spaces, where the anisotropy is in the $\Sigma$ direction.  From this, the second step is to gain regularity for $(A,\Phi)$ in the $\Sigma$ directions in a neighborhood of the boundary.  Here, we use the results from Section 2, namely Theorems \ref{ThmPathL1} and \ref{ThmPathL2}, that the nonlinear part of the chart maps for the space of paths through $\fL$ smooth in the $\Sigma$ directions.  Moreover, it is here that the complicated choice of topologies appearing in these theorems, particularly in Theorem \ref{ThmPathL2}, will serve their purpose. Using the anisotropic linear theory of Section 3, specifically Corollary \ref{CorSubspaceFred}, we then gain regularity for $(A,\Phi)$ in the $\Sigma$ directions, where the linear theory can be applied because the nonlinear contribution from the boundary condition is smoothing in the $\Sigma$ directions.  The third step is to gain regularity in the time direction and normal direction to $\Sigma$ using results from Banach space valued Cauchy-Riemann equations due to Wehrheim \cite{We} which we adapt to our needs in Section \ref{AppCR}.  Once we have gained some regularity in all the directions, then in our final step, we bootstrap to gain regularity to any desired order.\\

\textit{Step One:} In the previous section, we found a suitable gauge in which the Seiberg-Witten equations become a semilinear elliptic equation with quadratic nonlinearity, namely, we obtained the system (\ref{dtHeq}) for the equations in the interior.  We wish to do the same here, only now $(A,\Phi)$ is not smooth.  Furthermore, we must choose the base configuration $(A_0,\Phi_0)$ about which we linearize our configuration $(A,\Phi)$ carefully.

So choose $(A_0,\Phi_0) \in \fC(S^1 \times Y)$ a smooth configuration close to $(A,\Phi)$ and satisfying the Lagrangian boundary condition (\ref{SWeq3}), where we will define this more precisely in a moment.  In the usual way, write $A_0 = B_0(t) + \alpha_0(t)dt$ as a path of connections $B_0(t)$ on $Y$ plus its temporal part $\alpha_0(t)$, and write $\Phi_0 = \Phi_0(t)$ as a path of spinors on $Y$.  Then we can find a gauge transformation $g \in H^{2,p}\G_{\id}(S^1 \times Y)$, that places $A$ in Coulomb-Neumann gauge with respect to $A_0$, i.e.,
\begin{equation}
  d^*(g^*A - A_0) = 0, \qquad *(g^*A-A_0)|_{S^1\times \Sigma} = 0.
\end{equation}
This $g$ is determined by writing $g = e^f$, where $f \in \Omega^0(S^1\times Y; i\R)$, and solving the inhomogeneous Neumann problem
\begin{align}
  \Delta f & = d^*(A - A_0)\\
  *df|_{S^1\times\Sigma} &= *(A-A_0)|_{S^1\times\Sigma}
\end{align}
This equation has a unique solution $f \in H^{2,p}\Omega^0(S^1\times Y;i\R)$, up to constants, by the standard elliptic theory of the Neumann Laplacian.

Redefine $(A,\Phi)$ by the gauge transformation so obtained above, so that we have
\begin{equation}
  d^*(A-A_0) = 0, \qquad *(A-A_0)|_{S^1\times\Sigma}=0. \label{CNgaugeA_0}
\end{equation}
We want to gain regularity for the difference $(A,\Phi) - (A_0,\Phi_0)$, which we can write as the triple
\begin{equation}
  (b,\phi,\xi) \in H^{1,p}\Maps(S^1,\Omega^1(Y;i\R) \oplus \Gamma(\S) \oplus \Omega^0(Y;i\R)), \label{eq:bpx}
\end{equation}
where $b = b(t)$ is $B(t) - B_0(t)$, $\phi = \phi(t)$ is $\Phi(t) - \Phi_0(t)$, and $\xi(t) = \alpha(t) - \alpha_0(t)$.

Our goal is to show that $(b,\phi,\xi)$ is smooth.  As shown in Section 3, the configuration $(b,\phi,\xi)$ satisfies (\ref{dtHeq}) and (\ref{NBC}). We now have to add in the nonlinear Lagrangian boundary condition (\ref{SWeq3}) to these equations. To express this in terms of $(b,\phi,\xi)$ requires that we choose $(A_0,\Phi_0)$ sufficiently close to $(A,\Phi)$, as we now explain.  Recall from Theorem \ref{ThmPathL1} that for any path $\g \in \Maps^{1-1/p,p}(S^1,\fL)$, there is a local chart map $\E_\g$ which maps a neighborhood of $0$ in the tangent space $T_\g\Maps^{1-1/p,p}(S^1,\fL)$ diffeomorphically onto a neighborhood of $\g \in \Maps^{1-1/p,p}(S^1,\fL)$.  Furthermore, by construction, the chart maps contain a $C^0(S^1,B^{s',p}(\Sigma))$ neighborhood of $\Maps^{1-1/p,p}(S^1,\fL)$, for any $1/2 < s' \leq 1-2/p$, and the size of this neighborhood can be chosen uniformly on small $C^0(S^1,B^{s',p}(\Sigma))$ neighborhoods of $\g$. It follows that if the (smooth) $(A_0,\Phi_0)$ is sufficiently $H^{s,p}(S^1 \times Y)$ close to $(A,\Phi)$, with $s > 1/2+2/p$,  then on the boundary, the associated smooth path
$$\g_0:=\widehat{r_\Sigma}(B_0(t),\Phi_0(t)) \in \Maps(S^1,\fL)$$
is sufficiently $C^0(S^1,B^{s',p}(\Sigma))$ close to
$$\g=\w{r_\Sigma}(B(t),\Phi(t)) \in \Maps^{1-1/p,p}(S^1,\fL)$$
so that we can find a unique $z \in T_{\g_0}\Maps^{1-1/p,p}(S^1,\fL)$ near $0$ satisfying
\begin{equation}
  \g = \E_{\g_0}(z). \label{eq:defz}
\end{equation}
Rewriting this in terms of the map $\E_{\g_0}^1$ in Theorem \ref{ThmPathL1}, we thus have
\begin{equation}
  \g = \g_0+z+\E^1_{\g_0}(z), \qquad z \in T_{\g_0}\Maps^{1-1/p,p}(S^1,\fL). \label{gz}
\end{equation}
In other words, we have placed $\g$ in the range of the chart map $\E_{\g_0}$ centered at the smooth configuration $\g_0$.  Altogether, the interior equations (\ref{dtHeq}), the Neumann boundary condition on $b$ (\ref{NBC}), and the boundary condition (\ref{gz}) yield the following form for the full system of Seiberg-Witten equations with Lagrangian boundary conditions:
\begin{align}
   \left(\frac{d}{dt} + \tH_0\right)(b,\phi,\xi) &= N_{(A_0,\Phi_0)}(b,\phi,\xi) - SW_4(A_0,\Phi_0) \label{SWeqN1}\\
   (b,\phi)|_{S^1\times\Sigma} &= z+\E^1_{\g_0}(z), \qquad z \in T_{\g_0}\Maps^{1-1/p}(S^1,\fL) \label{SWeqN2}\\
   *b|_{S^1\times\Sigma} &= 0. \label{SWeqN3}
\end{align}
Recall that $N_{(A_0,\Phi_0)}$ is a quadratic multiplication map and $SW_4(A_0,\Phi_0)$ is a smooth term since $(A_0,\Phi_0)$ is smooth. Since we are only interested in regularity near the boundary, it suffices to gain regularity for a smooth truncation of $(b,\phi,\xi)$ with support near the boundary.  Thus, define
\begin{equation}
  (b_0,\phi_0,\xi_0) = \chi(b,\phi,\xi)
\end{equation}
where $\chi$ is a smooth cutoff function supported in a collar neighborhood $S^1\times[0,1]\times\Sigma$ of the boundary, with $\chi \equiv 1$ on $S^1 \times [0,1/2]\times\Sigma$ and $\chi \equiv 0$ on outside of $S^1 \times [0,3/4]\times\Sigma$.  Thus, via the notation of Section 3, we have
\begin{align*}
  (b_0,\phi_0,\xi_0) \in H^{1,p}\Maps(S^1,\tT_{[0,1]\times\Sigma}).
\end{align*}

\textit{Step Two:} We will gain regularity for $(b_0,\phi_0,\xi_0)$ in the $\Sigma$ directions from the equations (\ref{SWeqN1})-(\ref{SWeqN3}).
For this, we will use the linear theory on $L^2$ spaces developed in the previous section.  The main idea is simple.  The boundary condition (\ref{SWeqN2}) and (\ref{SWeqN3}) is essentially a perturbation of the linear boundary condition studied in Section 3.  Indeed, Theorems \ref{ThmSWlin} and \ref{ThmSWlinA} give us elliptic estimates when the nonlinear term $\E^1_{\g_0}(z)$ in (\ref{SWeqN2}) is absent, since then the boundary condition (\ref{SWeqN2}) satisfies the linear boundary conditions of Thereoms \ref{ThmSWlin} and \ref{ThmSWlinA}, where $L(t) = T_{\g_0(t)}\Maps(S^1,\fL)$.  Moreover, since $(b_0,\phi_0,\xi_0)$ is supported on $S^1\times[0,3/4]\times\Sigma$, the $(b_0,\phi_0,\xi_0)$ satisfy any boundary condition on $\{1\}\times\Sigma$, so that we may use Theorem \ref{ThmSWlinA} for any suitable ``dummy" boundary condition $L_1$. With the nonlinear term in (\ref{SWeqN2}) however, we use the fact that $\E^1_{\g_0}$ is smoothing in the $\Sigma$ directions, as given by Theorem \ref{ThmPathL1}.  Thus, we are able to gain regularity in the $\Sigma$ directions using the inhomogeneous elliptic estimate in Corollary \ref{CorSubspaceFred}.

In detail, we first have to embed the $p \neq 2$ Besov spaces into the $p=2$ Besov spaces of Section 3.  For this, we use the embedding $B^{s,p}(X) \hookrightarrow B^{s-\eps,2}(X)$ for any $s \in \R$, $p > 2$, and $\eps > 0$, on any compact manifold $X$.  In particular, we have $B^{1-1/p,p}(S^1\times\Sigma) \hookrightarrow B^{1-1/p-\eps,2}(S^1\times\Sigma)$.  Consequently, using Theorem \ref{ThmPathL1}, we have
\begin{equation}
  \cE^1_{\g_0}(z) \in \Maps^{(1-1/p,1-1/p-\eps),p}(S^1,\T_\Sigma) \hookrightarrow \Maps^{(1-1/p-\eps, 1-1/p-\eps),2}(S^1,\T_\Sigma). \label{E1embed}
\end{equation}
With $L(t) = T_{\g_0(t)}\Maps(S^1,\fL)$ and $L_1 \subset \T_{\Sigma_1}$ as in Corollary \ref{CorSubspaceFred}, we have a projection \begin{equation}
    \Pi_{\g_0}: \Maps^{(k+1/2,s),2}(S^1,\tT_{\Sigma_0}\oplus\tT_{\Sigma_1}) \circlearrowleft \label{eq4:Pi}
\end{equation}
defined by (\ref{eq3:Pi}) in Corollary \ref{CorSubspaceFred}.  By definition of $z$ in (\ref{eq:defz}) and the construction of $\Pi$, we have that $\Big((z,\xi|_{S^1\times\Sigma}),0\Big) \in \Maps^{1-1/p,p}(S^1,\tT_\Sigma \oplus \tT_{\Sigma_1})$ satisfies
\begin{equation}
  \Pi_{\g_0}\Big((z,\xi|_{S^1\times\Sigma}),0\Big)=0. \label{eq4:kerPi}
\end{equation}
Let $x = (b,\phi,\xi)$ and $x_0 = (b_0,\phi_0,\xi_0)$.  Thus applying the estimate (\ref{inhom-est}) to $x_0 = (b_0,\phi_0,\xi_0)$ and using equations (\ref{SWeqN1})-(\ref{SWeqN3}) and (\ref{eq4:kerPi}), we have
\begin{align}
     \|x_0\|_{B^{(k+1,s),2}((S^1 \times [0,1])\times\Sigma)} & \leq C\bigg(\left\|N_{(A_0,\Phi_0)}(x)\right\|_{B^{(k,s),2}((S^1 \times [0,1])\times\Sigma)} + \|\Pi \E^1_{\g_0}(z)\|_{B^{(k+1/2,s),2}(S^1\times(\Sigma \cup \Sigma_1))}\nonumber \\
       & \qquad + \|x\|_{B^{(k,s),2}((S^1 \times [0,1])\times\Sigma)} + \|SW_4(A_0,\Phi_0)\|_{B^{(k,s),2}((S^1 \times [0,1])\times\Sigma)}\bigg). \label{Sigmaboot}
\end{align}
for all $k \geq 1$ and $s \in [0,1]$ such that the right-hand side is finite (see also Remark \ref{RemER}).  

First, let $k=0$ and $s = 1$ in (\ref{Sigmaboot}).  Since $x \in H^{1,p}((S^1\times [0,1])\times\Sigma) \hookrightarrow B^{(0,1),2}(\nbhd)$, this means we always have control of the lower order third term of (\ref{Sigmaboot}).  Furthermore, since $H^{1,p}(S^1 \times Y)$ is an algebra, we have
\begin{equation}
  N(x) \in H^{1,p}(S^1 \times [0,1]\times\Sigma) \hookrightarrow B^{(0,1),2}(S^1\times[0,1]\times\Sigma). \label{N-est}
\end{equation}
since $p > 4$.  The final term to control is the boundary term, for which we have
\begin{equation}
  \Pi(\E^1_{\g_0}(z)) \in \Maps^{(1/2,1),2}(S^1,\tT). \label{Pibound1}
\end{equation}
Here, we used (\ref{E1embed}) and the embedding
$$\Maps^{(1-1/p-\eps, 1-1/p-\eps),2}(S^1,\T_\Sigma) \hookrightarrow \Maps^{(1/2,1),2}(S^1,\tT_\Sigma),$$
which follows since $p > 4$ and $\eps > 0$ can be chosen small.

Thus, $\|\Pi \E^1_{\g_0}(z)\|_{B^{(1/2,1),2}(S^1\times\Sigma)}$ is bounded thanks to the smoothing property of $\cE^1_{\g_0}(z)$.  Furthermore, the boundedness of $\Pi$ and the preceding analysis imply that
\begin{align}
  \|\Pi \E^1_{\g_0}(z)\|_{B^{(1/2,1),2}(S^1\times\Sigma)} & \leq \|\E^1_{\g_0}(z)\|_{B^{(1-1/p,1-1/p-\eps),p}(S^1\times\Sigma)} \nonumber \\
  & \leq \mu_{\g_0}( \|z\|_{B^{1-1/p,p}(S^1\times\Sigma) } ) \label{Pibound2.0} \\
  & \leq \mu_{\g_0}(\|x\|_{H^{1,p}(S^1\times\Sigma)}) \label{Pibound2}
\end{align}
where $\mu = \mu_{\g_0}: \R^+ \to \R^+$ is some continuous nonlinear function depending on $\g_0$ with $\mu(0)=0$.  It suffices to prove (\ref{Pibound2.0}), since the second line is just the trace theorem.  The main point is that even though Theorem \ref{ThmPathL1} tells us that $\E^1_{\g_0}(z) \in \Maps^{(1-1/p,1-1/p-\eps),p}(I,\T_\Sigma)$ given that $z \in \Maps^{1-1/p,p}(I,\T_\Sigma)$, we want an estimate in term of norms, as expressed by (\ref{Pibound2.0}).  However, it is plain to see that one can indeed obtain a norm estimate if one follows through all the various operators and constructions used in defining $\E^1_{\g_0}$.  In a few words, we obtain norm estimates because all estimates we perform along the way are derived from multiplication theorems, elliptic bootstrapping, interpolation, etc., all of which provide explicit norm dependent estimates. Hence, this proves (\ref{Pibound2.0}) and therefore (\ref{Pibound2}).\\

\noindent\textit{Notation: }In what follows, we write $\mu: \R^+\to\R^+$ to denote any continuous nonlinear function. Any subscripts on $\mu$ will be quantities which $\mu$ depends on which we wish to make explicit.  The precise form of $\mu$ is immaterial and may change from line to line.\\

Altogether, we have from (\ref{Sigmaboot}), (\ref{N-est}), and (\ref{Pibound2}) that
\begin{equation}
  \|(b_0,\phi_0,\xi_0)\|_{B^{(1,1),2}((S^1 \times [0,1])\times\Sigma)} \leq \mu(\|(b,\phi,\xi)\|_{H^{1,p}(S^1\times Y)}) + \|SW_4(A_0,\Phi_0)\|_{B^{(0,1),2}((S^1\times[0,1])\times\Sigma)}. \label{x0step2}
\end{equation}
Thus, we have gained a whole derivative in the $\Sigma$ direction, albeit with integrability exponent $2$ instead of $p$.\\

\textit{Step Three:} To gain regularity in the temporal and normal directions $S^1$ and $[0,1]$, respectively, we use the methods of \cite{We} which studies Cauchy-Riemann equations with values in a Banach space.  Here, the main results we need are summarized in Theorem \ref{ThmCR}, which is a refinement of \cite[Theorem 1.2]{We} to our situation.  Let us set up the notation for this analysis.

We have by definition
\begin{multline*}
  (b,\phi,\xi) \in H^{1,p}\Maps(S^1, \tT_{[0,1]\times\Sigma}) =\\ H^{1,p}\Maps\left(S^1, \Omega^1([0,1]\times\Sigma; i\R) \oplus \Gamma(\S_{[0,1]\times\Sigma}) \oplus \Omega^0([0,1]\times\Sigma;i\R)\right).
\end{multline*}
Let
\begin{equation}
  K = S^1 \times [0,1]
\end{equation}
and let $(t,v)$ be the corresponding temporal and normal coordinates on $K$.  Observe that we have the identification
\begin{align}
  \Maps(S^1, \tT_{[0,1]\times\Sigma}) & \cong \Gamma\left(K; \Omega^1(\Sigma; i\R) \oplus \Gamma(\S_\Sigma) \oplus
  \Omega^0(\Sigma; i\R) \oplus \Omega^0(\Sigma; i\R)\right)\\
  &= \Gamma(K,\tT_\Sigma),
\end{align}
via the restriction map
\begin{equation}
  r_v: \tT_{[0,1]\times\Sigma} \to \tT|_{\{v\}\times\Sigma} \cong \tT_\Sigma = \Omega^1(\Sigma; i\R) \oplus \Gamma(\S_\Sigma) \oplus \Omega^0(\Sigma; i\R) \oplus \Omega^0(\Sigma; i\R) \label{rv}
\end{equation}
induced by the restriction map $r$ in (\ref{restr}) for each $v \in [0,1]$.

Thus, we can regard $(b,\phi,\xi)$, a path from $S^1$ into $\tT_{[0,1]\times\Sigma}$, as a map from $K$ to $\tT_\Sigma$, in the appropriate function space topologies. Since we have the embeddings
$$H^{1,p}(K \times \Sigma) \hookrightarrow  C^0(S^1, B^{1-1/p,p}([0,1] \times \Sigma)) \hookrightarrow  C^0(K, B^{1-2/p,p}(\Sigma)),$$
and
$$B^{1-2/p,p}(\Sigma) \hookrightarrow L^p(\Sigma)$$
for $p > 2$, we can thus regard
\begin{equation}
  (b,\phi,\xi) \in H^{1,p}(K; L^p\tT_\Sigma). \label{H1pBanach}
\end{equation}
The space $H^{1,p}(K; L^p(\tT_\Sigma)$ is the Sobolev space of $H^{1,p}(K)$ functions with values in the Banach space $L^p\tT_\Sigma$.    As it turns out, we want to consider the larger space $H^{1,2}(K; L^2\tT_\Sigma) = B^{1,2}(K; L^2\tT_\Sigma)$, and we will instead regard
\begin{equation}
  (b,\phi,\xi) \in B^{1,2}(K; L^2\tT_\Sigma). \label{BesovBanach}
\end{equation}
We want to use (\ref{BesovBanach}) instead of (\ref{H1pBanach}) because the regularity in $\Sigma$ we gained in Step Two were with $p = 2$ Besov spaces.  This gain in regularity becomes essential when we reformulate the Seiberg-Witten equations as a nonlinear Cauchy-Riemann equation as we now explain.

Near the boundary, the Cauchy-Riemann operator occuring for us arises from the $t$ and $v$ derivatives of the operator $\frac{d}{dt} + \tH_0$ occurring in (\ref{SWeqN1}). Indeed, $\frac{d}{dt}+\tH_0$ is a Dirac operator, and so near the boundary where the metric is of the product form $g^2 = dt^2 + dv^2 + g_{\Sigma,v}^2$, we can write
\begin{equation}
  \frac{d}{dt}+\tH_0 = \frac{d}{dt}+J\frac{d}{dv} + D_\Sigma,
\end{equation}
where $J: \tT_\Sigma \to \tT_\Sigma$ is a smooth, bundle automorphism satisfying $J^2=-1$, and $D_\Sigma$ is a $v$-dependent differential operator acting on $\tT_\Sigma$.  Since we have gained regularity in the $\Sigma$ directions for $(b,\phi,\xi)$ in the previous step, the $\Sigma$ derivatives of $\frac{d}{dt} + \tH_0$ can be absorbed into $(b,\phi,\xi)$ and moved to the right-hand-side of (\ref{SWeqN1}).  Thus, (\ref{SWeqN1}) yields the semilinear Cauchy-Riemann equation
\begin{equation}
  \left(\frac{d}{dt} + J\frac{d}{dv}\right)(b,\phi,\xi) = -D_\Sigma(b,\phi,\xi) + N_{(A_0,\Phi_0)}(b,\phi,\xi) - SW_4(A_0,\Phi_0). \label{SWCR}
\end{equation}

In this setting, we reinterpret the boundary conditions (\ref{SWeqN2}) and (\ref{SWeqN3}) as follows. Recall that the configuration $(A,\Phi) = (A_0,\Phi_0) + (b,\phi,\xi)$ and the smooth configuration $(A_0,\Phi_0)$ both satisfy the Lagrangian boundary conditions (\ref{SWeq3}).  Thus, both $r_\Sigma(b(t),\phi(t)) + r_\Sigma(B_0(t),\Phi_0(t))$ and $r_\Sigma(B_0(t),\Phi_0(t))$ are elements of $L^2\fL \supset \L^{1-2/p,p}$ for every $t \in S^1$, by the above analysis.  Observe that for any pair of configurations $u,u_0 \in \fC^{0,2}(\Sigma)$, we have $u-u_0 \in \T^{0,2}_\Sigma$ since $\fC^{0,2}(\Sigma)$ is an affine space modeled on $\T^{0,2}_\Sigma$.  In particular, if $u_0 = (0,0)$ is the zero connection and zero spinor\footnote{For convenience, we assume the spinor bundle $\S_\Sigma$ on $\Sigma$ is trivial. This merely simplifies the notation since the reference configuration $u_0$ can be chosen to be zero.}, we can regard $\fC^{0,2}(\Sigma) = \T^{0,2}_\Sigma$.  In particular, we may regard $L^2\fL \subset \T^{0,2}_\Sigma$ and we may regard
\begin{equation}
L^2\wt{\fL} := L^2\fL \times L^2\Omega^0(\Sigma; i\R)
\end{equation}
as a subset of $\tT^{0,2}_\Sigma$.  Moreover, we may regard $r_\Sigma(B_0(t),\Phi_0(t))$ as a continuous path in $\tT^{0,2}_\Sigma$.

With these identification, the boundary conditions (\ref{SWeqN2}) and (\ref{SWeqN3}) can be expressed along the boundary of $K$, i.e., at $v = 0$, as
\begin{equation}
  r_\Sigma(B_0(t),\Phi_0(t)) + r_0(b(t),\phi(t),\xi(t)) \in L^2\wt{\fL}, \qquad\textrm{for all } t \in S^1.\label{CRBC}
\end{equation}
Thus, the boundary condition (\ref{CRBC}) captures the tangential Lagrangian boundary condition via the $L^2\fL$ factor of $L^2\wt{\fL}$, and it captures the Neumann boundary condition on $b$ via the remaining $L^2(\Omega^0(\Sigma; i\R))$ factor of $L^2\wt{\fL}$.

Altogether, we have a semilinear Cauchy-Riemann equation (\ref{SWCR}) with values in a Banach space $\tT_\Sigma^{0,2}$ and with boundary condition specified by (\ref{CRBC}). We can apply Theorem \ref{ThmCR} when the boundary condition (\ref{CRBC}) is given by a Banach manifold modeled on a closed subspace of an $L^p$ space.  In \cite{N1}, we studied the $L^p$ closure of $\fL$ for $p\geq 2$ and showed that, while we do not know if globally $L^p\fL$ is a manifold, we know that locally the chart maps $E_{u_0}$ for $\fL$ at a smooth configuration $u_0 \in \fL$ are bounded in the $L^p$ topology (see Corollary I.3.2).  More precisely, for every $2 \leq p < \infty$, there exists an $L^p(\Sigma)$ neighborhood $U$ of $0$ in $L^pT_{u_0}\fL$ containing an $L^p$ open ball, such that $E_{u_0}$ extends to a bounded map
\begin{equation}
  E_{u_0}: U \to L^p\fC(\Sigma) \label{ELp}
\end{equation}
which is a diffeomorphism onto its image. Moreover, because of the trace map
$$H^{1,p}(S^1 \times Y) \hookrightarrow C^0(K, B^{1-2/p,p}(\Sigma))$$
and because $\fL^{1-2/p,p} = B^{1-2/p,p}\fL$ is globally a smooth embedded submanifold of $\fC^{1-2/p,p}(\Sigma)$ (by \cite{N1} since $p > 4$), we know that the path
\begin{equation}
  \Big(t \mapsto r_\Sigma(b(t),\phi(t))+ r_\Sigma(B_0(t),\Phi_0(t))\Big) \in C^0(S^1, \fL^{1-2/p,p}) \label{CRpath}
\end{equation}
forms a continuous path in $\fL^{1-2/p,p}$, and hence on a small time interval $I \subset S^1$, the path lies in a single coordinate chart of a fixed configuration $u_0 \in \fL^{1-2/p,p}$ which we may take to be smooth.  In fact, we may as well take $u_0 = r_\Sigma(B_0(t_0),\Psi_0(t_0))$ for some fixed $t_0 \in I$.  Thus, we may replace (\ref{CRBC}), which may not be a manifold boundary condition in general, with
\begin{equation}
  r_0(b(t),\phi(t),\xi(t)) +  r_\Sigma(B_0(t),\Phi_0(t)) \in \widetilde{E_{u_0}(U)},  \qquad U \subset L^2T_{u_0}\fL, \textrm{ for all }t \in I\label{CRBCLp}
\end{equation}
where $U$ is an $L^2(\Sigma)$ open neighborhood of $0 \in L^2T_{u_0}\fL$ and
$$\widetilde{E_{u_0}(U)} := E_{u_0}(U) \times 0 \times L^2\Omega(\Sigma;i\R).$$
By the above remarks, (\ref{CRBCLp}) is a manifold boundary condition, since $E_{u_0}(U)$ is a submanifold of $L^2\fC(\Sigma)$. In effect, we have simply replaced a neighborhood of $u_0 \in \fL^{1-2/p,p} \subset \fC^{1-2/p,p}(\Sigma)$ with the larger $L^2$ neighborhood $E_{u_0}(U) \subset L^2\fC(\Sigma)$.
Moreover, since $\fL \subseteq \fC(\Sigma)$ is a Lagrangian submanifold, then $E_{u_0}(U)$ is a Lagrangian submanifold of $L^2\fC(\Sigma)$. Thus, $\wt{E_{u_0}(U)}$ is a product Lagrangian submanifold of
\begin{equation}
  L^2(\fC(\Sigma) \times \Omega^0(\Sigma; i\R) \times \Omega^0(\Sigma;i\R)), \label{eq4:prodspace}
\end{equation}
where the symplectic form on (\ref{eq4:prodspace}) is given by the product symplectic form (\ref{eq3:tomega}).  The time interval $I$ in (\ref{CRBCLp}) is chosen small enough so that the configuration $r_0(b(t),\phi(t),\xi(t)) + r_\Sigma(B_0(t),\Phi_0(t))$ remains inside the product chart $\widetilde{E_{u_0}(U)}$.  To simplify the below analysis, we can just suppose $I = S^1$. Otherwise, we can cover $S^1$ with small time intervals and sum up the estimates on each interval all the same.

Altogether, equations (\ref{SWCR}) and (\ref{CRBCLp}) form a Cauchy-Riemann equation for a configuration with values in a Banach space supplemented with Lagrangian boundary conditions.  Here, the Lagrangian submanifold $\widetilde{E_{u_0}(U)}$ is modeled on a closed subspace of an $L^2$ space.  Furthermore, it is an analytic Banach submanifold of $\tT_\Sigma^{0,2}$.  This is because the chart map $E_{u_0}$, by Theorem I.3.2, is constructed from the local straightening map $F_\co^{-1}$, where $\co \in \M$ satisfies $r_\Sigma\co = u_0$.  As discussed in the proof of Lemma \ref{LemmaFSigmaHat}, the map $F_\co^{-1}$ is analytic, which implies the analyticity of $E_{u_0}$.

With all our current function spaces being Sobolev spaces, we now apply Theorem \ref{ThmCR}, where so as to not confuse the value of $p$ in our present situation with that in Theorem \ref{ThmCR}, we let $p'$ denote what is $p$ in Theorem \ref{ThmCR}. So letting $X = L^2\tT_\Sigma$, $k=1$, $p' = 2$, and $q' = q = p > 4$, the hypotheses of Theorem \ref{ThmCR}(i) are satisfied, and we obtain
\begin{equation}
  (b,\phi,\xi) \in B^{2,2}(K, L^2\tT_\Sigma), \label{Omegareg}
\end{equation}
i.e., we have gained a whole derivative in the $K$ directions.  We can take $q' = q = p$, since
\begin{align}
  (b,\phi,\xi) & \in H^{1,p}(K, L^p\tT_\Sigma),
\end{align}
and $L^p\tT_\Sigma \subset X$, since $p > 4$. Furthermore, the elliptic estimate (\ref{CR-est}) implies
\begin{align}
  \|(b_0,\phi_0,\xi_0)\|_{B^{2,2}(K; L^2\tT_\Sigma)} & \leq \mu_{(A_0,\Phi_0)}\Big(\|D_\Sigma(b_0,\phi_0,\xi_0)\|_{B^{1,2}(K; L^2\tT_\Sigma)} + \|N_{(A_0,\Phi_0)}(b,\phi,\xi)\|_{B^{1,2}(K; L^2\tT_\Sigma)}\nonumber \\
   & \qquad + \|(b,\phi,\xi)\|_{B^{1,2}(K; L^2\tT_\Sigma)} + \|SW_4(A_0,\Phi_0)\|_{B^{1,2}(K; L^2\tT_\Sigma)}\Big)\nonumber \\
  & \leq \mu_{(A_0,\Phi_0)}\Big(\|(b,\phi,\xi)\|_{H^{1,p}(S^1 \times Y)} + \|SW_4(A_0,\Phi_0)\|_{B^{1,2}(K; L^2\tT_\Sigma)}\Big) \label{B22H1p}
\end{align}
for some nonlinear function $\mu_{(A_0,\Phi_0)}$.  Here $(b_0,\phi_0,\xi_0)$ plays the role of $u-u_0$ in (\ref{CR-est}). 

Summarizing, by using Theorem \ref{ThmCR} we deduced (\ref{B22H1p}) and gained regularity in the $K$ directions, i.e., we now have two derivatives in the $S^1 \times [0,1]$ directions in $L^2$.  Combined with the estimates from Step Two, where we gained regularity in just the $\Sigma$ directions, we see that we have gained a whole derivative in \textit{all} directions, i.e., $(b_0,\phi_0,\xi_0) \in B^{2,2}\tT_{[0,1]\times\Sigma}$.  Combined with interior regularity, altogether we have the elliptic estimate
\begin{equation}
  \|(b,\phi,\xi)\|_{B^{2,2}(S^1 \times Y)} \leq \mu_{(A_0,\Phi_0)}\Big(\|(b,\phi,\xi)\|_{H^{1,p}(S^1 \times Y)}\| + \|SW_4(A_0,\Phi_0)\|_{H^{1,p}(S^1 \times Y)}\Big). \label{xB22}
\end{equation}
on all of $S^1 \times Y$.\\

\textit{Step Four:} From the previous steps, our configuration $(A,\Phi) \in H^{1,p}(S^1 \times Y)$, which we redefined by a gauge transformation so that it is in Coulomb-Neumann gauge with respect to $(A_0,\Phi_0)$, is in $B^{2,2}(S^1\times Y)$.  Moreover, we have the elliptic estimate (\ref{xB22}).  Proceeding as in the previous two steps, we want to bootstrap and show that $(b,\phi,\xi) = (A,\Phi) - (A_0,\Phi_0)$ is in $B^{k,2}(S^1 \times Y)$ for all $k \geq 2$, which will prove the theorem.  Unfortunately, in our first step when we want to bootstrap from $B^{2,2}(S^1 \times Y)$ to $B^{3,2}(S^1 \times Y)$,  the space $B^{2,2}(S^1 \times Y)$ is not strictly stronger than the original space $H^{1,p}(S^1 \times Y)$, i.e., we do not have an embedding $B^{2,2}(S^1 \times Y) \hookrightarrow H^{1,p}(S^1 \times Y)$, since $p > 4$.  Thus, we will need to work with the mixed topology $H^{1,p}(S^1 \times Y) \cap B^{2,2}(S^1 \times Y)$.  This is the cause for the rather bizarre looking Thereom \ref{ThmPathL2}.

We first start off by increasing the regularity of $E^1_{\g_0}(z)$ in (\ref{gz}).  Indeed, since we now have $(A,\Phi) \in H^{1,p}(S^1 \times Y) \cap B^{2,2}(S^1 \times Y$), then $\g \in \Maps^{1-1/p,p}(S^1,\fL) \cap \Maps^{3/2,2}(S^1,\fL)$.  Thus, we can write
\begin{equation}
  \g = \g_0+z+\E^1_{\g_0}(z), \qquad z \in T_{\g_0}\Maps^{1-1/p,p}(S^1,\fL) \cap \Maps^{3/2,2}(S^1,\fL), \label{gz2}
\end{equation}
and by Theorem \ref{ThmPathL2} with $s_2 = 0$, we have
\begin{equation}
  \E^1_{\g_0}(z) \in \Maps^{(1-1/p,1-1/p-\eps),p}(S^1,\T_\Sigma) \cap \Maps^{(3/2,1/2),2}(S^1,\T_\Sigma), \label{eq:step4-1}
\end{equation}
Here, we have assumed that $\g_0$ is sufficiently close to $\g$ in $C^0(I, B^{s',p})$, using the \textsl{same} $s'$ we used in Step Two when we applied Theorem \ref{ThmPathL1} to (\ref{gz}), so that we can place $\g$ in the chart map (\ref{gz2}) using Theorem \ref{ThmPathL2}, which is stronger than (\ref{gz}).  (We could of course redefine $\g_0$ at this step by moving it closer to $\g$ if necessary.) Here, to perform the above step, it is crucial that in both these theorems, the size of the chart maps (the radius $\delta$ which appears) depends locally uniformly with respect to the $B^{s',p}$ topology, which is very weak since we can choose any $1/2 < s' < s - 2/p$ (since $s_2 = 0$).

Focusing on the second factor of (\ref{eq:step4-1}), then $\E^1_{\g_1}(z)$ smooths by $1/2$ a derivative in the $\Sigma$ directions.  Plugging this into (\ref{Sigmaboot}) and proceeding as before, we find
\begin{align}
  \|(b,\phi,\xi)\|_{B^{(2,1/2),2}((S^1\times [0,1])\times \Sigma)} & \leq \mu_{(A_0,\Phi_0)}\Big(\|(b,\phi,\xi)\|_{B^{2,2}(S^1 \times Y)}+\|(b,\phi,\xi)\|_{H^{1,p}(S^1 \times Y)}\nonumber\\
  & \hspace{1in} \|SW_4(A_0,\Phi_0)\|_{B^{2,2}(S^1\times Y)}\Big), \nonumber \\
  & \leq \mu_{(A_0,\Phi_0)}\Big(\|(b,\phi,\xi)\|_{H^{1,p}(S^1 \times Y)} + \nonumber \\
  & \hspace{1in} \|SW_4(A_1,\Phi_1)\|_{B^{2,2}(S^1\times Y)}\Big), \label{Sigmaboot1/2}
\end{align}
where we use that $\|(b,\phi,\xi)\|_{B^{2,2}(S^1 \times Y)}$ is controlled by $\|(b,\phi,\xi)\|_{H^{1,p}(S^1\times Y)}$ via (\ref{xB22}).  

Estimate (\ref{Sigmaboot1/2}) is insufficient however since we want to gain a full derivative in the $\Sigma$ direction, i.e., we want control of $B^{(2,1),2}((S^1 \times [0,1]) \times \Sigma)$ instead of $B^{(2,1/2),2}((S^1 \times [0,1])\times\Sigma)$.  Thus, we repeat the above steps again, where we replace $\Maps^{3/2,2}(S^1,\fL)$ with $\Maps^{(3/2,s_2),2}(S^1,\fL)$, with $s_2 = 1/2$. Using the same reasoning as before (and assuming $(A_0,\Phi_0)$ is sufficiently $C^0(S^1,B^{s',p}(\Sigma))$ close to $(A,\Phi)$ on the boundary), by Theorem \ref{ThmPathL2}, we obtain
\begin{equation}
  \E^1_{\g_0}(z) \in \Maps^{(1-1/p,1-1/p-\eps),p}(S^1,\T_\Sigma) \cap \Maps^{(3/2,1),2}(S^1,\T_\Sigma),
\end{equation}
thereby improving the gain in $\Sigma$ regularity from $1/2$ to $1$.  Doing this, we now gain a whole derivative in the $\Sigma$ directions:
\begin{align}
  \|(b,\phi,\xi)\|_{B^{(2,1),2}((S^1\times [0,1])\times \Sigma)} & \leq \mu_{(A_0,\Phi_0)}\Big(\|(b,\phi,\xi)\|_{H^{1,p}(S^1 \times Y)} + \nonumber\\
  & \hspace{1in}  \|SW_4(A_0,\Phi_0)\|_{B^{2,2}(S^1\times Y)}\Big). \label{Sigmaboot1}
\end{align}

Having gained a whole derivative in the $\Sigma$ directions, we can proceed to Step Three and gain regularity in the $S^1 \times [0,1]$ directions. Here, we need to choose our parameters in Theorem \ref{ThmCR} appropriately.  In the same way that we needed to proceed in two steps to gain a whole derivative in the $\Sigma$ directions, we will also need to proceed in two steps to gain a whole derivative in the $S^1\times[0,1]$ directions as well.  First, we let $k = 1$.  By the above, we have $x = (b,\phi,\xi)$ belongs to the space $H^{(2,1),2}(K \times \Sigma)$ in addition to belonging to the space $H^{1,p}(K \times \Sigma)$.  Thus, $D_\Sigma \in H^{2,2}(K \times \Sigma) \subset H^{1,4}(K \times \Sigma)$.  Using the multiplication theorem Theorem \ref{ThmMult}, we have that $H^{(2,1),2}(K \times \Sigma) \cap L^\infty$ is an algebra.  In particular, $x\# x \in H^{(2,1),2}(K \times \Sigma) \subset H^{1,4}(K \times \Sigma)$, where $\#$ denote any pointwise multiplication map.  It follows that $x$ satisfies a Cauchy-Riemann equation (\ref{CR-CR}), where $G$ has the same regularity as $D_\Sigma x + x\# x \in H^{1,4}(K \times \Sigma)$.  Thus, thinking of $G$ as taking values in
$$X := L^{p_0}\tT_\Sigma$$
for some $p_0 = 2 + \eps$, where $\eps > 0$ is small, we have
\begin{equation}
  G \in H^{1,4}(K, X).
\end{equation}
Thus, we may apply Theorem \ref{ThmCR} with $p' = p_0$, $q = 4$, and $q' = p > 4$ (recall that $p'$ is the value of the dummy variable ``$p$" in Theorem \ref{ThmCR}, to distinguish it from our present value of $p$).  Note that when we change from the Banach space $L^2\tT_\Sigma$ to $X = L^{p_0}\tT_\Sigma$, we must also consider for the Lagrangian boundary values of $x$ the locally embedded $L^{p_0}$ charts associated with $\fL$, instead of $L^2$ charts as before.  This is possible since $x \in C^0(K, L^p\tT_\Sigma)$, and $L^p\tT_\Sigma \subset L^{p_0}\tT_\Sigma$ since $p > p_0$.  In any event, we apply Theorem \ref{ThmCR} and obtain $x \in H^{2,p_0}(K,X)$.  On the other hand, since
$$x \in H^{(2,1),2}(K \times \Sigma) \subset H^{2,4}(\Sigma, L^4(K)) \subset H^{2,p_0}(\Sigma,L^{p_0}(K)),$$
it follows that $x \in H^{2,p_0}(K \times \Sigma)$.

This implies we have improved the integrability of $(A,\Phi)$ from $H^{1,p}(S^1\times Y) \cap H^{2,2}(S^1\times Y)$ to $H^{1,p}(S^1\times Y) \cap H^{2,p_0}(S^1\times Y)$, with $p_0 > 2$. This extra integrability now allows us to increase the regularity of $x$ by applying Theorem \ref{ThmCR} with $k=2$.  Here, we let $p' = q=2$, and $q' = p_0$.  Observe that $q' = 2$ does not work, which is why we needed the above step.  Doing this gives us a configuration in $H^{3,2}(S^1 \times Y) = B^{3,2}(S^1 \times Y)$ which is strictly stronger than $H^{1,p}(S^1 \times Y)$ for $p$ close to $4$.

We can now continue bootstrapping as above, using Theorem \ref{ThmPathL2} and estimate (\ref{inhom-est}) as above to gain $\Sigma$ regularity, and then Theorem \ref{ThmCR} to gain $S^1\times[0,1]$ regularity.  Each time, we apply Theorem \ref{ThmPathL2} to gain a full derivative in the $\Sigma$ directions, and then we apply Theorem \ref{ThmCR} once to gain a whole derivative in the $S^1\times[0,1]$ directions.  Indeed, our function spaces are now sufficiently regular that we can apply Theorem \ref{ThmPathL2} to gain a whole derivative in the $\Sigma$ directions ($s > 3/2$ so $\eps' = 0$ in the theorem), and we can apply Theorem \ref{ThmCR} to gain one whole derivative in the $S^1\times[0,1]$ directions in one step without having to first bootstrap the integrability of our configuration as in the above. Together, these steps gain for us a whole derivative in all directions.

Altogether, we have shown the following.  Pick any smooth reference connection $A_\rf$ and redefine $(A,\Phi)$ by a gauge transformation that places $A$ in Coulomb-Neumann gauge with respect to $A_\rf$.  Then finding smooth $(A_0,\Phi_0)$, satisfying the Lagrangian boundary conditions, that is sufficiently $H^{s,p}(S^1\times Y)$ close to $(A,\Phi)$, with $s > 1/2+2/p$, then for every $k \geq 2$, we have the estimate
\begin{align}
  \|(A - A_0,\Phi - \Phi_0)\|_{B^{k,2}(S^1 \times Y)} & \leq \mu_{k, (A_0,\Phi_0)}\Big(\|(A-A_0,\Phi-\Phi_0)\|_{H^{1,p}(S^1 \times Y)} + \nonumber\\
  & \hspace{1in} \|SW(A_0,\Phi_0)\|_{B^{k-1,2}(S^1\times Y)}\Big), \label{ellboot}
\end{align}
where $\mu_k = \mu_{k, (A_0,\Phi_0)}$ is a continuous nonlinear function depending on $k$ and $(A_0,\Phi_0)$.  This estimate proves the theorem.\End

There is no obstacle to extending the above result to the equations on $\R \times Y$:\\

\textit{Proof of Theorem A: } It is enough to show that for every compact interval $I \subset \R$, we can find a gauge transformation belonging to the identity component of $H^{2,p}\mathcal{G}(I \times Y)$ that makes $(A,\Phi)$ smooth on $I \times Y$. The result then follows by writing $\R \times Y$ as a union $\cup_{k \in \Z} I_k \times Y$, with $I_k = [k/2, k/2+1]$, and then patching together the gauge transformations on the overlaps.

So fix one of the intervals $I_k$, call it $I$. Embed $I$ inside the interior of a larger closed interval $\tilde I$ and regard the latter as being embedded into $S^1$. We will extend $(A,\Phi)|_{I \times Y}$ to $S^1 \times Y$, with $(A,\Phi)$ ``almost" satisfying Lagrangian boundary conditions and smooth outside of $\tilde I \times Y$. We do this as follows. Consider the boundary value configuration of $(A,\Phi)|_{\tilde I \times Y}$:
$$\gamma(t) = r_\Sigma(B(t),\Phi(t)) \in \Maps^{1-1/p,p}(\tilde I, \fL).$$
Pick a smooth path $\gamma_0 \in \Maps^{1-1/p,p}(\tilde I, \fL)$ such that $\gamma$ lies in a chart map centered at $\gamma_0$:
$$\gamma = \gamma_0 + z + \cE^1_{\gamma_0}(z), z \in T_{\gamma_0}\Maps^{1-1/p,p}(\tilde I,\fL).$$
We can extend $\gamma_0$ to a smooth element in $\Maps(S^1,\fL)$, again denoted $\gamma_0$, and then extend such a configuration smoothly into $S^1 \times Y$, obtaining $(A_0, \Phi_0)$ satisfying Lagrangian boundary conditions. Let $\psi: S^1 \to \R$ be a bump function identically one on $I$ and vanishing outside of $\tilde I$; this extends to a bump function on $S^1 \times Y \to \R$ that is constant in the $Y$-directions. Define an extension $(\tilde A, \tilde \Phi)$ of $(A,\Phi)|_{I \times Y}$ to $S^1 \times Y$ via
\begin{equation}
  (\tilde A, \tilde \Phi) = (A_0,\Phi_0) + (\psi(A - A_0), \psi(\Phi - \Phi_0))
\end{equation}
So $(\tilde A,\tilde \Phi) = (A_0,\Phi_0)$ outside of $\tilde I$. Moreover, on $S^1 \times Y$, we have
\begin{equation}
  SW_4(\tilde A, \tilde \Phi) = SW_4(A_0,\Phi_0) + P(A - A_0, \Phi-\Phi_0) \label{eqA:SW4}
\end{equation}
where $P(A-A_0,\Phi-\Phi_0)$ is a quadratic polynomial in $(A-A_0,\Phi-\Phi_0)$ (involving only pointwise multiplication). Indeed, this is because $(A,\Phi)$ satisfies the Seiberg-Witten equations on $\R \times Y$. It follows that for $(A,\Phi)$ in $H^{1,p}$, the right-hand side of (\ref{eqA:SW4}) belongs to $H^{1,p}(S^1 \times Y)$.
In what follows, we check that the boostrapping estimates of Theorem \ref{ThmA'} can be carried over to gain regularity for $(A,\Phi)|_{I \times Y}$ using the configuration $(\tilde A,\tilde \Phi)$ (after choosing a gauge which places the latter in Coulomb-Neumann gauge on $S^1 \times Y$). First, if we want to bootstrap the regularity of $(\tilde A, \tilde \Phi)$ on $S^1 \times Y$, we replace $-SW_4(A_0,\Phi_0)$ in (\ref{SWeqN1}) with $SW_4(\tilde A, \tilde \Phi) - SW_4(A_0,\Phi_0)$. This just follows from how equation (\ref{SWeqN1}) is defined, which expresses  $SW_4(\tilde A, \tilde \Phi) - SW_4(A_0,\Phi_0)$ in terms of a linear and quadratic functional of $(\tilde A - A_0, \tilde \Phi - \Phi_0)$. Because $SW_4(A_0,\Phi_0)$ is smooth and $SW_4(\tilde A, \tilde \Phi)$ always has the same regularity as $(A,\Phi)$ by (\ref{eqA:SW4}) (instead of one derivative less) the inhomogeneous term $SW_4(\tilde A, \tilde \Phi) - SW_4(A_0,\Phi_0)$ causes no problems with (interior) bootstrapping. The boundary value of $(\tilde A,\tilde \Phi)$ is given by
$$\gamma_0 + \psi z + \psi \cE_{\gamma_0}^1(z).$$
While not satisfying the Lagrangian boundary condition, we can still bootstrap the regularity of $(\tilde A,\tilde \Phi)$ in the $\Sigma$-directions near the boundary using Step Two in the proof of Theorem \ref{ThmA'}. This is because in that step, the precise nature of the error term $\psi \cE_{\gamma_0}^1(z)$ is immaterial; all that is required is that it is smoother than $\gamma_0 + \psi z$ lying inside the path space of the Lagrangian linearized about $\gamma_0$.

Next, in Step 3 of Theorem \ref{ThmA'}, we must restrict ourselves back to $I \times Y$, where we use the Banach-space valued Cauchy-Riemann equations which crucially needed Lagrangian boundary conditions on the nose (we take the domain of the Cauchy-Riemann equations to be $K = I \times [0,1]$, where $[0,1]$ is the normal direction to the boundary, in Theorem \ref{ThmCR}). But having gained regularity in the $\Sigma$-directions from Step 2, Step 3 goes through as before.

Finally, each time we bootstrap regularity in the final step of the proof of Theorem \ref{ThmA'}, we need to shrink the interval $I$ slightly (by an arbitrarily small amount) since we need the support of the bump function $\psi$ defining the extension $(\tilde A, \tilde \Phi)$ to be where $(A,\Phi)$ has gained regularity. Nevertheless, since our initial choice of $I$ was an arbitrary compact interval, the end result is that for any fixed compact interval $I$, after choosing a gauge that puts the extended configuration $(\tilde A,\tilde \Phi)$ in Coluomb-Neumann gauge on $S^1 \times Y$, we have that $(A,\Psi)|_{I \times Y}$ is smooth. More precisely, we established the following elliptic estimate
\begin{align}
  \|(A - A_0,\Phi - \Phi_0)\|_{B^{k,2}(I \times Y)} & \leq \mu_{k, (A_0,\Phi_0)}\Big(\|(A-A_0,\Phi-\Phi_0)\|_{H^{1,p}(I \times Y)} + \nonumber\\
  & \hspace{1in} \|P(A-A_0,\Phi-\Phi_0)\|_{B^{k-1,2}(I\times Y)}\Big), \label{ellboot2}
\end{align}
generalizing (\ref{ellboot}).\End

Next, we prove the analog of Theorem B in the periodic setting:

\begin{Theorem}\label{ThmB'} Let $p > 4$ and let $(A_i,\Phi_i) \in H^{1,p}\fC(S^1 \times Y)$ be a sequence of solutions to (\ref{SWeq1})-(\ref{SWeq3}), where $\fL$ is a fully gauge invariant monopole Lagrangian. More generally, we can suppose the right-hand side of (\ref{SWeq1}) is instead a quadratic pointwise-multiplication polynomial in $(A,\Phi)$. Suppose we have uniform bounds
\begin{equation}
  \|F_{A_i}\|_{L^p(S^1\times Y)}, \|\nabla_{A_i}\Phi_i\|_{L^p(S^1\times Y)}, \|\Phi_i\|_{L^p(S^1\times Y)} \leq C \label{energyboundS^1}
\end{equation}
for some constant $C$. Then there exists a subsequence of configurations, again denoted by $(A_i,\Phi_i)$, and a sequence of gauge transformations $g_i \in H^{2,p}\G(S^1\times Y)$ such that $g_i^*(A_i,\Phi_i)$ converges uniformly in $C^\infty(S^1\times Y)$.
\end{Theorem}

\Proof Fix any smooth reference connection $A_\rf$ and redefine the $(A_i,\Phi_i)$ by gauge transformations $g_i$ that place $A_i$ in Coulomb-Neumann gauge with respect to $A_\rf$.  The elliptic estimate for $d+d^*$ on $1$-forms with Neumann boundary condition implies that
\begin{equation}
  \|A_i-A_\rf\|_{H^{1,p}} \leq c(\|F_{A_i}\|_{L^p} + \|(A_i-A_\rf)^h\|_{L^p}), \label{controlA}
\end{equation}
where $(A_i-A_\rf)^h$ is the the orthogonal projection of $(A_i-A_\rf)$ onto the finite dimensional subspace
\begin{equation}
  \{a \in \Omega^1(S^1\times Y; i\R) : da = d^*a = 0, *a|_{S^1\times\Sigma}=0\} \cong H^1(Y; i\R).
\end{equation}
The above isomorphism is by the usual Hodge theory on manifolds with boundary.  From just the bounds (\ref{energyboundS^1}), we have no a priori control of $ \|(A_i-A_\rf)^h\|_{L^p}$.  However, we still have some gauge freedom left, namely, we can consider the following group of harmonic gauge transformations
\begin{equation}
  \G_{h,n} := \{g \in \G(S^1 \times Y) : d^*(g^{-1}dg) = 0, *dg|_{S^1\times \Sigma} = 0\}
\end{equation}
which preserve the Coulomb-Neumann gauge.  The map $g \mapsto g^{-1}dg$ maps $\G_{h,n}$ onto the lattice $H^1(Y; 2\pi i\Z)$ inside $H^1(Y;i\R)$.  Hence, by redefining the $A_i$ by gauge transformations in $\G_{h,n}$, we can arrange that the $(A_i-A_\rf)^h$ are bounded uniformly, which together with (\ref{energyboundS^1}) and  (\ref{controlA}) implies that we have a uniform bound
\begin{equation}
  \|A_i-A_\rf\|_{H^{1,p}} \leq C \label{controlA2}
\end{equation}
for some absolute constant $C$ (where $C$ denotes some constant independent of the $(A_i,\Phi_i)$, whose value may change from line to line).

From (\ref{energyboundS^1}) and (\ref{controlA2}), we have the control
\begin{align}
  \|\nabla_{A_\rf}\Phi_i\|_{L^p} & \leq \|\nabla_{A_i}\Phi_i\|_{L^p} + \|\rho(A_i-A_\rf)\Phi_i\|_{L^p} \nonumber \\
  & \leq \|\nabla_{A_i}\Phi_i\|_{L^p} + \|\rho(A_i-A_\rf)\|_{L^\infty}\|\Phi_i\|_{L^p} \nonumber \\
  & \leq C, \label{controlPhi}
\end{align}
due to the embedding $H^{1,p}(S^1\times Y) \hookrightarrow L^\infty(S^1\times Y)$ for $p > 4$.  The uniform bound (\ref{controlPhi}) and the uniform bound on $\|\Phi_i\|_{L^p}$ shows that we have the uniform bound
\begin{equation}
  \|\Phi_i\|_{H^{1,p}} \leq C.
\end{equation}
Thus, the configuration $(A_i,\Phi_i)$ is uniformly bounded in $H^{1,p}(S^1\times Y)$. Moreover, the $(A_i,\Phi_i)$ are smooth since they solve (\ref{SWeq1})-(\ref{SWeq3}) and $A_i$ is in Coulomb-Neumann gauge with respect to a smooth connection.  If we can show that the $(A_i,\Phi_i)$ are also uniformly bounded in $H^{k,2}(S^1\times Y)$ for each $k \geq 2$, then we will be done, due to the compact embedding $H^{k+1,2}(S^1 \times Y) \hookrightarrow H^{k,2}(S^1 \times Y)$ for all $k \geq 1$ and a diagonalization argument.

Since the $(A_i,\Phi_i)$ are uniformly bounded in $H^{1,p}(S^1\times Y)$, a subsequence converges strongly in $H^{s,p}(S^1\times Y)$ for any $s = 1-\eps$ with $\eps > 0$ arbitrarily small.  The limiting configuration $(A_\infty,\Phi_\infty)$, being a weak $H^{1,p}(S^1 \times Y)$ limit of the $(A_i,\Phi_i)$, belongs to $H^{1,p}(S^1 \times Y)$, and it solves (\ref{BVP}), since the equations are preserved under weak limits.  In the interior, this is easy to see; on the boundary, we use the fact that $\Maps^{s-1/p,p}(S^1,\fL)$ is a manifold, so that the Lagrangian boundary condition is preserved under weak limits.  Since Coulomb-Neumann gauge is also preserve under weak limits, then from Theorem \ref{ThmA'}, we know that $(A_\infty,\Phi_\infty)$ is smooth.

We now apply (\ref{ellboot}) with $(A_0,\Phi_0)$ replaced with the smooth configuration $(A_\infty,\Phi_\infty)$ and $(A,\Phi)$ replaced by the $(A_i,\Phi_i)$, for large $i$; or else in the more general situation where we have a quadratic map $Q(A,\Phi)$ on the right-hand side of (\ref{SWeq1}), we apply (\ref{ellboot}) with $SW_4(A_0,\Phi_0)$ replaced with $Q(A_i,\Phi_i) - SW_4(A_\infty,\Phi_\infty)$.  We can do this because the following are true: first, the $(A_i,\Phi_i)$ converge strongly to $(A_\infty,\Phi_\infty)$ in $H^{s,p}(S^1 \times Y)$, and $s > 1/2 + 2/p$; second, the proof of (\ref{ellboot}) shows that if $(A_0,\Phi_0)$ is any smooth configuration, then (\ref{ellboot}) holds for all $(A,\Phi)$ solving (\ref{SWeq1})--(\ref{SWeq3}) sufficiently $H^{s,p}(S^1 \times Y)$ close to $(A_0,\Phi_0)$.  It now follows that a subsequence of the $(A_i,\Phi_i)$ converges to $(A_\infty,\Phi_\infty)$ in $C^\infty$.\End

\textit{Proof of Theorem B: } We proceed as in the proof of Theorem A. Write $\R \times Y$ as a union $\cup_{l \in \Z} I_k \times Y$ with $I_k = [k/2,k/2+1]$. On each $I_k$, we can find gauge transformations $g_{k,i}$ on $I_k \times Y$ such that $g_{k,i}^*(A_i,\Phi_i)$ belong to $C^\infty(I_k \times Y)$. Moreover, with $k$ fixed, the $g_{k,i}^*(A_i,\Phi_i)$ are uniformly bounded since we can bound the extensions $(\tilde A_i,\tilde \Phi_i)$ on $S^1 \times Y \supset I_k \times Y$ in terms of (i) the smooth reference configuration $(A_0,\Phi_0)$ on $S^1 \times Y$ used to define the extensions and (ii) the uniform bounds on $\|F_{A_i}\|_{L^p}$, $\|\nabla_{A_i}\Phi_i\|_{L^p}$, and $\|\Phi_i\|_{L^p}$ on a sufficiently large compact subset of $\R \times Y$. (That an arbitrary extension $(A_0,\Phi_0)$ was used in (i) is immaterial; what matters is that it is independent of $i$.) This follows from Theorem \ref{ThmB'} and the modifications explained in the proof of Theorem A to accommodate for the extensions $(\tilde A_i, \tilde \Phi_i)$ being only ``approximate" solutions to the Seiberg-Witten equations with Lagrangian boundary conditions on $S^1 \times Y$.

So on the overlaps of the $I_k$, the transition maps $g_{k+1,i}^{-1}g_{k,i}$ are bounded in $C^\infty((I_k \cap I_{k+1}) \times Y)$ uniformly in $i$. Thus, passing to some subsequence in $i$, the homotopy class of  $g_{k+1,i}^{-1}g_{k,i}$ is constant. Without loss of generality, we can assume these maps are homotopic to the identity by modifying the $g_{k,i}$ or $g_{k+1,i}$ by some fixed gauge transformation. One can thus patch together the $g_{k,i}$'s for $k \in \Z$ and obtain a single gauge transformation $g_i$ on $\R \times Y$ such that $g_i^*(A_i,\Phi_i)|_K$  are uniformly bounded in $C^\infty(K)$ for any compact subset $K \subset \R \times Y$ (the bound depends on $K$). Thus, a subsequence of the $g_i^*(A_i,\Phi_i)$ converges in $C^\infty(K)$. Exhausting $\R \times Y$ by a sequence of compact subsets $K_j$ and taking a diagonal subsequence in $i$ and $j$ proves Theorem B.\End


\appendix

\section{Properties of Anisotropic Function Spaces}

Here, we collect some properties about anisotropic Besov spaces\footnote{Everything in this section also applies to anisotropic Bessel potential spaces, where Bessel potential spaces were defined as in the appendix of \cite{N1}. See also \cite{N3}.}.  All of these results are standard when we specialize to usual (isotropic) Besov spaces, i.e., when the anisotropic parameter is zero.  The goal of this section is to see how these standard isotropic results can be strengthened when there is anisotropy.  We will make some remarks to give the reader some intuition for the anisotropic spaces, but we will defer most proofs to \cite{N0} and \cite{N3}.

The first result is the extension of the usual trace and extension theorems \cite[Theorem C.3, C.5]{N1} to the anisotropic setting.  Recall that if $X$ is a compact manifold with boundary, then we have a (zeroth order) trace map $r: B^{s,p}(X) \to B^{s-1/p,p}(\partial X)$, $s > 1/p$, which ``costs" us $1/p$ derivatives and an extension map $B^{s-1/p,p}(\partial X) \to B^{s,p}(X)$ which gains us $1/p$ derivatives. Suppose now we have $X = X_1 \times X_2$, where $X_1$ and $X_2$ are compact manifolds (with or without boundary).  When taking a trace to the boundary, the anisotropy of a function can be either tangential or normal to the boundary. When the anisotropy of a function is tangential to the boundary, then the trace and extension operators preserve this anisotropy, since tangential operations commute with such operators (in the model Euclidean space when the boundary is a hyperplane).  On the other hand, if there is anisotropy in the normal direction, then if there is enough anisotropy, taking a trace costs us $1/p$ derivatives only in the anisotropic directions, when $p \geq 2$.  This is summarized in the following:
\begin{Theorem}\label{ThmATrace}(Anisotropic Traces and Extensions)
  Let $X = X_1 \times X_2$.
   \begin{enumerate}
     \item (Tangential anisotropy) Suppose $X_1$ has boundary $\partial X_1$.  Then for $s_1 > 1/p$ and $s_2 \geq 0$, the trace map satisfies
\begin{equation}
  r: B^{(s_1,s_2),p}(X_1 \times X_2) \to B^{(s_1-1/p,s_2),p}(\partial X_1 \times X_2).
\end{equation}
     Furthermore, for all $s_1 \in \R$, there exists a boundary extension map
     $$e: B^{(s_1-1/p,s_2),p}(\partial X_1 \times X_2) \to B^{(s_1,s_2),p}(X_1 \times X_2),$$
     and for $s_1 > 1/p$, we have $re = \id$.  Moreover, let $\tilde X_1$ be any closed manifold extending $X_1$.  Then for every $k \in \N$, we have an extension map
   $$E_k: B^{(s_1,s_2),p}(X_1\times X_2) \to B^{(s_1,s_2),p}(\tilde X_1\times X_2),$$
   for $|s_1| < k$.
     \item (Mixed anisotropy) Suppose $X_2$ has boundary $\partial X_2$.  Then for $s_1 \geq 0$, $s_2 > 1/p$, and $p\geq 2$, the trace map satisfies
     \begin{equation}
       r: B^{(s_1,s_2),p}(X_1 \times X_2) \to B^{(s_1, s_2-1/p-\eps),p}(X_1 \times \partial X_2),
     \end{equation}
         where $\eps > 0$ is arbitrary.  If $p = 2$, we can take $\eps = 0$.
   \end{enumerate}
\end{Theorem}
The above theorem readily generalizes to higher order trace maps.

Next, we recall that for isotropic Besov spaces, we have the embedding
\begin{equation}
  B^{s,p}(X) \hookrightarrow C^0(X) \label{C0embed}
\end{equation}
if $s > n/p$, where $n = \dim X$.  Thus, we have the following corollary:

\begin{Corollary}\label{CorEmbed}
  Let $X = X_1 \times X_2$ and $n_i = \dim X_2$.  Then if  $s_1 > n_1/p$ and $s_2 > n_2/p$, we have $B^{(s_1,s_2),p}(X_1 \times X_2) \hookrightarrow C^0(X_1 \times X_2)$.
\end{Corollary}

\Proof By (\ref{C0embed}) and Theorem \ref{ThmATrace}, we can take successive traces to conclude $B^{(s_1,s_2),p}(X_1 \times X_2) \hookrightarrow C^0(X_1; B^{s_1 + s_2 - n_1/p,p}(X_2))$ since $s_1 > n_1/p$.  By (\ref{C0embed}), we have $B^{s_1+s_2-n_1/p,p}(X_2) \hookrightarrow C^0(X_2)$ since $s_1+s_2-n_1/p > n_2/p$, whence the theorem follows.\End

Recall that we have a multiplication theorem for Besov spaces.  Such a theorem is proved using the paraproduct calculus.  By redoing this for anisotropic Besov spaces, one can also prove an anisotropic muliplication theorem.  We state one for $p = 2$.  This is all summarized in the below:

\begin{Theorem} (Multiplication Theorem) \label{ThmMult}
    \begin{enumerate}
     \item For all $s > 0$, we have $B^{s,p}(X) \cap L^\infty(X)$ is an algebra. Moreover, we have the estimate
     $$\|fg\|_{B^{s,p}} \leq C(\|f\|_{B^{s,p}}\|g\|_{L^\infty} + \|f\|_{L^{\infty}}\|g\|_{B^{s,p}}).$$
     In particular, if $s > n/p$, then $B^{s,p}(X)$ is an algebra.
     \item Let $s_1 \leq s_2$ and suppose $s_1 + s_2 > n \max(0,\frac{2}{p}-1)$. Then we have a continuous multiplication map
     \begin{align*}
       B^{s_1,p}(X) \times B^{s_2,p}(X) & \to B^{s_3,p}(X),
     \end{align*}
    where
    $$s_3 = \begin{cases}
      s_1 & \mathrm{if }\;\, s_2 > n/p\\
      s_1+s_2-n/p & \mathrm{if }\;\, s_2 < n/p.
    \end{cases}$$
    \item (Anisotropic Multiplication) Let $\dim X_i = n_i$, $i=1,2$ and suppose $s_1 > n_1/2$. Let $s_2', s_2'' \geq 0$ and let $s_2 \leq \min(s_2',s_2'')$ satisfy $s_2 < s_1+s_2'+s_2'' - \frac{n_1+n_2}{2}$.  Then we have a multiplication map
  \begin{equation}
    [B^{(s_1,s_2'),2}(X_1 \times X_2) \cap L^\infty] \times [B^{(s_1,s_2''),2}(X_1 \times X_2) \cap L^\infty] \to B^{(s_1,\max(s_2,0)),2}(X_1 \times X_2).
  \end{equation}
   \end{enumerate}
\end{Theorem}

We also have the following important Fubini property, which says that a function is in a Besov space on $X_1\times X_2$ if only if it lies in a Besov space in each of the $X_1$ and $X_2$ directions.  More precisely, we have the following theorem:

\begin{Theorem}\cite[Theorem 2.5.13]{Tr1} \label{ThmFubini}  (Fubini Property) For any $s > 0$, we have
$$B^{s,p}(X_1 \times X_2) = L^p(X_1, B^{s,p}(X_2)) \cap L^p(X_2, B^{s,p}(X_1)).$$
\end{Theorem}

\begin{Theorem}\label{ThmPSDO}
   Pseudodifferential operators are bounded on $B^{(s_1,s_2),p}(X_1 \times X_2)$, for all $1 < p <\infty$, $s_1 \in \R$, and $s_2 \geq 0$.
\end{Theorem}

\section{Elliptic Boundary Value Problems Revisited}

We want to generalize the general pseudodifferential setup of elliptic boundary value problems, as done in Appendix C of \cite{N1}, to more general boundary conditions and also to anisotropic function spaces.  Let $A: \Gamma(E) \to \Gamma(F)$ be an elliptic differential operator acting between sections of vector bundles $E$ and $F$ over a compact manifold $X$ with boundary $\partial X=\Sigma$.  For simplicity, let us take $A$ to be first order (the only case we will need in this paper), though what follows easily generalizes to elliptic operators of any order.  Observe then a boundary condition for $A$ is simply a choice subspace of $\U \subset \Gamma(E_\Sigma)$ of the boundary data space, where $E_\Sigma = E|_\Sigma$.  The desirable boundary conditions are those for which the operator
\begin{equation}
  A_\U: \{x \in \Gamma(E) : r(x) \in \U\} \to \Gamma(F) \label{AU}
\end{equation}
is a Fredholm operator in the appropriate function space topologies.  In \cite{N1} and in typical situations, the subspaces $\U$ are given by the range of pseudodifferential projections.  However, from the above viewpoint, one need only consider the functional analytic setup of subspaces and appropriate function space topologies in order to understand the operator (\ref{AU}).

Let $X$ be a manifold with boundary, and suppose it can be written as a product $X = X_1 \times X_2$, where $X_1$ is a manifold with boundary and $X_2$ is closed.  We have the anisotropic Besov spaces $B^{(s_1,s_2),p}(X_1\times X_2)$ on $X$ and $B^{(s_1,s_2),p}(\partial X_1\times X_2)$ on $\partial X$, for $s_1 \geq 1$, $s_2 \geq 0$, and $1<p<\infty$. These spaces induce topologies on vector bundles, and so in particular, we have the spaces $B^{(s_1,s_2),p}(E)$, $B^{(s_1,s_2),p}(F)$ and such.  We have a restriction map
$$r: B^{(s_1,s_2),p}(E) \to B^{(s_1-1/p,s_2)}(E_{\partial X_1 \times X_2}),$$
and given a subspace
$$\U \subset B^{(s_1-1/p,s_2)}(E_{\partial X_1 \times X_2}),$$
we get the space
$$B^{(s_1,s_2),p}_\U(E) =  \{x \in B^{(s_1,s_2),p}(E) : r(x) \in \U\}.$$
Thus, (\ref{AU}) yields the operator
\begin{equation}
  A_\U: B^{(s_1,s_2),p}_\U(E) \to B^{(s_1-1,s_2),p}(F). \label{AU1}
\end{equation}
We can now ask how properties of $\U$ correspond with properties of the induced map (\ref{AU1}).

For this, we must distinguish a special subspace of $B^{(s_1-1/p,s_2)}(E_{\partial X_1 \times X_2})$, namely $\im P^+$, where
$$P^+: B^{(s_1-1/p,s_2)}(E_{\partial X_1 \times X_2}) \to B^{(s_1-1/p,s_2)}(E_{\partial X_1 \times X_2})$$
is the Calderon projection of $A$ (see Definition I.D.3).  This is a projection onto $r(\ker A)$, the boundary values of the kernel of $A$.  As explained in \cite{N1}, $P^+$ is a pseudodifferental projection.  By Theorem \ref{ThmPSDO}, we know that pseudodifferential operators are bounded on anisotropic Besov spaces.  Thus, $\im P^+$ is a well-defined closed subspace of the boundary data space $B^{(s_1-1/p,s_2)}(E_{\partial X_1 \times X_2})$.  The relevant properties for the operator (\ref{AU1}) can now be understood via the relationship between $\U$ and $\im P^+$.  For this, recall the following definition for pairs of subspaces of a Banach space:

\begin{Def}\label{DefFred}
  A pair of complemented subspaces $(U,V)$ of a Banach space $Z$ is \textit{Fredholm} if $U \cap V$ is finite dimensional and the algebraic sum $U + V$ is closed and has finite codimension.  In this case, we say that $(U,V)$ form a Fredholm pair, or more simply, that $U$ and $V$ are Fredholm (in $Z$).
\end{Def}

Using this notion, we have the following theorem:

\begin{Theorem}\cite{N3} \label{ThmEBP2} Let\footnote{To keep matters simple, we state the hypotheses for $s_1 \geq 1$, which is all we need in this paper.  This is in contrast to \cite{N1}, where we needed to consider spaces of lower regularity than the order of the operator.} $s_1 \geq 1$, $s_2 \geq 0$, and $1<p<\infty$ and consider the operator $\A_\U$ given by (\ref{AU1}).
\begin{enumerate}
\item The operator $A_\U$ is Fredholm if and only if $\U$ and $\im P^+$ are Fredholm.
\item The kernel of $A_\U$ is spanned by finitely many smooth configurations if and only if $\U \cap \im P^+$ is spanned by finitely many smooth configurations.
\item The range of $A_\U$ is complemented by the span of finitely many smooth configurations if and only if $U + \im P^+$ is complemented by the span of finitely many smooth configurations.
\end{enumerate}
\end{Theorem}

Given a subspace $\U \subset B^{(s_1-1/p,s_2)}(E_{\partial X_1 \times X_2})$ Fredholm with $\im P^+$, we can thus construct a projection
\begin{equation}
  \Pi_{\cU}: B^{(s_1-1/p,s_2)}(E_{\partial X_1 \times X_2}) \to B^{(s_1-1/p,s_2)}(E_{\partial X_1 \times X_2})
\end{equation}
with $\im \Pi_{\cU} = \cU$ and $\ker \Pi_\cU$ equal to $\im P^+ = r(\ker A)$ up to some finite dimensional space.  The complementary projection $1-\Pi_{\cU}$ is then a (generalized) elliptic boundary condition, in the sense that $(1-\Pi_\cU): \im P^+ \to \im (1-\Pi_{\cU})$ is Fredholm.  Using this, we can obtain an elliptic estimate for the full mapping pair $(A_\cU, 1-\Pi_{\cU})$.

\begin{Theorem}\cite{N3}\label{ThmEBP3} Let $s_1 \geq 1$, $s_2 \geq 0$, $1<p<\infty$, and suppose $\cU \subseteq B^{(s_1-1/p,s_2)}(E_{\partial X_1 \times X_2})$ is Fredholm with $r(\ker A)$.  Then consider the full mapping pair
\begin{align}
  (A_{\cU}, (1-\Pi_{\cU})r): B^{(s_1,s_2),p}(E) \to B^{(s_1-1,s_2),p}(F) \oplus B^{(s_1-1/p,s_2),p}(E_{\partial X_1\times X_2}).  \label{fullmapA}
\end{align}
This operator is Fredholm and we have the elliptic estimate
\begin{equation}
  \|u\|_{B^{(s_1,s_2),p}(X_1 \times X_2)} \leq C\left(\|Au\|_{B^{(s_1-1,s_2),p}(X_1\times X_2)} + \|(1-\Pi_{\cU})ru\|_{B^{(s_1-1/p,s_2),p}(\partial X_1\times X_2)} + \|\pi u\|\right) \label{EBP-inhom}
\end{equation}
Here $\pi$ is any projection onto the finite dimensional kernel of (\ref{fullmapA}) and $\|\cdot\|$ is any norm on that space.
\end{Theorem}

When $s_2 = 0$, (\ref{EBP-inhom}) is the usual elliptic estimates on isotropic spaces (cf. Theorem I.D.1). Thus, the significance of the above theorem is that tangential anisotropy is preserved.

\section{Banach Space Valued Cauchy Riemann Equations}\label{AppCR}

In this section, we state a modified version of the results of \cite{We}, both to strengthen them for our needs and also to correct some subtle errors. Specifically, we need to make use of the elliptic estimates obtained in \cite{We} for Banach space valued Cauchy Riemann equations with totally real boundary conditions.  Namely, consider the following situation.  We have a Banach space $X$ endowed with a complex structure, i.e., an endomorphism $J: X \to X$ such that $J^2 = -\mathrm{id}$.  A subspace $Y \subset X$ is said to be totally real if $X \cong Y \oplus JY$.  A submanifold $\fL \subset X$ is said to be totally real if each of its tangent spaces is a totally real subspace of $X$.  In particular, for the situation that concerns us, if $X$ is a symplectic Banach space which is densely contained within a strongly symplectic Hilbert space and invariant under the complex structure, then the Lagrangian subspaces and Lagrangian submanifolds of $X$ are all totally real.  For simplicity, we assume we are in this symplectic situation, though everything we do generalizes to the general case.

Given a Lagrangian submanifold $\fL \subset X$ and some $1<p<\infty$, we assume the following hypotheses:

\begin{description}
  \item{$\mathrm{(I)}_p$} There exists a (finite dimensional) vector bundle $E$ over some closed manifold $M$, such that each tangent space to $\fL$ is isomorphic to a closed subspace of the Banach space $L^p(E)$, the space of all $L^p$ sections of $E$.
\end{description}

In the above hypothesis, we assume an inner product on $E$ is chosen so that an $L^p$ norm is defined.  In \cite{We}, the case where $E$ is a trivial bundle is considered, but one can see from the methods there that the more general case can easily be deduced from this latter case.

The next hypothesis, which is omitted from \cite{We}, is one concerning analyticity of the submanifold $\fL$.  First, recall the following definition of analyticity for maps between Banach spaces:

\begin{Def}
  Let $X$ and $Y$ be Banach spaces.  A map $F: X \to Y$ is said to be \textit{analytic} at $0 \in X$ if there exists a neighborhood $U$ of $0 \in X$ and multilinear maps $L_n: X^n \to Y$, $n \geq 0$, such that we have a power series expansion
  \begin{equation}
    F(x) = \sum_{n=0}^\infty L_n(x^n), \label{eq-app:F}
  \end{equation}
  where the series converges absolutely and uniformly for all $x \in U$.  A function $F$ is analytic on an open set $V$ if it is analytic at every point of $V$.
\end{Def}

From this, we can define the notion of an analytic Banach submanifold:

\begin{Definition}\label{DefASub}
  Let $X$ be a Banach space.  An \textit{analytic Banach submanifold} $M$ of $X$ is a subspace of $X$ (as a topological space) that satisfies the following.  There exists a closed Banach subspace $Z \subset X$ such that at every point $u \in M$, there exists an open set $V$ in $X$ containing $u$ and an analytic diffeomorphism $\Phi$ from $V$ onto an open neighborhood of $0$ in $X$ such that $\Phi(V \cap M) = \Phi(V) \cap Z$.  We say that $M$ is modeled on the Banach space $Z$.
\end{Definition}

We have the following additional hypothesis:

\begin{description}
  \item{$\mathrm{(II)}$} The Lagrangian submanifold $\fL \subset X$ is an analytic Banach submanifold of $X$.
\end{description}

Without hypothesis (II), Theorem 1.2 of \cite{We} is incorrect as stated.  We explain what modifications need to be made at the end of this section.

Given the above hypotheses, we consider the following situation.  Let $\Omega \subset \bH$ be a bounded open subset of the half-space
$$\bH = \{(t,v) \in \R^2: v \geq 0\}.$$
For $1 < p < \infty$, if $X$ is any Banach space, let $L^p(\Omega,X)$ denote the Bochner space of all strongly measurable functions $f:\Omega \to X$ such that $\|f\|_{L^p(\Omega,X)} := (\int_\Omega \|f\|_X^p)^{1/p}$ is finite.  We define $L^\infty(\Omega,X)$ similarly.  For $1<p<\infty$ and $k$ a nonnegative integer, let $H^{k,p}(\Omega, X)$ be the Sobolev space of functions\footnote{The notation $W^{k,p}(\Omega, X)$ is more appropriate here, but for all Banach spaces we consider, $H^{k,p}(\Omega,X)$ as traditionally defined (see \cite{N0}) coincides with $W^{k,p}(\Omega, X)$ for all nonnegative integers.} with values in $X$ whose derivatives up to order $k$ belong to $L^p(\Omega, X)$.  We equip $H^{k,p}(\Omega, X)$ with the norm
$$\|u\|_{H^{k,p}(\Omega, X)} = \left(\sum_{j = 0}^k\|\nabla^j u\|_{L^p(\Omega, X)}^p\right)^{1/p}.$$
Observe that if we have a bounded multiplication map
\begin{align}
  H^{k,p}(\Omega) \times H^{k',p'}(\Omega) \to H^{k'',p''}(\Omega)
\end{align}
on the usual scalar valued Sobolev spaces, then this induces a bounded multiplication map
\begin{align}
  H^{k,p}(\Omega, \mathrm{End}(X)) \times H^{k',p'}(\Omega, X) \to H^{k'',p''}(\Omega,X).
\end{align}

Suppose we are given a Lagrangian submanifold $\fL \subset X$ satisfying $\mathrm{(I)}_p$ and (II) for some $1<p<\infty$.  Let $u: \Omega \to X$ be a map that satisfies the boundary value problem
\begin{equation}
\begin{split}
\partial_t u + J\partial_v u &= G\\
u(t,0) & \in \fL, \textrm{ for all } (t,0) \in \partial\Omega \cap \partial \bH,
\end{split}\label{CR-CR}
\end{equation}
where $G: \Omega \to X$ is some inhomogeneous term.  Thus, the system (\ref{CR-CR}) is a Cauchy-Riemann equation for the Banach space valued function $u$, supplemented with a Lagrangian boundary condition.  In \cite{We}, the complex structure $J = J_{t,v}$ is allowed to vary with $t,v \in \Omega$.  For simplicity, and since we will not need to assume otherwise, we let $J$ be constant.

We have the following elliptic regularity theorem for the equations (\ref{CR-CR}), which is a refined and corrected version of \cite[Theorem 1.2]{We}:

\begin{Theorem}\label{ThmCR}
  Fix $1<p<\infty$, let $k \geq 1$, and let $K \subset \mathrm{int}\,\Omega$ be a compact subset.  Let $\fL \subset X$ be a Lagrangian submanifold satisfying $\mathrm{(I)}_p$ and $\mathrm{(II)}$.
  \begin{enumerate}
    \item Suppose $u \in H^{k,q'}(\Omega, X)$ solves (\ref{CR-CR}) with $G \in H^{k,q}(\Omega,X)$, for some $q$ and $q'$ satisfying $p \leq q \leq q'<\infty$.  Furthermore, suppose $q'$, $q$, and $p$ are such that we have bounded multiplication maps
        \begin{align}
          H^{k-1,q'}(\Omega) \times H^{k-1,q'}(\Omega) & \to H^{k-1,p}(\Omega) \label{CR-mult1}\\
          H^{k-1,q'}(\Omega) \times H^{k,q}(\Omega) & \to H^{k-1,p}(\Omega). \label{CR-mult2}
        \end{align}
        Then $u \in H^{k+1,p}(K,X)$.
    \item Furthermore, let $u_0 \in C^\infty(\Omega,X)$ be such that $u_0(t,0)\in\fL$ for all $(t,0) \in \partial\Omega\cap\partial\bH$.  Then there exists a $\delta > 0$ depending on $u_0$ such that if $\|u-u_0\|_{L^\infty(\Omega,X)} < \delta$ is sufficiently small, then
        \begin{equation}
          \|u-u_0\|_{H^{k+1,p}(K,X)} \leq C\bigl(\|G\|_{H^{k,q}(\Omega,X)} + \|u-u_0\|_{H^{k,q'}(\Omega,X)}\bigr). \label{CR-est}
        \end{equation}
        where $C$ is a constant bounded in terms of $\delta$, $u_0$, and $\|u - u_0\|_{H^{k,q'}(\Omega,X)}$.
  \end{enumerate}
\end{Theorem}

Let us note the main differences between Theorem \ref{ThmCR} and \cite[Theorem 1.2]{We}.  First, we allow for a more general range of $q$. In \cite{We}, only the case $q' = q$ is considered, in which case the range permissible of $q$ is such that
\begin{align}
  H^{k-1,q}(\Om) \times H^{k-1,q}(\Om) \to H^{k-1,p}(\Om)
\end{align}
is bounded, which is more restrictive than (\ref{CR-mult1})--(\ref{CR-mult2}). Second, we now explain why we assume the analyticity hypothesis (II).  If this is omitted, then the theorem above must be modified as follows: the constant $\delta$ appearing in Theorem \ref{ThmCR}(ii) a priori depends on $k$. While this seems innocuous, this implies that if one wishes to use (\ref{CR-est}) to bootstrap the regularity of $u$ to higher and higher regularity ($k$ increasing to infinity), one also needs $u$ to get closer and closer to $u_0$.  Analyticity, however, ensures that a small enough $\delta$ works for all $k$.  In a few words, this is because analytic maps, being expressible as a power series locally near any point, satisfy good estimates for all their derivatives within a \textit{fixed} neighborhood of any point. More precisely, for an analytic map, there is a fixed neighborhood about any point on which the $k$th Fr\'echet derivative of the map is Lipschitz, for every $k$. On the other hand, for a general smooth map on an infinite dimensional Banach space, its $k$th Fr\'echet derivative might only be Lipschitz on a neighborhood depending on $k$.\footnote{One of the errors in \cite[Theorem 1.2]{We} is that in the proof, it is asserted that we have an inequality $\|v\|_{W^{k,q}(\Omega)} \leq C\|u-u_0\|_{W^{k,q}(\Omega)}$ for $C$ independent of $u$ if $\delta$ sufficiently small.  This only holds if $\delta$ is sufficiently small depending on $k$.}   Thus, hypothesis (II) ensures that our analytic Banach submanifold has local chart maps obeying Lipschitz estimates on fixed neighborhoods, which gives us the uniformity of $\delta$ with respect to $k$ in Theorem \ref{ThmCR}.

We relegate the proof of Theorem \ref{ThmCR}, which involves minor technical modifications of the work of \cite{We}, to \cite{N0}.  

\end{document}